\documentclass[12pt,a4paper]{article}

\renewcommand{\baselinestretch}{2} 
\RequirePackage{amsthm,amsmath,amsfonts,amssymb}

\usepackage{times}
\usepackage{bm}
\usepackage{natbib}
\usepackage{amsmath}
\usepackage{extarrows}
\usepackage {mathtools}
\usepackage[ruled,vlined]{algorithm2e}
\usepackage{comment}
\usepackage{float}
\usepackage{subfigure}
\usepackage{geometry}
\usepackage{multicol}
\usepackage{color}
\usepackage{caption}
\usepackage{booktabs}

\usepackage{hyperref}
\hypersetup{
	colorlinks=true,
	linkcolor=blue,
	filecolor=magenta,      
	urlcolor=cyan,
	citecolor=black,
	pdfpagemode=FullScreen,
}

\newcommand{\vech}{{\text{vech}}}

\newcommand{\tr}{{\text{tr}}}
\newtheorem{lem}{Lemma}

\newtheorem{thm}{\hspace{1mm}Theorem}
\newtheorem{pro}{\hspace{1mm}Proposition}

\usepackage{authblk}

\title{
        Directional testing for one-way MANOVA in divergent dimensions}
\author[1]{Caizhu Huang\thanks{caizhu.huang@gdufe.edu.cn}}
\author[2]{Claudia Di Caterina\thanks{claudia.dicaterina@univr.it}}
\author[3]{Nicola Sartori\thanks{nicola.sartori@unipd.it}}
\affil[1]{School of Statistics and Data Science,Guangdong University of Finance and Economics}
\affil[2]{Department of Economics, University of Verona}
\affil[3]{Department of Statistical Sciences, University of Padova}

\date{}
\begin{document}
	\maketitle
	
\begin{abstract}
	\noindent
Testing the equality of   mean vectors across $g$ different groups plays an important role in many scientific fields. In regular frameworks, likelihood-based statistics under the normality assumption  offer a general solution to this task. However, the accuracy of standard asymptotic results is not reliable when the dimension $p$ of the data  is large relative to the sample size $n_i$ of each group.  We propose here  an exact directional test for  the equality of $g$  normal mean vectors with identical unknown covariance matrix in a high dimensional setting, provided that $\sum_{i=1}^g n_i \ge p+g+1$. In the case of two groups ($g=2$), the directional test coincides with the Hotelling's $T^2$ test. In the more general situation where the $g$ independent groups may have different unknown covariance matrices, although exactness does not hold, simulation studies show that the directional test is more accurate than  most commonly used likelihood{-}based solutions, at least in a moderate dimensional setting in which $p=O(n_i^\tau)$, $\tau \in (0,1)$.  Robustness of the directional approach and its competitors  under deviation from the assumption of  multivariate normality  is also numerically investigated.
Our proposal is here applied to data on blood characteristics of male athletes and to microarray data storing gene expressions in patients with breast tumors.
\end{abstract}

	
\section{Introduction}\label{section:introduction}
\begin{figure}[t]
	\centering
	\captionsetup{font=footnotesize}
	\subfigure{
		\begin{minipage}[b]{.35\linewidth}
			\centering
			\includegraphics[scale=0.2]{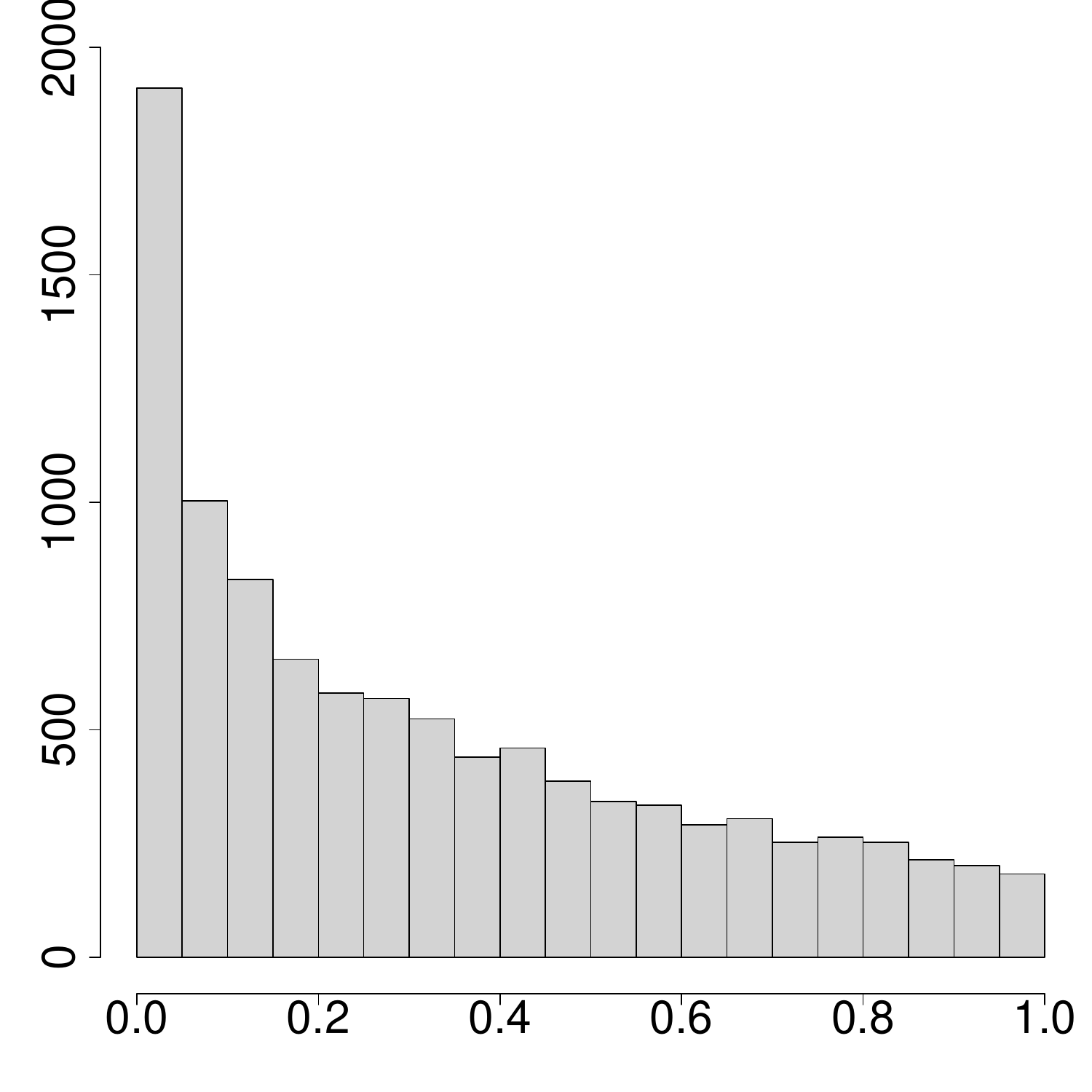}
		\end{minipage}
		\begin{minipage}[b]{.35\linewidth}
			\centering
			\includegraphics[scale=0.2]{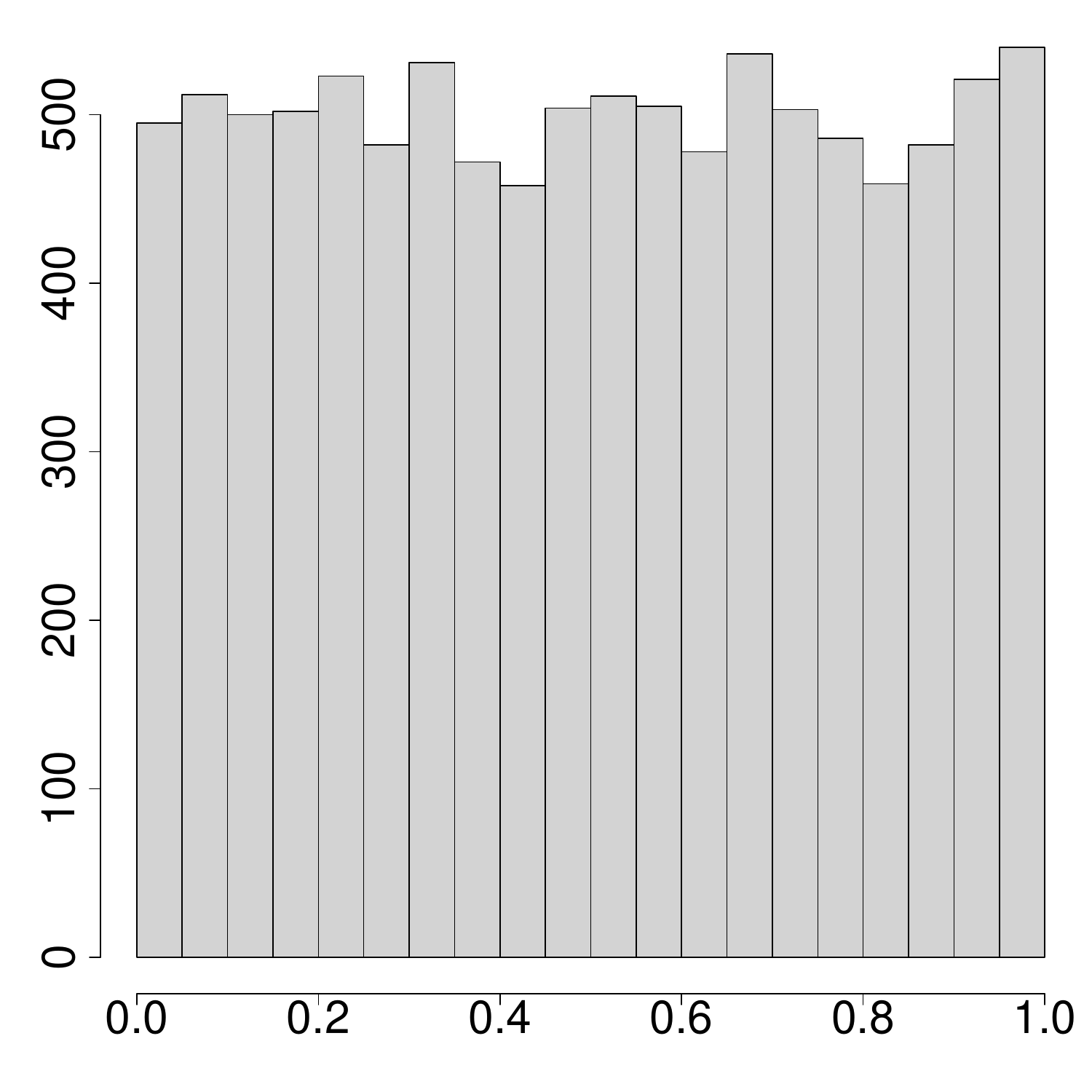}
		\end{minipage}
	}
	\caption{ Empirical null distribution of $p$-values from the likelihood ratio test (left) and directional test (right) for the hypothesis of  equality  of normal  mean vectors in $g=4$ groups with identical covariance matrix, based on $10,000$ Monte Carlo  simulations. Data are generated from a  $N_5(\mu, \Lambda^{-1})$ distribution with mean vector $\mu$ and covariance matrix $\Lambda^{-1}$  equal to the sample mean and sample covariance matrix, respectively, of the \texttt{Pottery} dataset. The total sample size is $n = \sum_{i=1}^{4}n_i=26$. 	  
	}
	\label{fig:empirical p-value example}
\end{figure}

Hypothesis testing for multivariate mean vectors is  a very important inferential problem in many applied research fields. 
Likelihood{-}based statistics and usual asymptotic results offer a general solution to this task in parametric models. Typically, such solutions are accurate when the model dimension $p$ and the sample size $n$ match the standard asymptotic setting, where the dimension of the parameter is considered fixed as the sample size increases. However, usual asymptotic results are no longer guaranteed  when $p$ is not negligible with respect to $n$ \citep[see for instance][]{jiang:2013,sur2019likelihood,tang2020modified,He2020B}. As a simple illustration, the left panel of Figure \ref{fig:empirical p-value example} shows the empirical null distribution based on 10,000 Monte Carlo samples  of $p$-values obtained using the asymptotic chi-square approximation for the  likelihood ratio statistic when testing the equality of $g=4$ normal mean vectors under the assumption of common unknown covariance matrix. This problem is  known as  homoscedastic one-way multivariate analysis of variance (MANOVA), and the simulation setup  is taken from the \texttt{Pottery} data in the R \citep{Rcode} package \texttt{car} \citep{car}: the group sizes  are  $5, 2, 5$ and $14$, respectively, with $n = \sum_{i=1}^4n_i =26$, and the dimension of the vectors is $p=5$. It is clear from  Figure \ref{fig:empirical p-value example} that the standard approximate $p$-values obtained from the likelihood ratio statistic  are far from being uniform, as opposed to the directional $p$-values proposed in this paper and  shown in the right panel. 

Settings like the former, in which the values of $p$ and $n$ are comparable, may be framed in a $(p, n)$-divergent dimensional asymptotic setting where both $p$ and $n$ are allowed to increase. Indeed, the data dimension $p$ is related to the number of parameters, e.g., for the homoscedastic one-way MANOVA case in Section \ref{section:hypothesis same}, the numbers of parameters is $gp+p(p+1)/2$.  Inspired by  \citet{battey2022some}, we distinguish between three asymptotic regimes: moderate dimensional, high dimensional and ultra-high dimensional. Here we do not deal with the ultra-high dimensional asymptotic regime, in which $p/n$ diverges or tends to a limit greater than one. Instead,  we focus on the  moderate dimensional asymptotic regime, in which $p/n \to 0$, in particular with $p = O(n^\tau)$, $\tau \in (0,1)$, and on the high dimensional asymptotic regime, in which  $p/n \to \kappa \in (0,1)$.  
In the moderate dimensional asymptotic regime,  \citet{He2020B} proved that in testing the equality of  normal mean vectors with identical covariance matrix the likelihood ratio test is valid if $p/n^{2/3} \to 0$, while the analogous condition for its Bartlett corrected version is $p/n^{4/5} \to 0$. To our knowledge, no similar results are available for the heteroscedastic case. 

Higher-order likelihood solutions based on saddlepoint approximations might generally  give substantial improvements over first{-}order solution, especially in high dimensions \citep[see, e.g.,][]{tang2020modified}. Among these, directional inference on a vector parameter of interest, as developed by \citet{sartori:2014} and \citet{sartori:2016}, has proven to be particularly accurate when testing canonical parameters in exponential families.   Its accuracy descends from  that of   the underlying saddlepoint approximation  to the conditional density of the canonical statistic of interest. 
Empirical results in \citet{sartori:2014} and \citet{sartori:2016}  showed that  directional tests are extremely accurate even in settings where the dimension of the parameter of interest, although lower than the sample size, has a comparable order.
The use of  a saddlepoint approximation, indeed, guarantees a fairly constant relative error.

The use of a directional test may be motivated by the fact that, 
in standard asymptotics,
the directional test is first-order equivalent to the likelihood ratio test. Yet, if  
first-order approximations are needed for the distribution of the latter, the directional test may be more convenient, given its improved accuracy  \citep{skogaard:1988,sartori2017introduction}.
Moreover, \citet{sartori:2019Ftest} showed that the directional test coincides with many well-known exact tests. For {example}, when testing a specific value for the mean of a multivariate normal distribution, the  directional test coincides with the exact  Hotelling's $T^2$ test.  \citet{huang2021directional} found other instances in which  the directional $p$-value is exactly uniformly distributed. Such examples are related to  several prominent inferential problems in which the independence of components or the equality of covariance matrices is tested in high dimensional multivariate normal models.  
Finally, \citet{claudia2023ss} showed  the accurate  properties of directional tests for covariance selection in high dimensional Gaussian graphical models.

Concentrating on tests for  the hypothesis of  equality of mean vectors in $g$ independent groups, this work makes a number of  contributions to the current literature. First, under the assumptions of normality and   identical unknown covariance matrix, we {prove} that the directional $p$-value is exactly uniformly distributed provided that  $n=\sum_{i=1}^g n_i \ge p+g+1$, and  coincides with the Hotelling's $T^2$ test {when} $g=2$.  
For the more general case with  $g$  unknown group  covariance matrices, known as the  Behrens-Fisher problem {if} $g=2$, the directional test is not exact. 
Still,  we show by means of extensive  simulation studies  that the directional approach overperforms  standard first-order solutions  as well as other higher-order modifications \citep{skovgaard:2001} in moderate dimensional settings.
In addition, the robustness of the available  solutions to the normality assumption {is empirically investigated}, considering  multivariate $t$, skew-normal or Laplace true generating processes in the Supplementary Material. Finally, the proposed directional test is applied for the analysis of two datasets: one from a study on the blood characteristics of male athletes, and another, presented in the Supplementary Material, from a microarray analysis on breast cancer patients.

The rest of the paper is organized as follows.  Section \ref{section:background} presents some background information.  In particular, Section \ref{section:notation} reviews some  likelihood-based statistics in exponential families, Section \ref{section:directional testing} {reports} the steps to compute the directional $p$-value,  and Section \ref{section:multivariate normal} details the necessary quantities for the multivariate normal model. The main results in Section \ref{section:main result} are for hypotheses concerning: (i) the equality  of  $g$ normal mean vectors with identical covariance matrix (Section \ref{section:hypothesis same}); (ii) the equality of $g$ normal mean vectors with different covariance matrices (Section \ref{section:hypothesis diff}).  For  hypothesis (i),  we prove the exact uniform distribution of the directional $p
$-value under the null.  In Section \ref{section:simulation studies} we report empirical results under the assumed model. An application of the considered methodology to real data, specifically blood characteristics of male athletes, is illustrated in Section \ref{sec:app}, while Section \ref{section:conclusion} concludes with a discussion.     
The Supplementary Material contains auxiliary computational results, proofs of the theorems, additional simulation studies, including larger group sizes and robustness analyses, as well as the second application involving breast cancer patient  microarray data.

\section{Background}\label{section:background}

\subsection{Notation and setup}\label{section:notation}
Suppose that data $y = (y_1, \dots, y_n)^\top$ are  generated from the model $f(y;\theta)$ with parameter $\theta = (\theta_1, \dots, \theta_q)^\top$. The log-likelihood function is $\ell(\theta) = \ell(\theta;y)= \log f(y; \theta)$, 
and we are interested in testing the null hypothesis $H_\psi: \psi(\theta) = \psi$ on the $d$-dimensional parameter of interest $\psi$. It will often be the case that $\psi$ is a component of $\theta$, possibly after a reparameterization. Assume the partition  $\theta = (\psi^\top, \lambda^\top)^\top$ holds, with $\lambda$ a  $(q-d)$-dimensional nuisance parameter. Let $\hat{\theta}$ denote the maximum likelihood estimate of $\theta$ and $\hat{\theta}_\psi$ its constrained maximum likelihood estimate  under $H_\psi$, i.e., $\hat{\theta}_\psi = (\psi^\top, \hat{\lambda}_\psi^\top)^\top$. 

Several likelihood-based statistics can be used to test the hypothesis $H_\psi$.  The likelihood ratio statistic, a parameterization-invariant measure, is
\begin{align}
	W = 2 \{\ell(\hat{\theta}) - \ell(\hat{\theta}_\psi)\}.
	\label{W}
\end{align}
When $q$ is fixed and $n \to \infty$, the statistic $W$ has a $\chi^2_d$ asymptotic  null distribution with relative error of order $O(n^{-1})$. A correction of $W$ proposed by \citet{bartlett:1937} rescales the likelihood ratio statistic by its expectation under $H_\psi$, that is
\begin{align}
	W_{BC} = \frac{d}{E(W)} W, 
	\label{WBC}
\end{align}
and  follows asymptotically a $\chi^2_d$ null distribution with relative error of order  $O(n^{-2})$ \citep[Section 7.4]{mccullage:2018}. More details on  the expectation $E(W)$ can be found in  \citet[][Section 2.1]{huang2021directional}.

\citet{skovgaard:2001} introduced two improvements on $W$, namely
\begin{align}
	W^* = W \left(1-\frac{\log \gamma}{W}\right)^2\quad \text{and} \quad W^{**} = W - 2 \log \gamma,
	\label{Wstar}
\end{align}
which also have approximate  $\chi^2_d$ distributions when the null hypothesis holds. The correction factor $\gamma$ can be found in \citet[][Equation (13)]{skovgaard:2001}, and is also reported in  S1.1  with the notation used here.

The test  statistics presented so far are omnibus measures of departure of the data from $H_\psi$: their $p$-value results from averaging  the deviations from $H_\psi$ in all the potential directions {of the parameter space}.   We {will} now introduce  the directional test developed by \citet{sartori:2014} and \citet{sartori:2016}, which measures the departure from $H_\psi$ along  the direction indicated by the observed data.

\subsection{Directional testing}\label{section:directional testing}

Let $\varphi = \varphi(\theta)$ be a reparameterization of the original model. Suppose we have an exponential family model with sufficient statistic $u = u(y)$ and canonical parameter $\varphi$, with density 
\begin{align}
	f(y;\theta) = \exp\left[\varphi(\theta)^\top u -K\{\varphi(\theta)\}\right]h(y),
	\label{exponential family}
\end{align}
maximum likelihood estimate $\hat{\varphi} = \varphi(\hat{\theta})$ and constrained maximum likelihood estimate $\hat{\varphi}_\psi = \varphi(\hat{\theta}_\psi)$. Henceforth, we shall use the 0 superscript to indicate quantities evaluated at the observed data point $y^0$. For computing the directional $p$-value, it is convenient to define a centered sufficient statistic at  $y^0$, $s = u -u^0$, with $u^0 = u(y^0)$. 
The tilted log-likelihood function of  model (\ref{exponential family}) takes the form
\begin{align}
	\ell(\varphi;s) = \varphi(\theta)^\top s + \ell^0(\theta),\nonumber
\end{align}
where  $\ell^0(\theta) = \ell(\theta;y^0)$ is the  observed log-likelihood function. 
The saddlepoint approximation \citep[see, e.g.,][Section 10]{salvan:1997} to the exponential model on $\mathbb{R}^q$ is 
\begin{align}
	f(s;\varphi) = c \exp\{\ell(\varphi;s) - \ell(\hat{\varphi};s)\} |J_{\varphi\varphi}(\hat{\varphi})|^{-1/2},\nonumber
\end{align}
where $J_{\varphi\varphi}(\hat{\varphi})$ is the observed Fisher information  $J_{\varphi\varphi}(\varphi)= - \partial^2 \ell(\varphi;s) / (\partial \varphi \varphi^\top)$ evaluated at $\hat \varphi$ and  $c$ is a normalizing constant.

The hypothesis $H_\psi$ specifies a value for the parameter $\psi = \psi(\varphi)$. Following \citet{sartori:2016}, the reduced model in $\mathbb{R}^d$ is given by
\begin{align}
	h(s;\psi) = c \exp\{\ell(\hat{\varphi}_\psi;s) - \ell(\hat{\varphi};s)\}|J_{\varphi\varphi}(\hat{\varphi})|^{-1/2} |\tilde{J}_{(\lambda \lambda)}|^{1/2}, \quad  s \in L_\psi^0, 
	\label{saddeL0}
\end{align}
where $L_\psi^0$ is the $d$-dimensional plane  obtained by setting $\hat \lambda_\psi = \hat{\lambda}_\psi^0$,  and the observed information for the nuisance parameter has been recalibrated to $\varphi$ as follows:
\begin{align}
	|\tilde{J}_{(\lambda \lambda)}| = |J_{\lambda \lambda}(\hat{\varphi}_\psi;s)| |\partial \varphi(\theta)/\partial \lambda|^{-2}_{\hat{\theta}_\psi}.
	\label{J for lambda}
\end{align}


The directional test is constructed by defining a line $L_\psi^*$ through the observed value $s^0=0_q$ and the expected value of $s$ under $H_\psi$ which depends on the  observed data point $y^0$, i.e. $s_{\psi} = - \ell_{\varphi}^0 \left(\hat{\varphi}^0_\psi \right) = - \partial \ell^0\left\{\theta(\varphi)\right\}/\partial\varphi\left|_{\varphi=\hat{\varphi}^0_\psi}\right.$. We parameterize this line by $t \in  \mathbb{R}$, namely $s(t) = s_{\psi} + t (s^0 - s_{\psi})$. In particular, $s(0) = s_\psi$, corresponding to $H_\psi$, and $s(1) = s^0$, corresponding to the observed data. Then, the directional test  for $s$ is conditional on being on the line $L^*_\psi$, i.e.  conditional on $s(t)/||s(t)||$.
The directional  $p$-value measuring the departure from $H_{\psi}$ along the line $L_{\psi}^*$  is defined as the probability that $s(t)$ is as far or farther from $s_\psi$ than is the observed value $s^0$:
\begin{eqnarray}
	p(\psi) = \frac{\int_{1}^{t_{\text{sup}}} {t^{d-1}h\{s(t);\psi\}} \text{d}t}{\int_{0}^{t_{\text{sup}}} {t^{d-1}h\{s(t);\psi\}}\text{d}t},
	\label{directed p-value}
\end{eqnarray}
where the denominator is a normalizing constant. The upper limit $t_{\text{sup}}$ of the integrals in (\ref{directed p-value}) is the largest value of $t$ for which the maximum likelihood estimate $\hat{\varphi}(t)$, corresponding to $s(t)$, exists. See \citet[][Section 3]{sartori:2016} for more details. 

In the particular case where $\psi$ and $\lambda$ are  linear functions of the canonical parameter of an exponential family, the quantity (\ref{J for lambda}) does not depend on $s$ and can  therefore be ignored in computing the directional $p$-value (\ref{directed p-value}). Moreover, the expected value $s_\psi$ simplifies to $s_\psi =\left\{-  \ell_{\psi}^0 \left(\hat{\varphi}^0_\psi \right)^\top, {{0}_{q-d}^\top}\right\}^\top$.

If the original model $f(y; \theta)$ {has} not the exponential family form, then a tangent exponential model  is used instead of (\ref{exponential family}). The construction of the tangent exponential model and its saddlepoint approximation are described in  \citet[][Appendix]{sartori:2016}; see also \cite{davison2022tangent}. {Here this step is not needed since we work with} the normal model which {belongs to the} exponential family, yet {the} saddlepoint approximation (\ref{saddeL0}) will be {used} when $\psi$ is not a linear function of the canonical parameter, as in the heteroscedastic one-way MANOVA of Section \ref{section:hypothesis diff}.

\subsection{Independent groups from  multivariate normal distributions}\label{section:multivariate normal}

Consider $g$ independent groups, and denote by  $y_{ij}$ the independent observations from the $i$th group with  multivariate normal  distribution $N_p(\mu_i, \Lambda_i^{-1})$, $i\in \{1,\dots, g\}; j\in \{ 1,\dots,n_i\}$. The mean vectors $\mu_i$ and the concentration matrices $\Lambda_i$, symmetric and positive definite, are assumed unknown. Let $\text{tr}(A)$ indicate the trace of a square matrix $A$, and $\text{vec}(A)$ be the  vector stacking the columns of $A$ one by one. We also define the vector  $\text{vech}(A)$, obtained from $\text{vec}(A)$ by eliminating all upper triangular elements of $A$ when this is symmetric. These two vectors satisfy the relationship  $D_p \text{vech}(A) = \text{vec}(A)$, where $D_p$ is the so-called duplication matrix \citep[Section 3.8]{magnus1999}. 

We rewrite the data from the $i$th group as  $y_i = [y_{i1} \cdots y_{in_i}]^\top$, which is a $n_i \times p$ matrix. Then, the log-likelihood for the parameter $\theta = \{\mu_1^\top, \dots, \mu_g^\top, \text{vech}(\Lambda_1^{-1})^\top, \dots, \text{vech}({\Lambda}_g^{-1})^\top\}^\top$ is
\begin{align}
	\ell(\theta)
	&= \sum_{i=1}^{g} n_i\mu_i^\top \Lambda_i  \bar{y}_i  -\frac{1}{2} \tr (\Lambda_i y_i^\top y_i) + \frac{n_i}{2} \log|\Lambda_i| - \frac{n_i}{2} \mu_i^\top \Lambda_i \mu_i,\nonumber
\end{align}
where $\bar{y}_i = y_i^\top 1_{n_i}  / n_i$ with $1_{n_i}$ a $n_i$-dimensional vector of ones, $i\in \{1,\dots,g\}$. 
In this exponential family model, the canonical parameter is   $\varphi = \{\xi_1^\top,\dots, \xi_g^\top, $ $  \vech(\Lambda_1)^\top, \dots, \vech(\Lambda_g)^\top\}^\top = \{\mu_1^\top\Lambda_1,\dots, \mu_g^\top\Lambda_g, $ $ \vech (\Lambda_1)^\top,\dots, \vech (\Lambda_g)^\top\}^\top$. 
The log-likelihood as a function of $\varphi$ is $\ell(\varphi) = \sum_{i=1}^{g} \ell_i(\varphi_i)$ with $\varphi_i = \{\xi_i^\top, \vech(\Lambda_i)^\top\}^\top$, and the $i$-th group's contribution is 
\begin{align}
	\ell_i(\varphi_i)
	&=  n_i\xi_i^\top \bar{y}_i  -\frac{1}{2} \tr (\Lambda_i y_i^\top y_i) + \frac{n_i}{2} \log|\Lambda_i| - \frac{n_i}{2} \xi_i^\top \Lambda_i^{-1} \xi_i\nonumber\\
	& = \xi_i^\top n_i\bar{y}_i  - \text{vech}(\Lambda_i)^\top \left\{\frac{1}{2} D_p^\top D_p \text{vech}(y_i^\top y_i) \right\} + \frac{n_i}{2} \log|\Lambda_i| - \frac{n_i}{2} \xi_i^\top \Lambda_i^{-1} \xi_i.\nonumber
\end{align}
The score function is $\ell_{\varphi}(\varphi) =\left.{\partial} \ell(\varphi) \right/ {\partial \varphi} = \{\ell_{\varphi_1}(\varphi_1),\dots, \ell_{\varphi_g}(\varphi_g)\}^\top$, with 
\begin{align}
	\ell_{\varphi_i} (\varphi_i) &=\left\{\ell_{\xi_i}(\varphi_i)^\top, \ell_{\text{vech}(\Lambda_i)}(\varphi_i)^\top\right\}^\top\nonumber\\
	&= \left\{n_i \bar{y}_i^\top - n_i \xi_i^\top \Lambda_i, \;\frac{n_i}{2} \text{vech}(\Lambda_i^{-1} - y_i^\top y_i/n_i + \Lambda_i^{-1} \xi_i \xi_i^\top \Lambda_i^{-1})^\top\right\}^\top.\nonumber
\end{align}
The maximum likelihood estimates for $\mu_i$ and $\Lambda_i^{-1}$ are  $\hat{\mu}_i = \bar{y}_i$ and  $\hat{\Lambda}_i^{-1}=y_i^\top y_i /n_i -\bar{y}_i\bar{y}_i^\top$,  respectively; hence, $\hat{\xi}_i = \hat{\Lambda}_i \hat{\mu}_i$, $i\in \{1,\dots,g\}$. Moreover, the observed information matrix $J_{\varphi \varphi}(\varphi)$  can be written in a block-diagonal form, with the group-specific diagonal block $J_{\varphi_i \varphi_i}(\varphi_i)$, $i\in \{1,\dots,g\}$, equal to \citep[Section 2.3]{huang2021directional}
\begin{eqnarray}
	J_{\varphi_i \varphi_i} (\varphi_i)=  \left[\begin{array}{cc}
		n_i\Lambda_i^{-1}   & \quad-n_i(\xi_i^T \Lambda_i^{-1} \otimes \Lambda_i^{-1} ) D_p \\
		-n_i D_p^T (\Lambda_i^{-1} \xi_i \otimes \Lambda_i^{-1})    & \quad \frac{n}{2} D_p^T \{ \Lambda_i^{-1} (I_p + 2\xi_i \xi_i^T \Lambda_i^{-1}) \otimes \Lambda_i^{-1}\} D_p
	\end{array} \right],\nonumber
\end{eqnarray}
where $ \otimes $ denotes the Kronecker product  \citep[see, e.g.,][Section 5.1]{Lauritzen:1996}. 
Then, the  determinant of  $J_{\varphi \varphi}$ is such that $|J_{\varphi \varphi} (\varphi)|=\prod_{i=1}^{g}|	J_{\varphi_i \varphi_i} (\varphi_i)| \propto \prod_{i=1}^{g} |{\Lambda}_i^{-1}|^{p+2}$.

If the covariance matrices of the $g$ groups are the same, $\Lambda^{-1}_i = \Lambda^{-1}$ for $i \in \{1,\dots,g\}$, the canonical parameter is $\varphi = \{\xi_1,\dots, \xi_g,$ $\text{vech}(\Lambda)^\top \}^\top = \left\{\mu_1^\top\Lambda,\dots,\mu_g^\top\Lambda,\text{vech}(\Lambda)^\top\right\}^\top$, and the maximum likelihood estimate for $\mu_i$ and $\Lambda^{-1}$ are, respectively, $\hat{\mu}_i=\bar{y}_i$ and $\hat{\Lambda}^{-1} = B/n$ with $B = \sum_{i=1}^g y_i^\top y_i - n_i \bar{y}_i \bar{y}_i^\top$ and $n= \sum_{i=1}^{g}n_i$. 
In this setting, the observed information matrix $J_{\varphi\varphi} (\varphi)$ can be computed in block form (see Section S1.2).

\section{One-way MANOVA problems}\label{section:main result}

\subsection{Homoscedastic one-way MANOVA}\label{section:hypothesis same}
\noindent
Suppose that $y_{ij}$ are independent observations from distributions $N_p(\mu_i, \Lambda^{-1})$, for $i\in \{1,\dots, g\}\; (g \ge 2)$ and $j \in \{ 1,\dots,n_i\}$. We are interested in testing the equality of the $g$ mean vectors:
\begin{align}
	H_\psi: \mu_1 = \dots = \mu_g.
	\label{hypothesis:dirmean}
\end{align}
The hypothesis problem (\ref{hypothesis:dirmean}) is equivalent to testing
\begin{align}
	H_\psi: \Lambda \mu_1 = \dots = \Lambda \mu_g.\nonumber
\end{align}
In this framework, the $q$-dimensional canonical parameter is   $\varphi = \{\xi_1^\top,\dots, \xi_g^\top,  \vech(\Lambda^\top)\}^\top =\{\mu_1^\top\Lambda,\dots,  \mu_g^\top\Lambda, $ $  \vech (\Lambda^\top)\}^\top$, with $d$-dimensional parameter of interest $\psi = (\xi_2^\top - \xi_1^\top, \dots, \xi_g^\top - \xi_1^\top)^\top$ and $(q-d)$-dimensional  nuisance parameter  $\lambda =\{\xi_1^\top, $ $ \vech(\Lambda^\top)\}^\top$. The parameter of interest $\psi$ is therefore a linear function of the canonical parameter $\varphi$. The maximum likelihood estimates of parameters $\mu_i$ and $\Lambda^{-1}$ are  given in  Section \ref{section:multivariate normal}. The constrained maximum likelihood estimate are $\hat{\mu}_0 = \bar{y}$ and $\hat{\Lambda}^{-1}_0 = \left(\sum_{i=1}^{g} y_i^\top y_i - n \bar{y}\bar{y}^\top\right)/n = (A+B)/n$ with $\bar{y} = \sum_{i=1}^{g} n_i\bar{y}_i/n$, $A= \sum_{i=1}^{g} n_i \bar{y}_i^\top \bar{y}_i - n_i \bar{y}\bar{y}^\top$, and $B$ {as defined} in Section \ref{section:multivariate normal}. 

There are several likelihood-based tests for the hypothesis problem (\ref{hypothesis:dirmean}) 
when the dimension $p$ is fixed and the group size $n_i$ goes to infinity. 
However, as discussed in Section \ref{section:introduction}, our focus here is on improving the existing standard approximations when $p$ is relatively large with respect to $n_i$.
Below we provide  the necessary and sufficient conditions for the validity of the directional test in the high dimensional regime where $p/ n_i \to \kappa \in (0,1), i\in \{1,\dots,g\}$. 

First, we summarize here  the key methodological steps to compute the directional $p$-value (\ref{directed p-value}) for testing hypothesis (\ref{hypothesis:dirmean}) \citep[see also][for more details]{sartori:2014}.
Under $H_\psi$, we define the expected value of $s$ under $H_\psi$ as   $s_\psi = (n_1\bar{y}^\top - n_1 \bar{y}_1^\top, \dots, $ $ n_g\bar{y}^\top-n_g\bar{y}_g^\top, 0_{q-d})^\top$, and the line $s(t) = (1-t) s_\psi$.
The tilted log-likelihood along the line $s(t)$ is  then
\begin{align}
	\ell (\varphi;t) = \sum_{i=1}^g \ell_i(\varphi_i;t) =\sum_{i=1}^g \ell_i(\varphi_i) + \varphi^\top_i s_i(t),\nonumber 
\end{align}
where 
\begin{align}
	\ell_i(\varphi_i;t)
	=& {n_i}  \xi_i^\top \left\{  (1-t)\bar{y} + t\bar{y}_i \right\} - \frac{1}{2} \tr\left(\Lambda y_i^\top y_i\right) +  \frac{n_i}{2} \log|\Lambda| -\frac{n_i}{2}\xi_i^\top\Lambda^{-1}\xi_i.   \nonumber
\end{align}
The corresponding saddlepoint approximation is 
\begin{align}
	h\{s(t);\psi\}
	=& c \exp\left\{\frac{(n-p-g-1)}{2} \log|\hat{\Lambda}(t)|\right\},\nonumber
\end{align}
where $c$ is a normalizing constant. 
The following lemma states that $t_{\text{sup}}$ in  (\ref{directed p-value}), i.e. the largest $t$ such that $\hat{\Lambda}(t)^{-1}$  is positive definite, is equal to  $1/\sqrt{{\nu_{(p)}}}$,  where $\nu_{(p)}$ is the largest eigenvalue of the matrix $(B_0^\top)^{-1} (A/n) B_0^{-1}$ with $B_0$ the square root  of $\hat{\Lambda}_0^{-1}$, namely $\hat{\Lambda}^{-1}_0 = B_0^\top B_0$.  The function \texttt{chol()} in the \texttt{R}  package \texttt{Matrix} \citep{Matrix} can be used to compute $B_0$. 

\begin{lem}\label{lemma1}
	The estimator $\hat{\Lambda}^{-1}(t)$ is positive definite if and only if $t \in \left[0, 1/\sqrt{\nu_{(p)}}\right]$, where $\nu_{(p)}$ is the largest eigenvalue of $(B_0^\top)^{-1} (A/n) B_0^{-1}$, with $\hat{\Lambda}^{-1}_0 =  B_0^\top B_0$.
\end{lem}

Thanks to the favorable properties of the multivariate normal distribution, we can show that the  saddlepoint approximation to the conditional density of $s$ along the line $s(t)$ is exact, and so {when the null hypothesis holds} the directional $p$-value is exactly uniformly distributed even in the high dimensional asymptotic  regime. The condition for the validity of this result is given in the following theorem. 

\begin{thm}\label{theorem1}
	Assume that $p=p_n$ is such that $n \ge p + g +1$, with $n = \sum_{i=1}^g n_i$ and fixed $g$. Then, under the null hypothesis (\ref{hypothesis:dirmean}), the directional $p$-value is exactly uniformly distributed.
\end{thm}

When the number of independent groups is $g=2$, the  Hotelling's $T^2$ test statistic with exact $F$ distribution can be used for hypothesis  (\ref{hypothesis:dirmean}).  We also prove that in this case the directional $p$-value coincides with  the one {from the} Hotelling's $T^2$ test.

\begin{pro}\label{proposition1}
	When $g=2$, the directional test is equivalent to the Hotelling's $T^2$ test.
\end{pro}
\noindent
The proofs of Lemma \ref{lemma1}, Theorem \ref{theorem1} and Proposition \ref{proposition1}  are given in  Section S2 of the Supplementary Material.

For comparison, we derive the expressions {in this framework} of the  likelihood ratio test $W$, its Bartlett corrected version $W_{BC}$, and the two modifications $W^*$ and $W^{**}$ proposed by \citet{skovgaard:2001}.  The likelihood ratio test statistic (\ref{W}) is 
\begin{align}
	W = n (\log |\hat{\Lambda}^{-1}_0| - \log |\hat{\Lambda}^{-1}|).\nonumber
\end{align}
Under $H_\psi$, its distribution is approximated by a $\chi^2_d$ with degrees of freedom  $d=p(g-1)$  if and only if $p = o(n^{2/3})$ \citep{He2020B}. 
In this framework, the expectation $E(W)$ in the Bartlett corrected test (\ref{WBC})  can be calculated exactly as done by \citet[][Section A.1]{He2020B}. The $\chi^2_d$ approximation for the distribution of $W_{BC}$ holds if and only if  $p=o(n^{4/5})$ \citep{He2020B}.
The statistics  $W^{*}$ and $W^{**}$ in (\ref{Wstar})  for hypothesis (\ref{hypothesis:dirmean}) can be computed  explicitly, based on the formula (S1) in the Supplementary Material for the  correction factor $\gamma$.
The quantities required  are 
$(\hat{\varphi} - \hat{\varphi}_\psi)^\top (s-s_\psi) = \tr(\hat{\Lambda}A)$, $\log(|{J}_{\varphi\varphi}(\hat{\varphi}_\psi)|/|{J}_{\varphi\varphi}(\hat{\varphi})|)=(p+g+1)(\log|\hat{\Lambda}| - \log|\hat{\Lambda}_0|)$ and 
\begin{align}
	(s-s_\psi)^\top J_{\varphi\varphi}(\hat{\varphi}_{\psi})^{-1}(s-s_\psi) =& \left\{\sum_{i=1}^g n_i(\bar{y}_i - \bar{y})^\top \hat{\Lambda}_0 (\bar{y}_i -\bar{y})\right\} \nonumber\\
	&+ \left\{\sum_{i=1}^g n_i(\bar{y}_i-\bar{y})^\top\right\} \frac{(\bar{y}^\top\hat{\Lambda}_0\bar{y})\hat{\Lambda}_0+\hat{\Lambda}_0\bar{y}\bar{y}^\top\hat{\Lambda}_0}{n}\left\{\sum_{i=1}^g n_i(\bar{y}_i-\bar{y})\right\}.\nonumber
\end{align}




\subsection{Heteroscedastic one-way MANOVA}\label{section:hypothesis diff}

Suppose that  $y_{ij}$ are independent observations from distributions $N_p(\mu_i, \Lambda_i^{-1})$ for $i\in \{1,\dots,g\}$ and $j \in  \{1,\dots,n_i\}$.   We are again interested in testing the equality of the $g$ mean mean vectors $\mu_i$, i.e. 
\begin{eqnarray}
	H_\psi: \mu_1 = \dots = \mu_g.
	\label{hypothesis:meansdiff}
\end{eqnarray} 
In this framework, the constrained maximum likelihood estimates  are denoted by $\hat{\mu}_{0i}=\hat{\mu}_0$ and $\hat{\Lambda}_{0i}^{-1}$, where $\hat{\Lambda}_{0i}^{-1} = \hat{\Lambda}_i^{-1} +  (\bar{y}_i - \hat{\mu}_0) (\bar{y}_i - \hat{\mu}_0)^\top$. We compute the constrained maximum likelihood estimate $\hat{\mu}_0$ numerically by maximization of the profile log-likelihood $\ell_P(\mu) = -\sum_{i=1}^g (n_i/2)\log|\hat{\Lambda}_i^{-1} +$ $  (\bar{y}_i - \mu) (\bar{y}_i - \mu)^\top |. $

To develop the directional test under the null hypothesis (\ref{hypothesis:meansdiff}), following Fraser et al. (2016) we consider the parameterization $(\psi, \lambda)$ with parameter of interest $\psi = (\mu_2^\top-\mu_1^\top, \dots,\mu_g^\top-\mu_1^\top)^\top$ and  nuisance parameter $\lambda = \{\mu_1^\top, {\text{vech}}(\Lambda_1)^\top, $ $\ldots, {\text{vech}}(\Lambda_g)^\top\}^\top$. This parameterization places nonlinear constraints on the canonical parameter $\varphi$. 
The tilted log-likelihood is $\ell (\varphi;t) = \sum_{i=1}^g \ell_i(\varphi_i;t)$ 
with $\varphi_i = \{\mu_i^\top\Lambda_i, \text{vech}(\Lambda_i)^\top\}^\top$ and the $i$-th group's contribution is
\begin{align}
	\ell_i(\varphi_i;t)
	=& \; {n_i} \mu_i^\top\Lambda_i \left\{ t\bar{y}_i + (1-t)\hat{\mu}_0\right\} -\frac{1}{2}{\text{tr}\;} \left( \Lambda_i \left[  y_i^\top y_i + (1-t)\left\{n_i\hat{\mu}_0 (\bar{y}_i-\hat{\mu}_0)^\top  + n_i(\bar{y}_i - \hat{\mu}_0)\hat{\mu}_0^\top\right\} \right] \right) \nonumber \\
	& \; +\frac{n_i}{2} \log|\Lambda_i| -\frac{n_i}{2}\mu_i^\top\Lambda_i\mu_i .\nonumber
\end{align}
The expected value $s_\psi$ of the corresponding sufficient statistic $s$ under $H_\psi$ has components $[n_i \hat{\mu}_0^\top - n_i \bar{y}_i^\top, $  $ \frac{n_i}{2} {\text{vech}}\{ \hat{\mu}_0 (\bar{y}_i-\hat{\mu}_0)^\top + (\bar{y}_i - \hat{\mu}_0)\hat{\mu}_0^\top\}^\top]$, $i\in \{1,\dots,g\}$. The maximum likelihood estimates along the line $s(t) = (1-t)s_\psi$ are $\hat{\mu}_i(t) = t\bar{y}_i + (1-t)\hat{\mu}_0$ and $\hat{\Lambda}_i^{-1}(t) = \hat{\Lambda}_{0i}^{-1} - t^2(\bar{y}_i-\hat{\mu}_0)(\bar{y}_i-\hat{\mu}_0)^\top$. The maximum likelihood estimate $\hat{\Lambda}_i^{-1}(t)$ exists for  $t \le t_{\text{sup}} = \mathop{\min}\limits_{1\le i \le g} \left[\{(\bar{y}_i-\hat{\mu}_0)^\top \hat{\Lambda}_{0i}(\bar{y}_i-\hat{\mu}_0) \}^{-1/2} \right]$, with $(\bar{y}_i-\hat{\mu}_0)^\top \hat{\Lambda}_{0i}(\bar{y}_i-\hat{\mu}_0) \ne 0$. Therefore, the saddlepoint approximation along the line $s(t)$ is 
\begin{align}
	h\{s(t);\psi\}
	=&  \prod_{i=1}^g |\hat{\Lambda}_i^{-1}(t)|^{\frac{n_i-p-2}{2}}   \left|\sum_{i=1}^{g} n_i\hat{\Lambda}_{0i} \left[I_p-t^2\left\{ (p+1)I_p - \tr(\hat{\Lambda}_i^{-1}\hat{\Lambda}_{0i})I_p - \hat{\Lambda}^{-1}_i\hat{\Lambda}_{0i}\right\} \right]\right|^{1/2}.
	\label{density for diffdir}
\end{align}
In this  case,   the directional $p$-value is not expected to be  exactly uniformly distributed under $H_\psi$, since there is no exact conditional density of the sufficient statistic that can be  approximated by (\ref{density for diffdir}).

We also  consider  the likelihood ratio test and its modifications proposed by \citet{skovgaard:2001}.  The likelihood ratio statistic (\ref{W}) takes the form  
\begin{align}
	W  =  \sum_{i=1}^g n_i \left(\log|\hat{\Lambda}_{0i}^{-1}| - \log|\hat{\Lambda}_i^{-1}|\right).\nonumber
\end{align}
The statistic $W$ has  approximate $\chi^2_d$ null distribution with degrees of freedom $d=p(g-1)$, when $p$ is fixed.
The expression of the correction factor $\gamma$, given in (S1),  for hypothesis (\ref{hypothesis:meansdiff}) in \citeauthor{skovgaard:2001}'s modifications (\ref{Wstar}) becomes 
\begin{align}
	\gamma = \frac{\left\{\sum_{i=1}^g n_i (\hat{\mu}_0 - \bar{y}_i)^\top \hat{\Lambda}_{0i} (\hat{\mu}_0 - \bar{y}_i) \right\}^{d/2}}{W^{d/2-1} \sum_{i=1}^g n_i (\hat{\mu}_0 - \bar{y}_i)^\top \hat{\Lambda}_i (\hat{\mu}_0 - \bar{y}_i)}\left\{ \prod_{i=1}^g \frac{|\hat{\Lambda}_i|}{|\hat{\Lambda}_{0i}|}\right\}^{\frac{p+2}{2}} \left\{\frac{|\tilde{C}_1|}{|\sum_{i=1}^g n_i\hat{\Lambda}_{0i}|} \right\}^{1/2},\nonumber
\end{align}
where  $\tilde{C}_1 = \sum_{i=1}^g n_i\hat{\Lambda}_{0i} \left\{{\text{tr}\;}\left(\hat{\Lambda}_i^{-1}\hat{\Lambda}_{0i} - I_p\right)I_p + \hat{\Lambda}^{-1}_i\hat{\Lambda}_{0i}\right\}$.

If $g=2$,  hypothesis (\ref{hypothesis:meansdiff})  reduces to the  multivariate  Behrens-Fisher problem. 
Then, we can compare the directional test also  with the procedure proposed by \citet{nel1986}, which  was shown to maintain a reasonable empirical type I error.  The  test by \citet{nel1986} is based on the quantity 
\begin{align}
	T^{*2} = (\bar{y}_1-\bar{y}_2)^\top \left(\frac{S_1}{n_1} + \frac{S_2}{n_2}\right)^{-1} (\bar{y}_1 - \bar{y}_2),\nonumber
\end{align}
where $S_i = \frac{n_i}{n_i-1} \hat{\Lambda}_i^{-1}$, $i=1,2$. 
The  statistic $\frac{\upsilon - p+1 }{p\upsilon} T^{*2}$ under $H_\psi$ has approximate  $F$-distribution with degrees of freedom $(p,\upsilon-p+1)$, where
\begin{align}
	\upsilon = \frac{{\text{tr}}\left\{\left(\frac{S_1}{n_1} + \frac{S_2}{n_2}\right)\left(\frac{S_1}{n_1} + \frac{S_2}{n_2}\right) \right\} + \left\{{\text{tr}}\left(\frac{S_1}{n_1} + \frac{S_2}{n_2}\right)\right\}^2}{ \left.\left[ {\text{tr}}\left\{\left(\frac{S_1}{n_1}\right) \left(\frac{S_1}{n_1}\right)\right\} + \left\{{\text{tr}}\left(\frac{S_1}{n_1}\right)\right\}^2 \right]\right/{(n_1-1)} + \left[{\text{tr}}\left\{\left(\frac{S_2}{n_2}\right) \left.\left(\frac{S_2}{n_2}\right)\right\} + \left\{{\text{tr}}\left(\frac{S_2}{n_2}\right)\right\}^2 \right]\right/{(n_2-1)} }.\nonumber
\end{align} 
See  \citet[][Section 3.9]{rencher1998multivariate} for more details.

All the different solutions will be evaluated by means of simulation studies in Section \ref{section:simulation dirmeandiff}.

\section{Simulation studies}\label{section:simulation studies}

\subsection{Homoscedastic one-way MANOVA}	\label{section simulation:dirmean}

The  performance  of the directional test  for hypothesis (\ref{hypothesis:dirmean}) in the high dimensional multivariate normal framework is here assessed via Monte Carlo simulations based on $10,000$ replications. The exact directional test is compared with the chi-square approximations for $W$, $W_{BC}$, $W^*$, $W^{**}$, and with the normal approximation for the central limit theorem test  proposed by \citet{He2020B}, specifically for the high dimensional setting. The  six tests are evaluated in terms of empirical size. 

Groups of size $n_i, i \in \{1,\dots, g\}$,  are  generated from a $p$-variate standard normal distribution $N_p(0_p, I_p)$ under the null hypothesis.  For each simulation experiment, we show  results for $p = \kappa n_i$ with $\kappa \in \{0.05, 0.1,0.3,$ $0.5,0.7,0.9\}$  and $n_i \in \{100, 500\}$. Throughout, we set $g = 3$ and $n_1= n_2 =n_3$.  In addition, we also consider some extreme settings with $\kappa \in \{1.0,1.5,2.0,2.5,2.9\}$ in which the chi-square approximations for  $W$, $W_{BC}$, $W^*$, $W^{**}$ break down very fast (see Table \ref{table:type I error same}).

Additional empirical results for the large sample size $n_i=1000$, and different values of $p$ and $g=30$ are reported in the Supplementary Material, which shows that the directional $p$-value maintains high accuracy.    


\begin{figure}[t]
	\centering
	\captionsetup{font=footnotesize}
	\subfigure{
		\begin{minipage}[b]{.5\linewidth}
			\centering
			\includegraphics[scale=0.32]{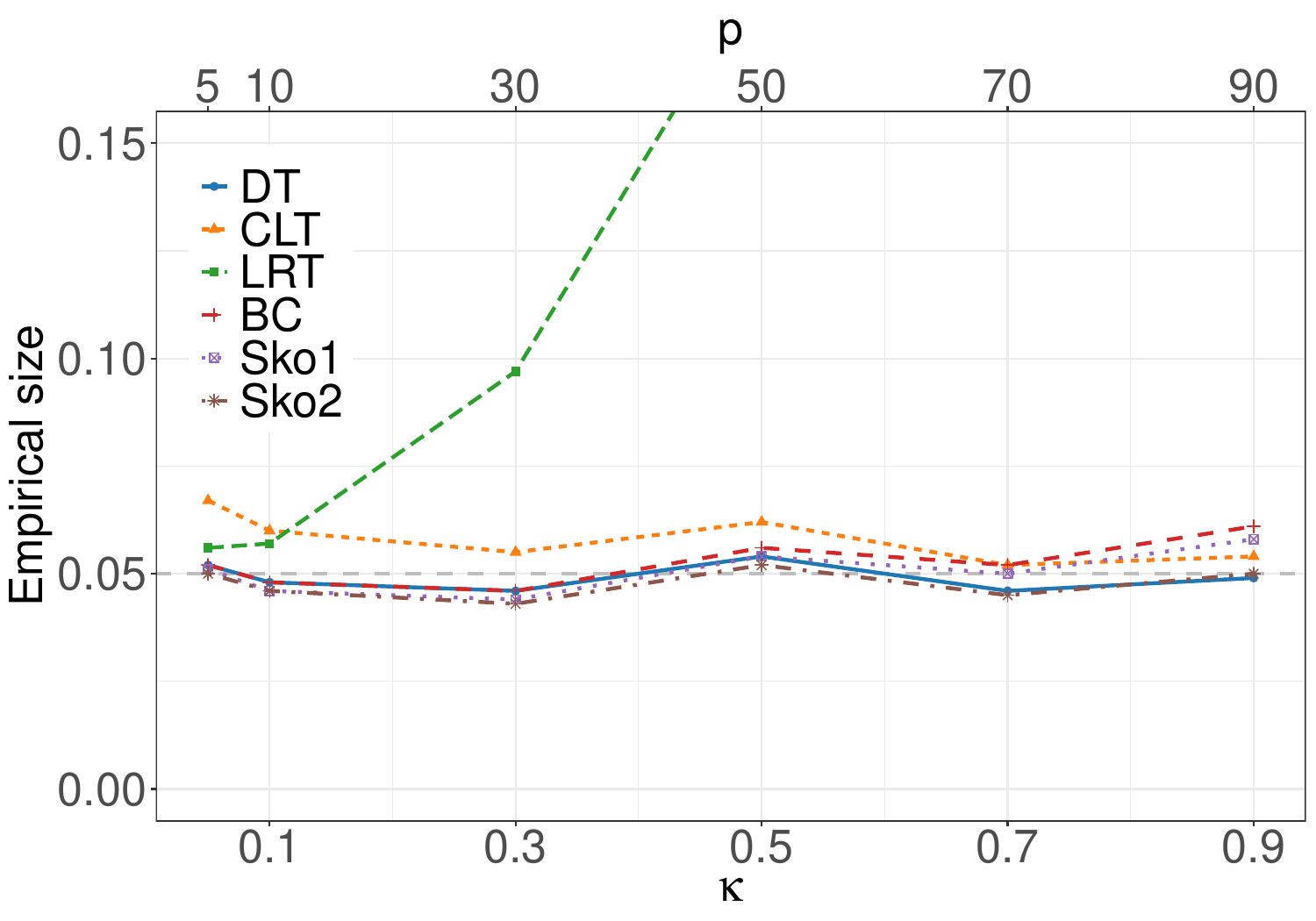}
		\end{minipage}
		\begin{minipage}[b]{.5\linewidth}
			\centering
			\includegraphics[scale=0.32]{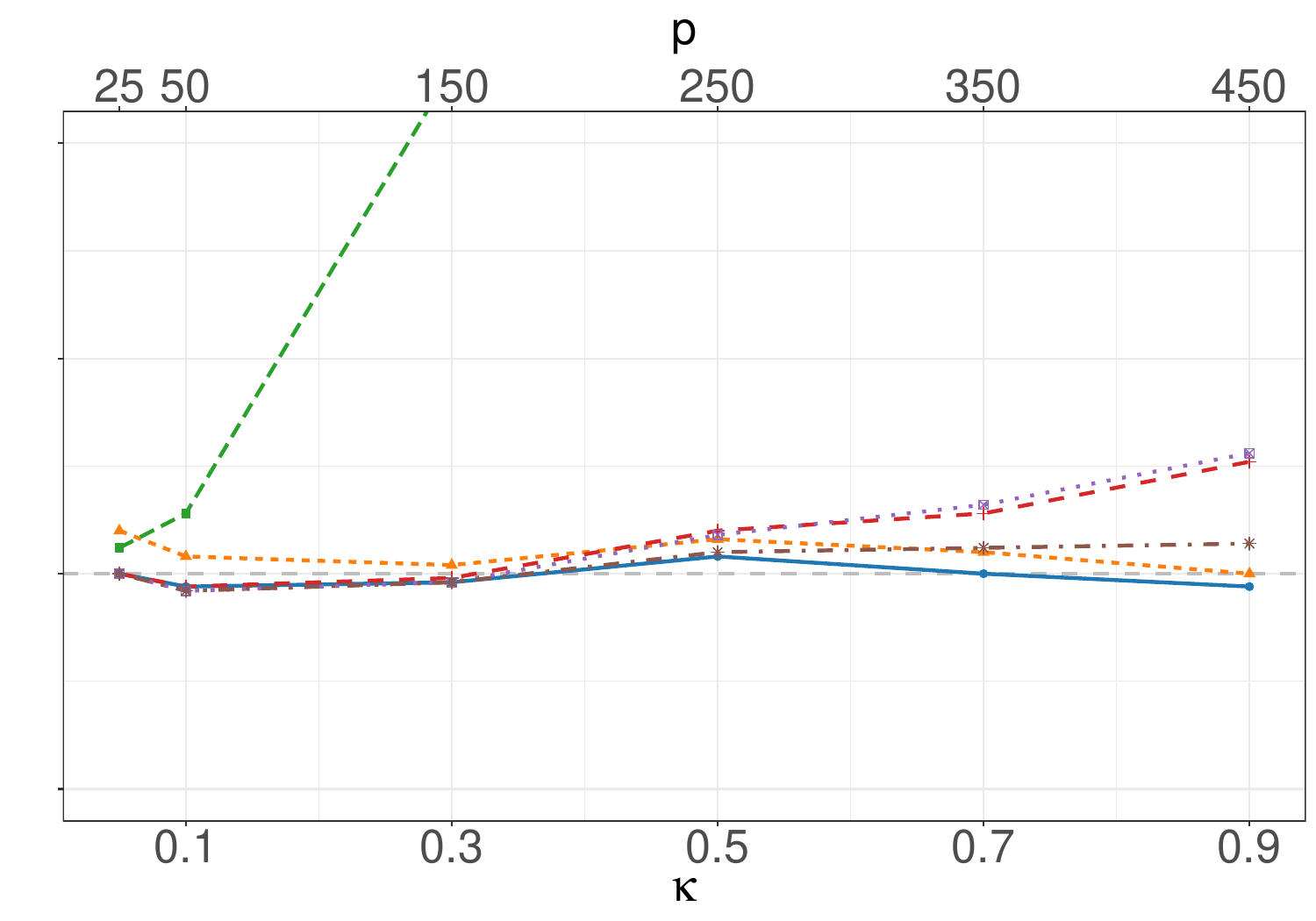}
		\end{minipage}
	}
	\caption{
		Empirical size of  the directional test (DT), central limit theorem test (CLT), likelihood ratio test (LRT), Bartlett corrected test (BC) and two \citeauthor{skovgaard:2001}'s modifications (Sko1 and Sko2) for hypothesis (\ref{hypothesis:dirmean}) with $g=3$, at nominal level $\alpha = 0.05$ given by the dashed gray horizontal line. The left  and right panels correspond to $n_i = 100, 500$, respectively. 
	}
	\label{fig:type I error same}
\end{figure}


The empirical size, i.e. the actual probability of type I error, at nominal  level $\alpha =0.05$ based on the null distribution of the various statistics is reported in Figure \ref{fig:type I error same}.
The  directional $p$-value performs very well across all different values of $p$, confirming the exactness result in Theorem \ref{theorem1},  while the central limit theorem test is less accurate when $p$ is small. The test based on  $W$ breaks down in all settings. Instead,  the Bartlett corrected test $W_{BC}$ proves accurate for moderate values of $p$, as seen in \citet{He2020B}.  
However, the chi-square approximations for $W_{BC}$, $W^*$ and $W^{**}$ get unreliable as $p$ grows. Table  \ref{table:type I error same} reports the empirical size in some extreme settings with $n_i \le p \le \sum_{i=1}^g n_i -g-1$. The results show that the chi-square approximations for $W$, $W_{BC}$, $W^*$ and $W^{**}$  do not work in such extreme settings, while the directional test is still  accurate. The central limit theorem test, instead, {becomes too liberal} with large $p$.

\setlength{\tabcolsep}{1.7mm}{
	\begin{table}[t]
		\centering
		\renewcommand\arraystretch{0.5}
		\caption{Empirical size of  the directional test (DT), central limit theorem test (CLT), likelihood ratio test (LRT), Bartlett corrected test (BC) and two \citeauthor{skovgaard:2001}'s modifications  (Sko1 and Sko2) for hypothesis (\ref{hypothesis:dirmean})  with $g=3$, $p = \kappa n_i$ and $n_i = 100$, at nominal level $\alpha = 0.05$. }
		\medskip
		\begin{tabular}{ccccccc}
			\toprule[0.09 em]
			\midrule[0 em]
			$\kappa$ ($p$) & DT & CLT & LRT & BC & Sko1 & Sko2 \\ 
			\midrule
			1.0 (100) & 0.048 & 0.053 & 0.677 & 0.066 & 0.064 & 0.050 \\ 
			1.5 (150) & 0.048 & 0.055 & 0.991 & 0.133 & 0.131 & 0.073 \\ 
			2.0 (200) & 0.046 & 0.052 & 1.000 & 0.401 & 0.402 & 0.145 \\ 
			2.5 (250) & 0.052 & 0.057 & 1.000 & 0.967 & 0.966 & 0.567 \\ 
			2.9 (290) & 0.048 & 0.060 & 1.000 & 1.000 & 1.000 & 0.999 \\ 
			\bottomrule[0.09 em]
		\end{tabular}
		\label{table:type I error same}
\end{table}}

\subsection{Heteroscedastic one-way MANOVA}	
\label{section:simulation dirmeandiff}

The performance  of the directional test for hypothesis (\ref{hypothesis:meansdiff}) in a moderate dimensional multivariate normal framework with $p=O(n_i^\tau)$, $\tau \in (0,1)$ is here evaluated  via Monte Carlo simulations based on 10,000 replications. The directional test is again compared with the chi-square approximations for the likelihood ratio test and its modifications proposed by \citet{skovgaard:2001}. When $g=2$, we also  consider  the $F$-approximation for the Behrens-Fisher test $T^{*2}$ by \citet{nel1986}.  The different testing approaches are evaluated in terms of empirical size.  The results for $g=5, 30$ are reported in Section S2.3 of 
the Supplementary Material and are in line with the ones  discussed below. 

We generate the data matrix $y_i$ as $n_i$ independent replications from a  multivariate normal distributions $N_p(\mu_i,\Lambda_i^{-1}), i\in \{1,\dots,g\}$. Under the null hypothesis $H_\psi$, we  set $\mu_i =0_p$ and use an  autoregressive structure for the covariance matrices, i.e. $\Lambda_i^{-1} = (\sigma_{jl})_{p \times p} = (\rho_i^{|j-l|})_{p\times p}$, with the $\rho_i$ chosen as  an  equally-spaced sequence from 0.1 to 0.9 of length $g$. In particular, when $g=2$, $\Lambda_1^{-1} = (0.1^{|j-l|})_{p\times p}$ and $\Lambda^{-1}_2 = (0.9^{|j-l|})_{p \times p}$.  We show results for $p = \lceil n_i\rceil^\tau$ with $\tau = j/24$, $j \in \{6,\dots,22\}$, $n_i \in \{100, 500\}$ and  $g = 2$. 


\begin{figure}[t]
	\centering
	\captionsetup{font=footnotesize}
	\subfigure{
		\begin{minipage}[b]{.5\linewidth}
			\centering
			\includegraphics[scale=0.32]{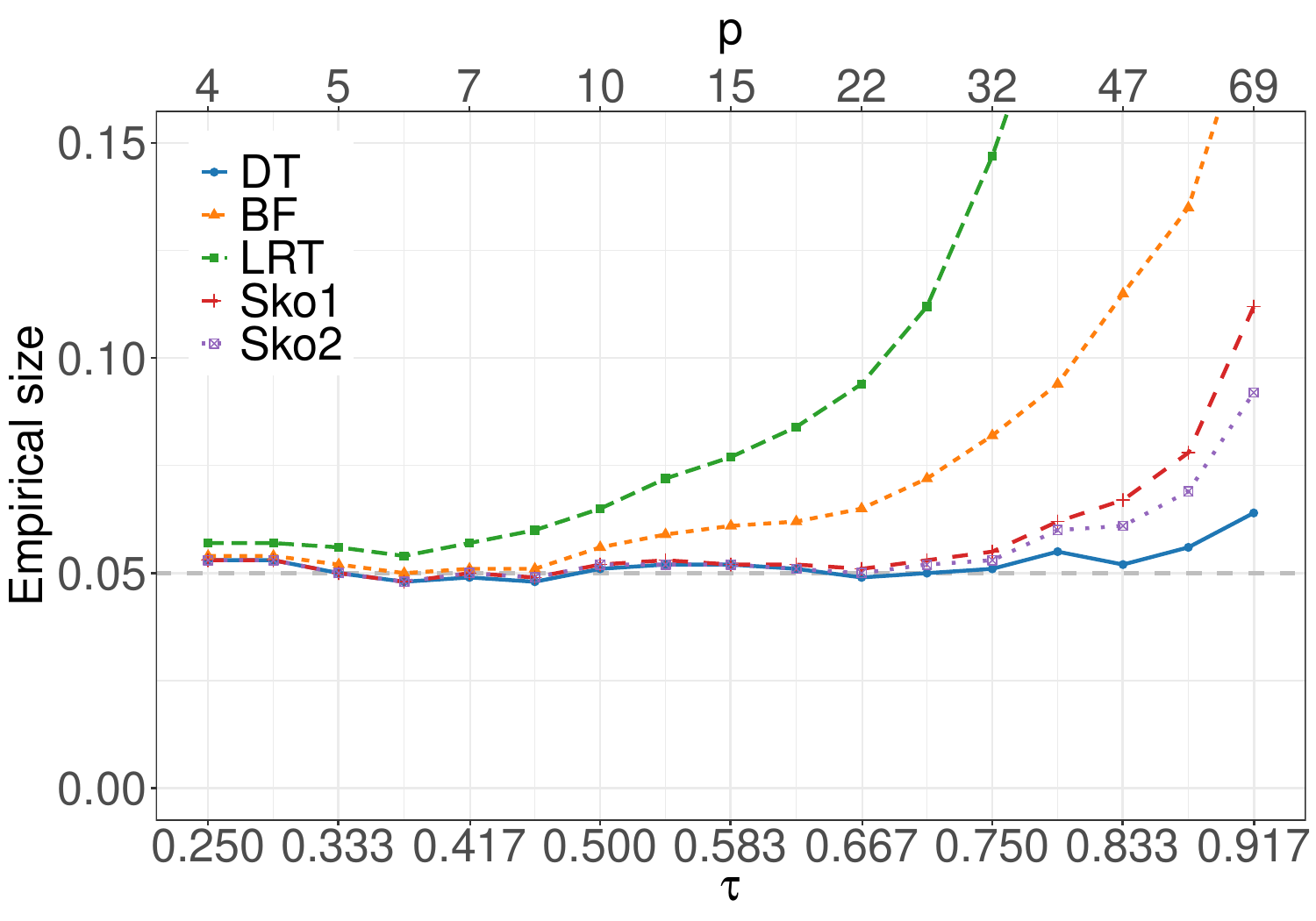}
		\end{minipage}
		\begin{minipage}[b]{.5\linewidth}
			\centering
			\includegraphics[scale=0.32]{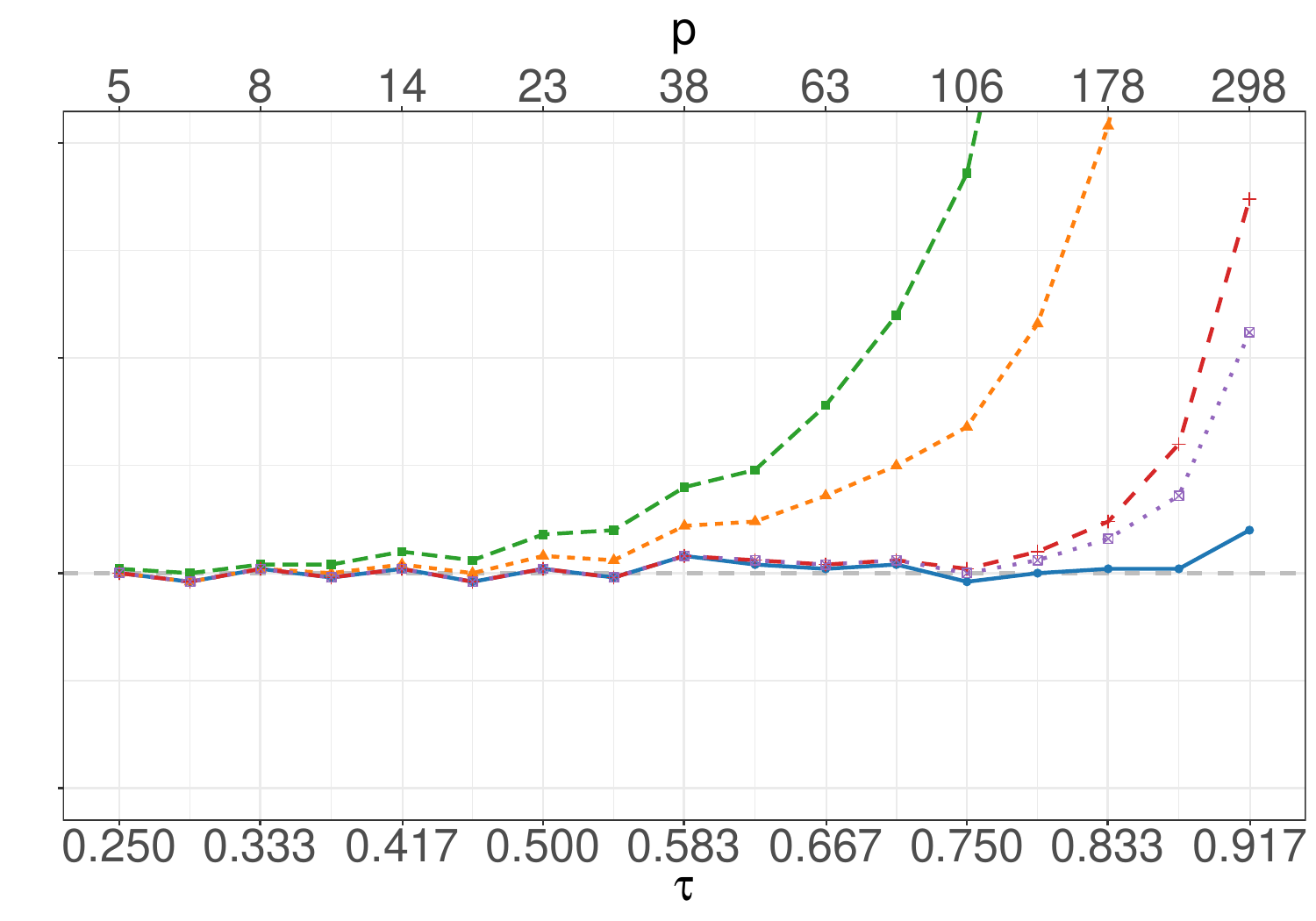}
		\end{minipage}
	}
	
	\caption{
		Empirical size of the directional test (DT), Behrens-Fisher test (BF) \citep{nel1986}, likelihood ratio test (LRT),  and two \citeauthor{skovgaard:2001}'s modifications (Sko1 and Sko2) for hypothesis (\ref{hypothesis:meansdiff}) with {$g=2$}, at nominal level $\alpha = 0.05$  given by the dashed gray horizontal line. The left  and right panels correspond to $n_i = 100, 500$, respectively.
	}
	\label{fig:type I error full}
\end{figure}


Figure \ref{fig:type I error full}	  reports the empirical size at the nominal level $\alpha = 0.05$ under the null hypothesis for various statistics. The directional test is always more reliable than its competitors in terms of the empirical size, even if in this framework it is not exact.  \citeauthor{skovgaard:2001}'s modifications are not as accurate when $p$ is large. Moreover,  
the likelihood ratio test and the Behrens-Fisher test  
break down quickly when $p = \lceil n_i\rceil^\tau$, even for moderately large values of $\tau$.
Additional simulation results are shown in Table \ref{table:type I error} for different structures of the covariance matrices: (I) identity matrix, $\Lambda_i^{-1} = I_p$;
(II) compound symmetric matrix, $\Lambda^{-1}_i = (1-\rho_i) I_p + \rho_i 1_p 1_p^\top$, with  the same values of $\rho_i$ as in the autoregressive  structure.
Table \ref{table:type I error}  reports the empirical size at the nominal level $\alpha = 0.05$ for $n_i = 100$ and $p = \lceil n_i\rceil^\tau$, with covariance structures (I) and (II) and $g=2, 30$. As expected, the likelihood ratio test performs in general very poorly, especially for larger values of $g$.  \citeauthor{skovgaard:2001}'s modifications are much more accurate but not as much as the directional test, whose excellent performance is confirmed. When $g=2$, the approximate solution for the Behrens-Fisher test seems reliable in the identity covariance matrix case, but not under compound symmetry.   

\setlength{\tabcolsep}{0.7mm}{
	\begin{table}[t]
		\centering
		\renewcommand\arraystretch{0.5}
			\caption{Empirical size of the directional test (DT), likelihood ratio test (LRT), two \citeauthor{skovgaard:2001}'s modifications (Sko1 and  Sko2), and Behrens-Fisher test (BF) (only for $g=2$) for hypothesis (\ref{hypothesis:meansdiff})  at nominal level $\alpha = 0.05$, with $n_i = 100$, $p = n_i^\tau$,  $\Lambda_i^{-1}=I_p$ (I) and  $\Lambda^{-1}_i = (1-\rho_i) I_p + \rho_i 1_p 1_p^\top$(II), with the $\rho_i$ chosen as  an  equally-spaced sequence from 0.1 to 0.9 of length $g$.}
		\medskip
		\begin{tabular}{cccccccccccc}
			\toprule[0.09 em]
			$\Lambda_i^{-1}$	&$\tau$ ($p$)&  \multicolumn{5}{c}{$g=2$}&& \multicolumn{4}{c}{$g=30$}\\
			\cline{3-7} \cline{9-12}
			\midrule[0 em]
			&& DT & BF & LRT & Sko1 & Sko2& &DT & LRT &Sko1 & Sko2 \\
			\midrule
			(I)	& 10/24 (7) & 0.051 & 0.051 & 0.058 & 0.051 & 0.051 &  & 0.048 & 0.118 & 0.049 & 0.048 \\ 
			& 12/24 (10) & 0.050 & 0.050 & 0.059 & 0.050 & 0.050 &  & 0.048 & 0.188 & 0.050 & 0.048 \\ 
			& 14/24 (15) & 0.050 & 0.050 & 0.069 & 0.050 & 0.050 &  & 0.048 & 0.382 & 0.051 & 0.048 \\ 
			& 16/24 (22) & 0.048 & 0.048 & 0.079 & 0.049 & 0.048 &  & 0.049 & 0.754 & 0.060 & 0.052 \\ 
			& 18/24 (32) & 0.053 & 0.052 & 0.110 & 0.054 & 0.053 &  & 0.050 & 0.993 & 0.087 & 0.061 \\ 
			& 20/24 (47) & 0.046 & 0.045 & 0.180 & 0.050 & 0.047 &  & 0.054 & 1.000 & 0.223 & 0.099 \\
			& 22/24 (69) & 0.054 & 0.052 & 0.377 & 0.066 & 0.059 &  & 0.048 & 1.000 & 0.928 & 0.404\\
			\vspace{-0.8em}
			&&&&&&&&&&&\\
			
			(II)	&   10/24 (7) & 0.050 & 0.054 & 0.059 & 0.050 & 0.050 &  & 0.049 & 0.118 & 0.049 & 0.048 \\ 
			& 12/24 (10) & 0.053 & 0.058 & 0.065 & 0.053 & 0.053 &  & 0.048 & 0.187 & 0.049 & 0.048 \\ 
			& 14/24 (15) & 0.048 & 0.062 & 0.076 & 0.050 & 0.049 &  & 0.049 & 0.383 & 0.052 & 0.050 \\ 
			& 16/24 (22) & 0.049 & 0.075 & 0.104 & 0.051 & 0.049 &  & 0.049 & 0.755 & 0.060 & 0.051 \\ 
			& 18/24 (32) & 0.050 & 0.104 & 0.171 & 0.057 & 0.054 &  & 0.050 & 0.993 & 0.085 & 0.059 \\ 
			& 20/24 (47) & 0.049 & 0.154 & 0.345 & 0.071 & 0.062 &  & 0.053 & 1.000 & 0.225 & 0.099 \\ 
			& 22/24 (69) & 0.081 & 0.296 & 0.766 & 0.174 & 0.136 &  & 0.048 &   1.000  & 0.932 & 0.408 \\
			\bottomrule[0.09 em]
		\end{tabular}
		\label{table:type I error}
\end{table}}

\section{Application}\label{sec:app}
We consider two real-data applications of the proposed methodology. One regards microarray data of patients suffering from  breast cancer, and can be found in Section S5 of the Supplementary Material, while the other is described in detail here.
The dataset  is taken from a study on how  $p=5$ characteristics of the blood (red cell count, hemoglobin concentration, hematocrit, white cell count, and plasma ferritin concentration) vary with sport, body size, and sex of the athletes \citep{telford1991sex}. Our goal is to investigate whether the blood characteristics of athletes differ between  endurance-related events and  power-related events. The data can be found in the function $\texttt{ais}$ of the R  package  \texttt{DAAG}  \citep{daag24}. We performed several tests on this dataset to study the difference in blood parameters between athletes performing  endurance-related  and  power-related sports. Typically, the conclusions drawn by the directional approach align with those from its best competitors. As an illustration, we report here some results involving two hypotheses:  (i) testing the  equality of blood characteristics between athletes competing in an  endurance-related event (400-meter track event, $T_{400}$) and those in a power-related event (sprint, $T_{Sprint}$), i.e. $H_{0}: \mu_{T_{400}} = \mu_{T_{Sprint}}$; and (ii) testing the equality of blood characteristics for three groups of athletes that perform endurance-related sports ($Row$, $Swim$, and $T_{400}$), i.e. $H_{0}: \mu_{T_{400}} = \mu_{Row} = \mu_{Swim}$. Here we focus exclusively on data from  male athletes as each blood variable  is significantly different between genders \citep{telford1991sex}. We also consider the variable plasma ferritin concentration on  the logarithmic scale. The sample sizes of $T_{400}$, $T_{Sprint}$, $Row$ and $Swim$  groups are respectively 18, 11, 15, 13. 
All the testing procedures under analysis can then be used to perform one-way MANOVA, and the corresponding  $p$-values are shown in Table \ref{table:p-value application ais}.

For  hypothesis (i), the  $p$-values in Table \ref{table:p-value application ais} of the directional test (homoscedastic: 0.060; heteroscedastic: 0.097), the Bartlett corrected test (homoscedastic: 0.059), the two  Skovgaard's modifications (homoscedastic: 0.078 and 0.083; heteroschedastic: 0.089 and 0.093) and the  Behrens-Fisher test (heteroscedastic: 0.179) 
claim that the average blood characteristics levels for male athletes competing in  $T_{400}$ and $T_{Sprint}$  events are not significantly different at a 5\% level.   Instead,   the $p$-values  of the central limit theorem test (homoscedastic: 0.036) and the likelihood ratio test (homoscedastic: 0.027; heteroscedastic: 0.034) both lead to a different conclusion.  Note that the central limit theorem test is not reliable with small $p$ and $n_i$, while the likelihood ratio test is not valid if  $\kappa$ and $\tau$ are large. To assess the validity of the various testing procedures in this setting, we can also examine the empirical size of the various tests.  We generate  $R=10,000$ samples with sizes  $n_{T_{400}} = 18$ and $n_{T_{Sprint}} = 11$ from a normal distribution with mean vector $\mu_{T_{400}} = \mu_{T_{Sprint}}=(n_{T_{400}} \bar{y}_{T_{400}} + n_{T_{Sprint}}\bar{y}_{T_{Sprint}}) / (n_{T_{400}} + n_{T_{Sprint}})$ where $\bar{y}_{T_{400}}$ and $\bar{y}_{T_{Sprint}}$ are the group sample means, and covariance matrix $\Lambda^{-1} = \hat{\Lambda}^{-1}_0$,  the pooled sample covariance matrix  for the homoscedastic case, and covariance matrices $\Lambda^{-1}_i = \hat{\Lambda}_i^{-1}$ sample covariance matrices of each group, for the heteroscedastic case. The simulation results for  hypothesis (i) reported in Table \ref{table:p-value simulation ais} show that the directional test is more reliable than its competitors, having empirical size  closer to the nominal level $0.05$. The central limit theorem test and likelihood ratio test suffer from large bias.  These results confirm the simulation outputs in Section \ref{section simulation:dirmean}.

\setlength{\tabcolsep}{1.5mm}{
	\begin{table}[t]
		\centering
		\renewcommand\arraystretch{0.5}
		\caption{$p$-values of  the directional test (DT), central limit theorem test (CLT),  likelihood ratio test (LRT), Bartlett corrected test (BC) , two \citeauthor{skovgaard:2001}'s modifications (Sko1 and  Sko2), and Behrens-Fisher test (BF)  for testing whether the blood characteristics of athletes differ between endurance-related and power-related events  at nominal level $\alpha = 0.05$, with $n_i\in\{18,11,15,13\}$ and $p = 5$. }
		\medskip
		\begin{tabular}{cccccccccccccc}
			\toprule[0.09 em]
			Hypothesis	&$g$ & \multicolumn{6}{c}{Homoscedastic}& & \multicolumn{5}{c}{Heteroscedastic}\\
			\cline{3-8} \cline{10-14}
			\midrule[0 em]
			&& DT & CLT & LRT & BC & Sko1 & Sko2 &&  DT & BF & LRT  & Sko1 & Sko2  \\
			\midrule
			(i)	& 2 & 0.060 & 0.036 & 0.027 & 0.059 & 0.078 &   0.083 && 0.097 & 0.179 & 0.034 & 0.089 & 0.093\\ 
			(ii)& 3 & 0.092 & 0.069 & 0.045 & 0.084 & 0.101 &   0.105 && 0.157 &-& 0.062 & 0.152 & 0.158\\ 
			\bottomrule[0.09 em]
		\end{tabular}
		\label{table:p-value application ais}
\end{table}}

\setlength{\tabcolsep}{1.5mm}{
	\begin{table}[t]
		\centering
		\renewcommand\arraystretch{0.5}
		\caption{Empirical size of  the directional test (DT), central limit theorem test (CLT),  likelihood ratio test (LRT), Bartlett corrected test (BC) , two \citeauthor{skovgaard:2001}'s modifications (Sko1 and  Sko2), and Behrens-Fisher test (BF)  for testing whether the blood characteristics of athletes differ between endurance-related and power-related events  at nominal level $\alpha = 0.05$, with $n_i\in\{18,11,15,13\}$ and $p = 5$. }
		\medskip
		\begin{tabular}{cccccccccccccc}
			\toprule[0.09 em]
			Hypothesis	&$g$ & \multicolumn{6}{c}{Homoscedastic}& & \multicolumn{5}{c}{Heteroscedastic}\\
			\cline{3-8} \cline{10-14}
			\midrule[0 em]
			&& DT & CLT & LRT & BC & Sko1 & Sko2 &&  DT & BF & LRT  & Sko1 & Sko2  \\
			\midrule
			(i)	& 2 &  0.049 &  0.070 & 0.092 &  0.050 &0.038 &  0.034 && 0.053 &0.017 & 0.111 & 0.055 & 0.053\\ 
			(ii)& 3&  0.045 &  0.062 & 0.087 & 0.046 & 0.039 &  0.038  && 0.049 &-& 0.125 & 0.050 & 0.047\\ 
			\bottomrule[0.09 em]
		\end{tabular}
		\label{table:p-value simulation ais}
\end{table}}

For hypothesis (ii),  all the tests used under heteroscedasticity in Table  \ref{table:p-value application ais} give $p$-values directional: 0.097; likelihood ratio: 0.062; Skovgaard's: 0.152 and 0.158) 
concluding that differences in mean blood characteristics among endurance-related events ($T_{400}$, $Row$ and $Swim$) for male athletes are not statistically significant at a 5\% level.  In the   homoscedastic case, only the likelihood ratio test yields a significant $p$-value (0.045). Conversely, the other test statistics  have $p$-values (directional: 0.092; central limit theorem: 0.069; Bartlett corrected: 0.084; Skovgaard's: 0.101 and 0.105)  that agree with the hypothesis of no significant difference 
among the three groups.  Even in this case, a simulation study  can be performed in order to  verify the reliability of the considered methods in the specific setting. Under the null hypothesis, we generate $10,000$ samples with sizes $n_{T_{400}}=18$, $n_{Row}=15$ and $n_{Swim}=13$   from a normal  distribution with mean vector  $\mu_{T_{400}} = \mu_{Row} = \mu_{Swim} =(n_{T_{400}} \bar{y}_{T_{400}} + n_{Row}\bar{y}_{Row} + n_{Swim}\bar{y}_{Swim}) / (n_{T_{400}} + n_{Row} + n_{Swim})$ where $\bar{y}_{T_{400}}$, $\bar{y}_{Row}$ and $\bar{y}_{Swim}$ are the group sample means, and covariance matrix $\Lambda^{-1} = \hat{\Lambda}_0^{-1}$,  the pooled sample covariance matrix under the homoscedasticity, and covariance matrices $\Lambda^{-1}_i = \hat{\Lambda}_i^{-1}$,  sample covariance matrices under the heteroscedasticity. Simulation results are reported in Table \ref{table:p-value simulation ais}. In the homoscedastic case, the two Skovgaard's modifications are  conservative. The central limit theorem test and the likelihood test, on the other hand, have slightly inflated empirical size. Hence, the $p$-values of these four methods  are not reliable in this scenario. In contrast, the directional test and the Bartlett corrected test well control the  type I error.  In the heteroscedastic case, the likelihood ratio test is  overly inflated, whereas the remaining procedures, including the directional test, exhibit a good performance.

\section{Conclusions} \label{section:conclusion}

In this paper, we have developed the directional test for one-way MANOVA problems  when the data dimension is comparable with the sample size. The directional $p$-value has been proved exactly uniformly distributed under the null hypothesis in a high dimensional setting, provided that $\sum_{i=1}^{g}n_i \ge p+g+1$,  when testing the equality of normal mean vectors with identical covariance matrix. This finding is supported by our numerical studies. 
In addition, simulations in moderate dimensional scenarios under the  heteroscedastic one-way MANOVA framework  indicate that the directional test outperforms its competitors  in terms of empirical null distribution,  even in the more general setting with  different covariance matrices. 
While formal conditions for the validity of the various test statistics in this scenario could perhaps be developed using
recent results on the accuracy of the saddlepoint approximation in moderate dimensional regimes   \citep{tang2021laplace}, we believe further work is needed in this area to refine such methods under more complex settings.

It is also worth mentioning that, given the inconclusive outcome of the power analysis of likeli-hood-based tests for high dimensional normal distributions included in \citet[Sect. 5.3]{huang2021directional}, here we have deliberately chosen not to compare the performance of the tests under a few arbitrary alternative hypotheses. 
	Finally, although additional numerical investigations in Sections S3.1--S3.4 and S4.3 of the Supplementary Material show that  all the testing solutions examined are reasonably robust to misspecification of the underlying multivariate normal model, different and/or more substantial deviations from normality could be considered in future research for gaining a more comprehensive picture on this aspect.

\section*{Supplementary Material}

Supplementary Material provides auxiliary computational details, formal proofs, and extensive additional simulation studies covering both homoscedastic and heteroscedastic high and moderate dimensional one-way MANOVA, along with an applicaiton results from breast cancer microarray data. 
The \texttt{R} code to reproduce all numerical results in the paper are available at \url{https://github.com/stat-cz/DirTestMANOVA}.

\bibliographystyle{apalike}
\bibliography{dirmean}

\newpage
\setcounter{page}{1}
\setcounter{section}{0}
\setcounter{figure}{0}
\setcounter{table}{0}

\renewcommand{\thefigure}{S\arabic{figure}}
\renewcommand{\thetable}{S\arabic{table}}
\renewcommand{\theequation}{S\arabic{equation}}
\renewcommand{\thesection}{S\arabic{section}}
\renewcommand{\thesubsection}{S\arabic{section}.\arabic{subsection}}

\section*{Supplementary material to Directional testing for one-way MANOVA in divergent dimensions}

\begin{abstract}
Section S1 provides auxiliary computational results for  the correction factor $\gamma$ in  \citeauthor{skovgaard:2001}'s  modifications, as well as  the  observed information matrix operator used for testing hypothesis (8) of the main text.  Section S2 presents the proofs for Lemma 1, Theorem 1 and Proposition 1. Section S3 reports extensive simulation studies for the homoscedastic  one-way MANOVA; specifically, Sections S3.1--S3.4 examine varing dimenions $p$ and the group numbers $g$. Sections S3.5--S3.6 assess robustness under model misspecification for the high and moderate dimensional setups. Furthermore, Section S4 reports  extensive simulations for the  heteroscedastic one-way MANOVA, including a robustness analysis in Section S4.3. Finally, Section S5 provides results from an additional application to   breast cancer microarray data.
\end{abstract}

\renewcommand{\baselinestretch}{1} 

\section{Computational results} \label{sec:appA}

\subsection{Correction factor $\gamma$ in \citeauthor{skovgaard:2001}'s modifications}\label{section:correction factor}

The general expression for $\gamma$ appearing in \citeauthor{skovgaard:2001}'s modified versions of W (3) in the main text  can be written as
\begin{eqnarray}
	\gamma =\frac{\left\{(s-s_\psi)^\top J_{\varphi\varphi}(\hat{\varphi}_\psi)^{-1} (s-s_\psi)\right\}^{d/2}}{W^{d/2-1}(\hat{\varphi}-\hat{\varphi}_\psi)^\top (s-s_\psi)} \left\{\frac{|{J}_{\varphi \varphi}(\hat{\varphi}_\psi)|}{|{J}_{\varphi \varphi}(\hat{\varphi})|}\right\}^{1/2} \left\{\frac{|{J}_{\lambda \lambda}(\hat{\varphi}_\psi)|}{|{I}_{\lambda \lambda}(\hat{\varphi}_\psi)|}\right\}^{1/2},
	\label{skovgaard's correction factor}
\end{eqnarray}
where $J_{\lambda \lambda}$ and $I_{\lambda \lambda}$ are the blocks of  observed and expected information matrices, respectively, relative to the nuisance parameter $\lambda$.
For calculating the $p$-value, the quantity (\ref{skovgaard's correction factor}) is evaluated at $s = s^0 = 0_q$, corresponding to the observed data point $y^0$. In particular, if $\lambda$ is a component or a linear function of the canonical parameter, as in Section 3.1 of  the main text, the last factor in (\ref{skovgaard's correction factor}) equals 1.

\subsection{Observed information matrix }\label{section:J for same}
When testing hypothesis (8) in the main text, the saddlepoint approximation \citep[][Section 10]{salvan:1997} and  the correction factor $\gamma$ in (3) require the computation of the observed information matrix $J_{\varphi\varphi} (\varphi)$. This has the  block form
\begin{align}
	J_{\varphi\varphi} (\varphi)
	&=\left(\begin{array}{cc}
		J_{\xi\xi} & J_{\xi \text{vech}(\Lambda)}\\
		J_{\text{vech}(\Lambda)\xi} & J_{ \text{vech}(\Lambda) \text{vech}(\Lambda)}
	\end{array}\right), \nonumber
\end{align}
\begin{align}
	\text{with} \; J_{\xi\xi}
	=\left(\begin{array}{cccc}
		n_1\Lambda^{-1} & 0_{p\times p}   & \cdot & 0_{p\times p} \\
		0_{p\times p}   & n_2\Lambda^{-1} & \cdot & 0_{p\times p} \\
		\vdots          & \vdots          & \vdots& \vdots\\
		0_{p\times p}   &  0_{p\times p}  & \cdot & n_g\Lambda^{-1}
	\end{array}\right), \;
	J_{\text{vech}(\Lambda)\xi}^\top = J_{\xi \text{vech}(\Lambda)} = \left(\begin{array}{c}
		-n_1(\xi_1^\top\Lambda^{-1}\otimes\Lambda^{-1})D_p\\
		-n_2(\xi_2^\top\Lambda^{-1}\otimes\Lambda^{-1})D_p\\
		\vdots  \\
		-n_k(\xi_g^\top\Lambda^{-1}\otimes\Lambda^{-1})D_p
	\end{array}\right)  \nonumber
\end{align}
and $J_{ \text{vech}(\Lambda) \text{vech}(\Lambda)} = \sum_{i=1}^g \frac{n_i}{2} D_p^\top \Lambda^{-1} (I_p + 2\xi_i\xi_i^\top\Lambda^{-1})\otimes\Lambda^{-1}D_p,$ where $\otimes$ denotes the Kronecker product \citep[see, for instance,][]{Lauritzen:1996}. Then, the determinant of $J_{\varphi \varphi}(\varphi)$  can be computed as
\begin{align}
	|J_{\varphi \varphi}(\varphi)| = |J_{\xi\xi}||C_2|,\nonumber
\end{align}
where $C_2 =J_{ \text{vech}(\Lambda) \text{vech}(\Lambda)} - J_{\text{vech}(\Lambda)\xi} J_{\xi\xi}^{-1}J_{\xi \text{vech}(\Lambda)} =  \sum_{i=1}^g \frac{n_i}{2} D_p^\top \Lambda^{-1} (I_p + 2\xi_i\xi_i^\top\Lambda^{-1})\otimes\Lambda^{-1}D_p - \sum_{i=1}^g n_i D_p^\top \Lambda^{-1} $ $\xi_i \xi_i^\top$ $ \Lambda^{-1} \otimes \Lambda^{-1}D_p$ $=\frac{n}{2} D_p^\top  \left(\Lambda^{-1} \otimes \Lambda^{-1}\right) D_p$.

After some algebra, we get 
\begin{align}
	|J_{\varphi \varphi}(\varphi)|  =\left(\prod_{i=1}^g n_i^p\right) |\Lambda^{-1}|^g  n^{\frac{p(p+1)}{2}} 2^{-p} |\Lambda^{-1}|^{p+1} \propto |\Lambda^{-1}|^{p+g+1}.\nonumber
\end{align}

\section{Proofs} \label{sec:appB}

\subsection{Proof of Lemma 1 }\label{section:proof lemma1}
Recall that  $\hat \Lambda^{-1}(t)= \hat{\Lambda}_0^{-1} - t^2 A / n = \hat{\Lambda}^{-1}+ (1-t^2) A/n $ with $\hat \Lambda_0^{-1} = \hat{\Lambda}^{-1} + A/n$, where $A = \sum_{i=1}^g n_i (\bar{y}_i - \bar{y})(\bar{y}_i - \bar{y})^\top$. 
If $t\in [0,1]$ the result is straightforward, since, for all $x\in \mathbb{R}^p$, $x^\top \hat \Lambda^{-1}(t) x = x^\top \hat{\Lambda}^{-1}x/n +(1-t^2) x^\top Ax/n>0$.
Let us focus on the case $t>1$.  We rewrite the estimator $\hat{\Lambda}^{-1}(t)$  as
\begin{align}
	\hat{\Lambda}^{-1}(t) = B_0^\top \left\{ I_p - t^2 (B_0^\top)^{-1} (A/n) B_0^{-1} \right\} B_0,\nonumber
\end{align}
with  $B_0$  such that $\hat{\Lambda}_0^{-1} = B_0^\top B_0$.  According to the eigen decomposition, the matrix $(B_0^\top)^{-1} (A/n) B_0^{-1} = PQP^\top$ with  an orthogonal  matrix $P$ whose columns are eigenvectors of $(B_0^\top)^{-1} (A/n) B_0^{-1}$ and a diagonal matrix $Q$ whose diagonal elements are the eigenvalues of $(B_0^\top)^{-1} (A/n) B_0^{-1}$.  Then, we have $\hat{\Lambda}^{-1}(t) = B_0^\top P \{I_p - t^2 Q\} P^\top B_0$. Therefore, checking that $\hat{\Lambda}^{-1}(t)$ is positive definite is equivalent to checking that $I_p -t^2Q$ is positive definite. Indeed, for all $x\in \mathbb{R}^p, x\ne 0$, then 
\begin{eqnarray}
	x^\top \hat{\Lambda}^{-1} x &=& x^\top B^\top_0 P \{I_p - t^2Q\} P^\top B_0 x\nonumber\\
	&=& \tilde{x}^\top \{I_p - t^2Q\}\tilde{x} >0,\nonumber
\end{eqnarray}
where $\tilde{x} = P^\top B_0x$, with $\tilde{x}\ne 0$ if $x\ne 0$.  

Next, the positive definiteness of $I_p -t^2Q$ should be proved.  It is equivalent to checking that all elements of the diagonal matrix  $I_p -t^2Q = \text{diag}(1-t^2\nu_l)$ are positive, where $\nu_l, l \in \{1,\dots,p\}$, are the eigenvalues of the matrix $(B_0^\top)^{-1} (A/n) B_0^{-1}$. We now need to find out the largest $t$ such that $1-t^2 \nu_l>0$, $l \in \{1,\dots,p\}$. Then, the largest value  of $t$ for which $\hat{\Lambda}^{-1}(t)$ is positive definite equals $t_{\text{sup}} = \sqrt{1/\nu_{(p)}}$, where $\nu_{p}$ is the maximum eigenvalue of the matrix $(B_0^{-1})^\top A B_0^{-1}/n$. 
Therefore, $\hat{\Lambda}^{-1}(t)$ is positive definite in $t \in [0,\sqrt{1/\nu_{(p)}}]$.

\subsection{Proof of Theorem 1} \label{section:proof theorem1}

Suppose $y_i=[y_{i1} \cdots y_{in_i}]^\top$, with $y_{ij} \sim N_p(\mu_i, \Lambda^{-1})$, $i\in \{1,\dots,g\}$ and $j\in \{1,\dots,n_i\}$. The log-likelihood for the canonical parameter $\varphi$ is 
\begin{eqnarray}
	\ell(\varphi;s) = \sum_{i=1}^g -\frac{n_i}{2} \log|\Lambda^{-1}| - \frac{1}{2} \tr(\Lambda B_i) - \frac{n_i}{2} (\bar{y}_i -\mu_i)^\top \Lambda (\bar{y}_i -\mu_i),
	\label{proof:likelihood dirmean}
\end{eqnarray}
where $B_i = y_i^\top y_i - n_i \bar{y}_i\bar{y}_i^\top$. 
The  maximum likelihood estimate  and the constrained  maximum likelihood estimate are respectively $\hat{\varphi} = \{\hat{\xi}^\top_1, \ldots, \hat{\xi}^\top_k, \text{vech}(\hat{\Lambda})^\top\}^\top$ $ = \{\bar{y}_1^\top\hat{\Lambda}, \ldots, \bar{y}_g^\top\hat{\Lambda}, \text{vech}(\hat{\Lambda})^\top\}^\top$  and  $\hat{\varphi}_{\psi} = \{\bar{y}^\top{\hat{\Lambda}_0}, \ldots, \bar{y}^\top{\hat{\Lambda}_0}, $ $ \text{vech}({\hat{\Lambda}_0})^\top\}^\top$.
Evaluating (\ref{proof:likelihood dirmean}) at the unconstrained and constrained maximum likelihood estimates for $\varphi$, we have $\ell(\hat{\varphi};s) =-n/2 \log|\hat{\Lambda}^{-1}| - {np}/2$ and $\ell(\hat{\varphi}_\psi;s) = 2^{-1}\sum_{i=1}^g -{n_i} $ $\log|{\hat{\Lambda}_0}^{-1}| - {n_i}p$, respectively. Then, under the null hypothesis $H_\psi$, 
using the fact that  $|J_{\varphi \varphi}(\hat{\varphi})| \propto |\hat{\Lambda}^{-1}|^{p+g+1}$, the saddlepoint approximation (5) in the main text   is 
%
\begin{eqnarray}
	h(s;\psi) &=&  \left[ \prod_{i=1}^g c_{i1} |\hat{\Lambda}_0^{-1}|^{-\frac{1}{2}} \exp\left\{-\frac{n_i}{2}(\bar{y}_i - \hat\mu_\psi)^\top\hat\Lambda_0 (\bar{y}_i -\hat\mu_\psi)\right\}\right] \nonumber\\
	&& \;\;  \times c_{2} |\hat{\Lambda}_0^{-1}|^{\frac{n-1}{2}} \exp\left\{-\frac{n}{2} \tr(\hat \Lambda_0\hat{\Lambda}^{-1})\right\} |\hat{\Lambda}^{-1}|^{\frac{n-p-g-1}{2}}.
	\label{proof:saddle dirmean}	
\end{eqnarray}
Expression (\ref{proof:saddle dirmean}) equals the exact joint distribution of $\bar{y}_1, \dots, \bar{y}_g$ and $\hat{\Lambda}^{-1}$ if $c_{i1} = (2\pi)^{-p/2}$ and $c_{2} = \left(\frac{n}{2}\right)^{p(n-1)/2}$ $\Gamma_p\left(\frac{n-g}{2}\right)^{-1}$ and with fixed values of $\hat{\mu}_\psi$ and $\hat{\Lambda}_0^{-1}$. Indeed, the constrained maximum likelihood estimates are fixed and equal to their observed values when  we consider  the saddlepoint approximation along the line $s(t)$ under $H_\psi$.  Moreover, we have the unconstrained  maximum likelihood estimates  $\hat{\mu}(t) = t\bar{y}_i + (1-t)\bar{y}$ and  $\hat{\Lambda}^{-1}(t) = \hat{\Lambda}_0^{-1} - t^2A/n$.   Then, the saddlepoint approximation to the conditional distribution of $s(t)$ under $H_{\psi}$ follows from (\ref{proof:saddle dirmean}) and is equal to
\begin{eqnarray}
	h\{s(t);\psi\} &=&  \prod_{i=1}^g c_{i1} |\hat{\Lambda}_0^{-1}|^{-\frac{1}{2}} \exp\left[-\frac{n_i}{2}\{\hat{\mu}_i(t) - \hat\mu_\psi\}^\top\hat\Lambda_0 \{\hat{\mu}_i(t) -\hat\mu_\psi\}\right] \nonumber\\
	&& \;\;\times c_{2} |\hat{\Lambda}_0^{-1}|^{\frac{n-1}{2}} \exp\left[-\frac{n}{2} \tr\{\hat \Lambda_0\hat{\Lambda}(t)^{-1}\}\right] |\hat{\Lambda}(t)^{-1}|^{\frac{n-p-g-1}{2}}	\nonumber\\
	&\propto& \prod_{i=1}^g \exp\left[-\frac{n_i}{2}\{t\bar{y}_i + (1-t)\bar{y} - \bar{y}\}^\top\hat{\Lambda}_0\{t\bar{y}_i + (1-t)\bar{y} - \bar{y}\}\right] \nonumber\\
	&& \;\; \times \exp\left[-\frac{n}{2} \tr\{\hat{\Lambda}_0 (\hat{\Lambda}_0^{-1}- t^2A/n ) \}\right] |\hat{\Lambda}^{-1}(t)|^\frac{n-p-g-1}{2}\nonumber\\
	&\propto&  \exp\left\{-\frac{1}{2}t^2 \sum_{i=1}^g(\bar{y}_i -\bar{y} )^\top\hat{\Lambda}_0(\bar{y}_i -\bar{y})\right\} \exp\left\{\frac{1}{2} t^2\tr( \hat{\Lambda}_0 A ) \right\} |\hat{\Lambda}^{-1}(t)|^\frac{n-p-g-1}{2}\nonumber\\
	&\propto& |\hat{\Lambda}^{-1}(t)|^\frac{n-p-g-1}{2}.\nonumber
\end{eqnarray}

Since the saddlepoint approximation $h\{s(t);\psi\}$ is in fact exact, up to a normalizing  constant,  the integral in the denominator of the directional $p$-value (7) of the main text   is  the normalizing constant of the conditional distribution of $||s||$ given the direction $s/||s||$.  The directional $p$-value is then the exact probability of $||s|| > ||s^0||$ given the direction $s/||s||$ under the null hypothesis $H_\psi$, and hence is  exactly uniformly distributed.

\subsection{Proof of Proposition 1 }	\label{section:proof propo1}

	First, we express the estimate $\hat{\Lambda}^{-1}(t) = \hat{\Lambda}_0^{-1} - t^2A/n$, where $n=n_1+n_2$, as
	\begin{align}
		\hat{\Lambda}^{-1}(t) =&  \hat{\Lambda}_0^{-1} -  t^2 \sum_{i=1}^{2} n_i (\bar{y}_i - \bar{y}) (\bar{y}_i - \bar{y})^\top/n,\nonumber
	\end{align}
	where $\bar{y} = (n_1\bar{y}_1 + n_2\bar{y}_2)/n$. Moreover, we have
	\begin{align}
		\sum_{i=1}^{2}  n_i (\bar{y}_i - \bar{y}) (\bar{y}_i - \bar{y})^\top  
		=& \; n_1 (\bar{y}_1 - \bar{y}) (\bar{y}_1 - \bar{y})^\top +  n_2 (\bar{y}_2 - \bar{y}) (\bar{y}_2 - \bar{y})^\top \nonumber \\
		=& \; \frac{n_1n_2(\bar{y}_1-\bar{y}_2)(\bar{y}_1-\bar{y}_2)^\top }{n}, \nonumber 
	\end{align}
	since 
	\begin{align}
		n_1\bar{y}_1\bar{y}_1^\top + n_2\bar{y}_2\bar{y}_2^\top  =& \frac{n_1^2\bar{y}_1\bar{y}_1^\top + n_1n_2\bar{y}_1\bar{y}_1^\top + n_1n_2\bar{y}_2\bar{y}_2^\top + n_2^2\bar{y}_2\bar{y}_2^\top}{n}, \nonumber\\
		n\bar{y}\bar{y}^\top =& \frac{n_1^2\bar{y}_1\bar{y}_1^\top + n_1n_2\bar{y}_1\bar{y}_2^\top + n_1n_2\bar{y}_2\bar{y}_1^\top + n_2^2\bar{y}_2\bar{y}_2^\top}{n}. \nonumber 
	\end{align}
	Then 
	\begin{align}
		\hat{\Lambda}^{-1}(t) =&  \hat{\Lambda}_0^{-1} -  t^2 \frac{n_1n_2}{n^2} (\bar{y}_1-\bar{y}_2)(\bar{y}_1-\bar{y}_2)^\top, \nonumber \\
		|\hat{\Lambda}^{-1}(t)|  =&  |\hat{\Lambda}_0^{-1}| \left\{1 -  t^2 \frac{n_1n_2}{n^2} (\bar{y}_1-\bar{y}_2)^\top\hat{\Lambda}_0 (\bar{y}_1-\bar{y}_2) \right\}.\nonumber 
	\end{align}
	The integrand function in the directional $p$-value (\ref{directed p-value}) along the line $s(t)$ can then be simplified to
	\begin{align}
		g(t) = t^{d-1} |\hat{\Lambda}^{-1}(t)|^{\frac{n-p-3}{2}} \propto t^{d-1} \left\{1 -  t^2 \frac{n_1n_2}{n^2} (\bar{y}_1-\bar{y}_2)^\top\hat{\Lambda}_0 (\bar{y}_1-\bar{y}_2)\right\}^\frac{n-p-3}{2}.\nonumber
	\end{align}
	Thus the directional $p$-value under the null hypothesis $H_\psi$ takes the form 
	\begin{align}
		p(\psi) = \frac{\int_{1}^{t_{\text{sup}}} t^{d-1} \left\{1 -  t^2 \frac{n_1n_2}{n^2} (\bar{y}_1-\bar{y}_2)^\top\hat{\Lambda}_0 (\bar{y}_1-\bar{y}_2)\right\}^\frac{n-p-3}{2} \text{d}t }{\int_{0}^{t_{\text{sup}}} t^{d-1} \left\{1 -  t^2 \frac{n_1n_2}{n^2} (\bar{y}_1-\bar{y}_2)^\top\hat{\Lambda}_0 (\bar{y}_1-\bar{y}_2)\right\}^\frac{n-p-3}{2} \text{d}t},\nonumber
	\end{align}
	where  $t_{\text{sup}} = \left\{\frac{n_1n_2 }{n^2}(\bar{y}_1-\bar{y}_2)^\top\hat{\Lambda}_0 (\bar{y}_1-\bar{y}_2) \right\}^{-1/2}$.
	In order to make the notation more compact, we define $C=\frac{n_1n_2 }{n^2}(\bar{y}_1-\bar{y}_2)^\top\hat{\Lambda}_0 (\bar{y}_1-\bar{y}_2)$, so that $t_\text{sup} = C^{-1/2}$. We can now rewrite the directional $p$-value as 
	\begin{align}
		p(\psi) = \frac{\int_{1}^{C^{-1/2}} t^{d-1} \left\{1 -  t^2 C \right\}^\frac{n-p-3}{2} \text{d}t }{\int_{0}^{C^{-1/2}} t^{d-1} \left\{1 -  t^2 C\right\}^\frac{n-p-3}{2} \text{d}t}. 
		\label{proof:directed pvalue}
	\end{align}
	Since the Hotelling's $T^2$ statistic has Hotelling's $T^2$ distribution, i.e., $T^2 \sim T^2(p, n-p-1)$ with  degrees of freedom $df_1 = p$ and $df_2 = n-p-1$, we change variable from $t$ to $\left\{C \left(x{df_2}/{df_1} + 1\right)\right\}^{-1/2}$. 
	
	The following steps are used to compute   the numerator in the directional $p$-value (\ref{proof:directed pvalue}).
	
	\noindent
	\textit{Step 1}. Change of the integration  interval:
	\begin{align}
		& 1 \le t \le C^{-1/2} \Leftrightarrow   1 \le \left\{C \left(\frac{df_2}{df_1}x + 1\right)\right\}^{-1/2} \le C^{-1/2} \nonumber\\
		& \Leftrightarrow   C \le  \left(\frac{df_2}{df_1}x + 1\right)^{-1} \le 1   \Leftrightarrow   1 \le  \frac{df_2}{df_1}x + 1 \le C^{-1}\nonumber \\
		& \Leftrightarrow   0 \le  \frac{df_2}{df_1}x  \le C^{-1} -1 \Leftrightarrow   0 \le  x  \le \frac{df_1}{df_2} \frac{1-C}{C}\nonumber 
	\end{align}
	hence the new integral is on $\left[0, \frac{df_1}{df_2} \frac{1-C}{C} \right]$. 
	
	\noindent
	\textit{Step 2}. Change of variable from $t$ to $x$:
	\begin{align}
		&\int_{1}^{C^{-1/2}} t^{d-1} \left\{1 -  t^2 C \right\}^\frac{n-p-3}{2} \text{d}t \nonumber\\
		=& \int_{0}^{\frac{df_1}{df_2} \frac{1-C}{C}} \left\{C \left(\frac{df_2}{df_1}x + 1\right)\right\}^{-\frac{d-1}{2}}  \left\{1 -  C^{-1} \left(\frac{df_2}{df_1}x + 1\right)^{-1} C \right\}^\frac{n-p-3}{2} \text{d}\left\{C \left(\frac{df_2}{df_1}x + 1\right)\right\}^{-\frac{1}{2}},\nonumber
	\end{align}
	where $d= p(k-1) = p$. Then we get
	\begin{align}
		{\text{d}\left\{C \left(\frac{df_2}{df_1}x + 1\right)\right\}^{-\frac{1}{2}}} = C^{-1/2} (-\frac{1}{2})  \left(\frac{df_2}{df_1}x + 1\right)^{-\frac{3}{2}} \frac{df_2}{df_1}{\text{d}x},\nonumber
	\end{align}
	and 
	\begin{align}
		&\int_{1}^{C^{-1/2}} t^{d-1} \left\{1 -  t^2 C \right\}^\frac{n-p-3}{2} \text{d}t \nonumber\\
		=& c \int_{0}^{\frac{df_1}{df_2} \frac{1-C}{C}} \left(\frac{df_2}{df_1}x + 1\right)^{-\frac{p-1}{2}}  \left\{1 - \left(\frac{df_2}{df_1}x + 1\right)^{-1}  \right\}^\frac{n-p-3}{2}  \left(\frac{df_2}{df_1}x + 1\right)^{-\frac{3}{2}} \text{d}x \nonumber\\
		= & c \int_{0}^{\frac{df_1}{df_2} \frac{1-C}{C}} \left(\frac{df_2}{df_1}x + 1\right)^{-\frac{p-1+3}{2}}  \left\{  \frac{\frac{df_2}{df_1}x+1-1 }{\frac{df_2}{df_1}x+1} \right\}^\frac{n-p-3}{2}  \text{d}x \nonumber\\
		= & c \int_{0}^{\frac{df_1}{df_2} \frac{1-C}{C}} \left(\frac{df_2}{df_1}x + 1\right)^{-\frac{n-p-1+p}{2}} \left(\frac{df_2}{df_1} x \right)^{\frac{n-p-1}{2} -1}  \text{d}x \nonumber\\
		= & c \int_{0}^{\frac{df_1}{df_2} \frac{1-C}{C}} \left(\frac{df_2}{df_1}x + 1\right)^{-\frac{df_2+df_1}{2}} \left(\frac{df_2}{df_1} x \right)^{\frac{df_2}{2} -1}  \text{d}x \nonumber\\
		= & c'F\left(\frac{df_1}{df_2} \frac{1-C}{C}, df_2, df_1\right),\nonumber
	\end{align}
	where $F(X, df_2, df_1)$ is the cumulative distribution function of a  random variable $X$ following a  $F$-distribution with degrees of freedom $df_2$ and $df_1$. Since if  $X \sim F(df_2, df_1)$, then $X^{-1} \sim F(df_1, df_2)$, we can express the directional $p$-value  as
	\begin{align*}
		p(\psi) = 1 - F\left(\frac{df_2}{df_1} \frac{C}{1-C}, df_1, df_2\right).
	\end{align*}
	According to the Sherman-Morrison formula \citep[see also][Section 4.2]{sartori:2019Ftest}, we have that \\ $C =  \left\{\left. {n_1n_2/n} (\bar{y}_1 - \bar{y}_2)^\top \hat{\Lambda}(\bar{y}_1 - \bar{y}_2)  \right\}\right/ \left\{ 1 + {n_1n_2/{n}} (\bar{y}_1 - \bar{y}_2)^\top \hat{\Lambda}(\bar{y}_1 - \bar{y}_2) \right\}$ and 
	\begin{align}
		\frac{df_2}{df_1} \frac{C}{1-C} =&  \frac{n-p-1}{p} \cdot \frac{n_1n_2}{n^2} (\bar{y}_1-\bar{y}_2)^\top\hat{\Lambda}(\bar{y}_1-\bar{y}_2)\nonumber\\
		=& \frac{n-p-1}{p} \cdot \frac{T^2}{n_1+n_2-2}, \nonumber
	\end{align}
	with $T^2 = \left.{n_1n_2}(\bar{y}_1 - \bar{y}_2)^\top S^{-1} (\bar{y}_1 - \bar{y}_2)\right/{n}$ and $S =\left.{n} \hat{\Lambda}^{-1}\right/{(n_1 - n_2- 2)}$. Since $\left.{(n-p-1)} {T^2}\right/ \{p(n_1+n_2-2)\}  \sim F(p,n-p-1)$, the directional test is identical to the Hotelling $T^2$ test.

\section{Simulation studies for homoscedastic one-way MANOVA}	\label{section simulation:dirmean}

This section reports additional empirical result for homoscedastic one-way MANOVA in the multivariate normal framework. We compare the performance of exact directional test (DT) with other five approximate approaches: the central limite theorem test (CLT), log-likelihood ratio test (LRT), Bartlett correction (BC) and two \citeauthor{skovgaard:2001}'s modifications (Sko1 and Sko2). The six tests are evaluated in terms of empirical size.  The simulation results are computed via Monte Carlo simulation based on 10,000 replications.

\begin{figure}[t]
	\centering
	\captionsetup{font=footnotesize}
	\subfigure{
		\begin{minipage}[b]{.9\linewidth}
			\centering
			\includegraphics[scale=0.4]{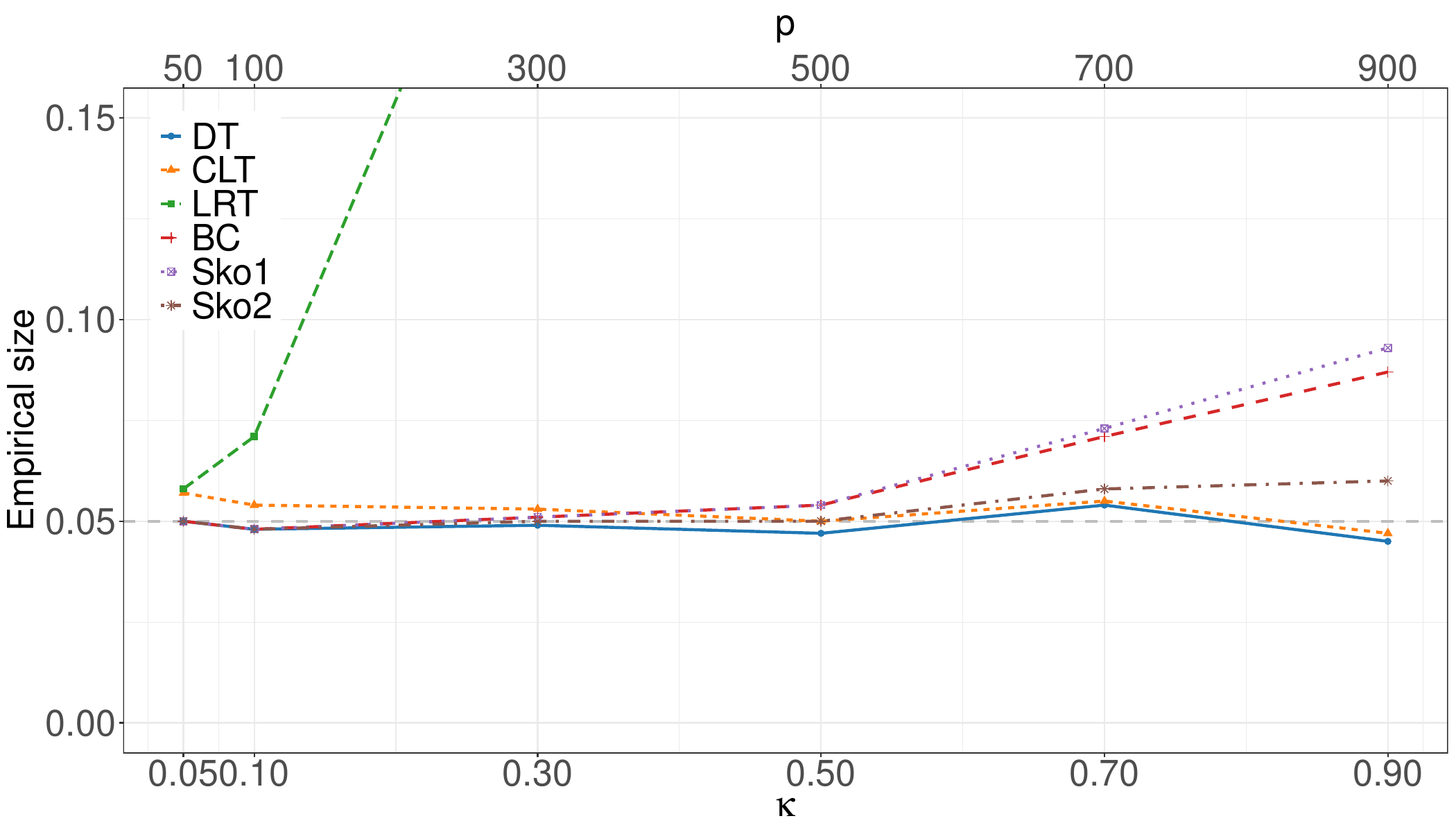}
		\end{minipage}
	}
	\caption{
		Empirical size of  the directional test (DT), central limit theorem test (CLT), likelihood ratio test (LRT), Bartlett corrected test (BC) and two \citeauthor{skovgaard:2001}'s modifications (Sko1 and Sko2) for hypothesis (8) with $g=3$, at nominal level $\alpha = 0.05$ given by the dashed gray horizontal line. The left, middle and right panels correspond to $n_i = 1000$, respectively. 
	}
	\label{fig:type I error same 1000}
\end{figure}

\subsection{Empirical results for the high dimentional setup}
Groups of size $n_i, i \in \{1,\dots, g\}$,  are  generated from a $p$-variate standard normal distribution $N_p(0_p, I_p)$ under the null hypothesis.  For each simulation experiment, we show  results for $p = \kappa n_i$ with $\kappa \in \{0.05, 0.1,0.3,$ $0.5,0.7,0.9\}$  and $n_i = 1000$. Throughout, we set $g = 3$ and $n_1= n_2 =n_3$.  Figure \ref{fig:type I error same 1000} shows the empirical size of the directional test, central limit theorem test, likelihood ratio test, Bartlett corrected test and two \citeauthor{skovgaard:2001}'s modifications for hypothesis (8) with $g=3$, at nominal level $\alpha =0.05$.  The directional $p$-value performs better than the other.  The log-likelihood ratio test is disastrous in the large sample size setting with $n_i=1000$.


\setlength{\tabcolsep}{3.5mm}{
	\begin{table}[H]
		\centering
		\caption{Empirical size of  the directional test (DT), central limit theorem test (CLT), likelihood ratio test (LRT), Bartlett correction (BC) and two \citeauthor{skovgaard:2001}'s modifications  (Sko1 and Sko2) for hypothesis (8) in the paper with $n_i \in \{ 100, 500, 1000\}$, at nominal level $\alpha = 0.05$.}
		\medskip
		{\begin{tabular}{cccccccc}
				\toprule[0.09 em]
				$n_i$	&$\tau (p)$ &  DT & CLT & LRT &BC & Sko1 & Sko2\\
				\hline
				100& 0.250 (3) & 0.054 & 0.072 & 0.057 & 0.054 & 0.053 & 0.053 \\ 
				& 0.333 (4) & 0.047 & 0.063 & 0.051 & 0.047 & 0.046 & 0.046 \\ 
				& 0.417 (6) & 0.050 & 0.066 & 0.056 & 0.050 & 0.049 & 0.049 \\ 
				&0.500 (10) & 0.048 & 0.060 & 0.057 & 0.048 & 0.046 & 0.046 \\ 
				&0.580 (14) & 0.052 & 0.064 & 0.067 & 0.052 & 0.050 & 0.050 \\ 
				&0.667 (21) & 0.050 & 0.061 & 0.080 & 0.049 & 0.048 & 0.047 \\ 
				&0.750 (31) & 0.048 & 0.059 & 0.104 & 0.049 & 0.046 & 0.046 \\ 
				&0.833 (46) & 0.049 & 0.057 & 0.161 & 0.051 & 0.049 & 0.048 \\ 
				&0.917 (68) & 0.050 & 0.057 & 0.311 & 0.056 & 0.052 & 0.048 \\ 
				\hline
				500 & 0.250 (4) & 0.051 & 0.068 & 0.052 & 0.051 & 0.051 & 0.051 \\ 
				&0.333 (7) & 0.053 & 0.070 & 0.055 & 0.053 & 0.053 & 0.053 \\ 
				&0.417 (13) & 0.051 & 0.063 & 0.054 & 0.051 & 0.051 & 0.051 \\ 
				&0.500 (22) & 0.050 & 0.061 & 0.056 & 0.050 & 0.050 & 0.050 \\ 
				&0.580 (37) & 0.048 & 0.056 & 0.058 & 0.048 & 0.048 & 0.048 \\ 
				&0.667 (62) & 0.055 & 0.063 & 0.078 & 0.055 & 0.054 & 0.054 \\ 
				&0.750 (105) & 0.050 & 0.054 & 0.105 & 0.051 & 0.050 & 0.049 \\ 
				&0.833 (177) & 0.054 & 0.058 & 0.222 & 0.056 & 0.055 & 0.054 \\ 
				&0.917 (297) & 0.047 & 0.050 & 0.609 & 0.057 & 0.056 & 0.049 \\ 
				\hline
				1000	&0.250 (5) & 0.054 & 0.070 & 0.055 & 0.054 & 0.054 & 0.054 \\ 
				&0.333  (9) & 0.052 & 0.066 & 0.053 & 0.053 & 0.052 & 0.052 \\ 
				&0.417 (17) & 0.049 & 0.061 & 0.050 & 0.049 & 0.049 & 0.049 \\ 
				&0.500 (31) & 0.054 & 0.063 & 0.058 & 0.054 & 0.054 & 0.054 \\ 
				&0.580 (56) & 0.052 & 0.059 & 0.061 & 0.052 & 0.052 & 0.052 \\ 
				&0.667 (99) & 0.049 & 0.055 & 0.072 & 0.049 & 0.049 & 0.049 \\ 
				&0.750 (177) & 0.053 & 0.056 & 0.116 & 0.053 & 0.053 & 0.052 \\ 
				&0.833 (316) & 0.048 & 0.050 & 0.261 & 0.049 & 0.049 & 0.048 \\ 
				&0.917 (562) & 0.053 & 0.055 & 0.794 & 0.064 & 0.065 & 0.056 \\ 
				\bottomrule[0.09 em]            
			\end{tabular}
		}
		\label{table type I normal same he}
\end{table}}

\subsection{Empirical results for moderate dimentional setup}
Groups of size $n_i, i \in \{1,\dots, g\}$,  are  generated from a $p$-variate standard normal distribution $N_p(0_p, I_p)$ under the null hypothesis.  For each simulation experiment, we show  results for $ p = \lfloor n_i ^ \tau \rfloor$ with $\tau = j/24$, $j \in \{6,8,\cdots,22\}$ and $n_i \in \{100, 500, 1000\}$. Throughout, we set $g = 3$ and $n_1= n_2 =n_3$.

\begin{figure}[t]
	\centering
	\captionsetup{font=footnotesize}
	\subfigure{
		\begin{minipage}[b]{.3\linewidth}
			\centering
			\includegraphics[scale=0.2]{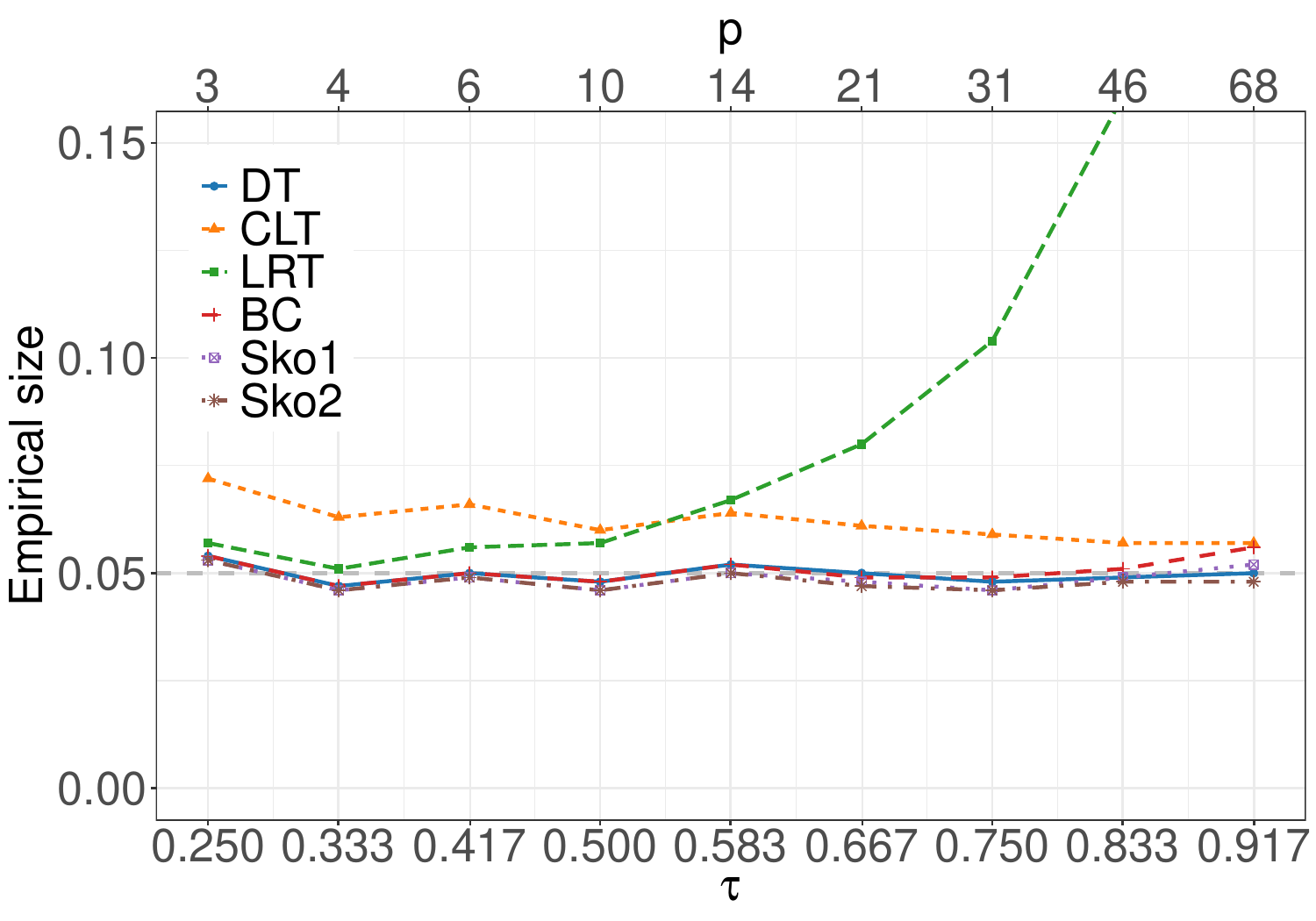}
		\end{minipage}
		\begin{minipage}[b]{.3\linewidth}
			\centering
			\includegraphics[scale=0.2]{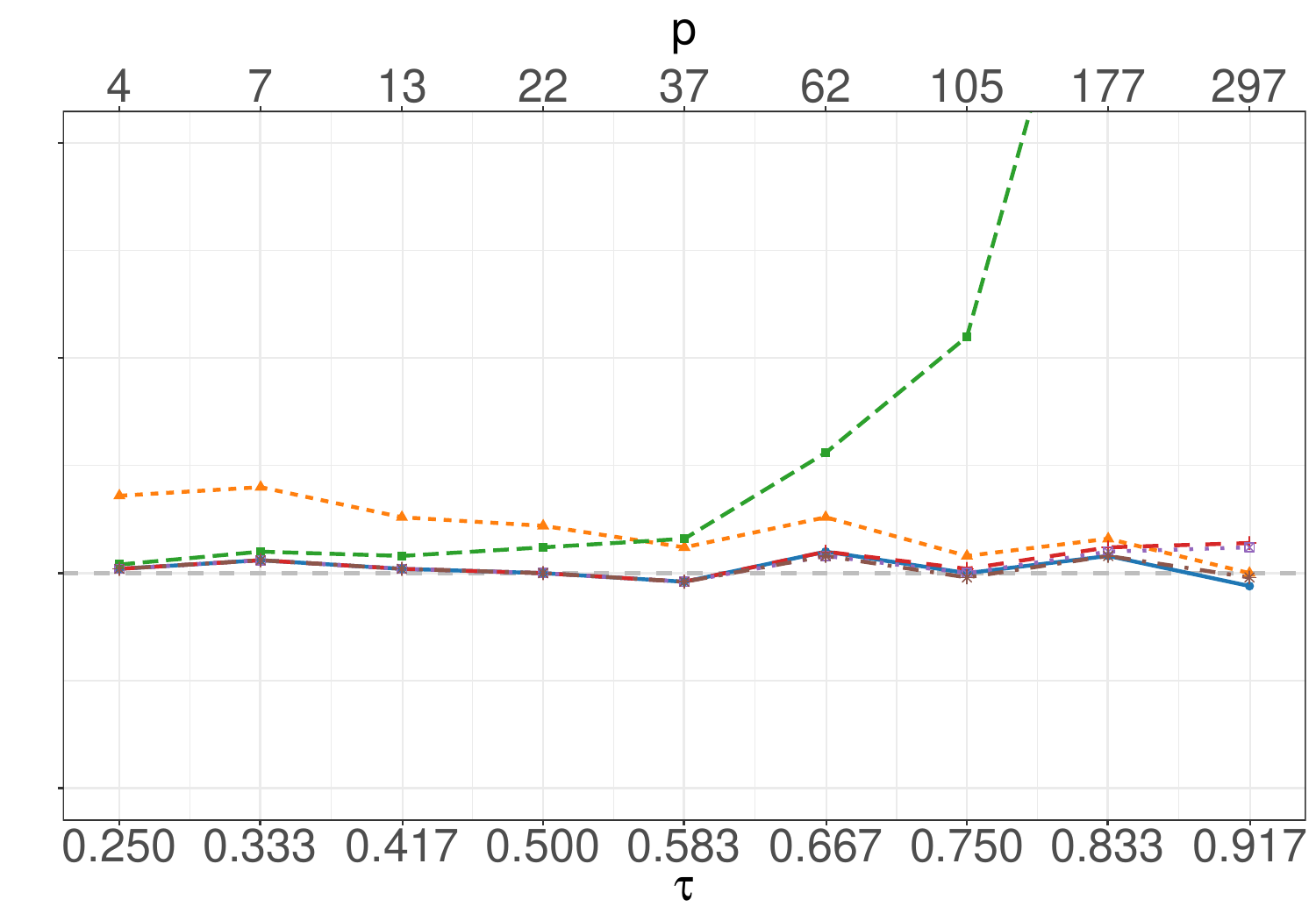}
		\end{minipage}
		\begin{minipage}[b]{.3\linewidth}
			\centering
			\includegraphics[scale=0.2]{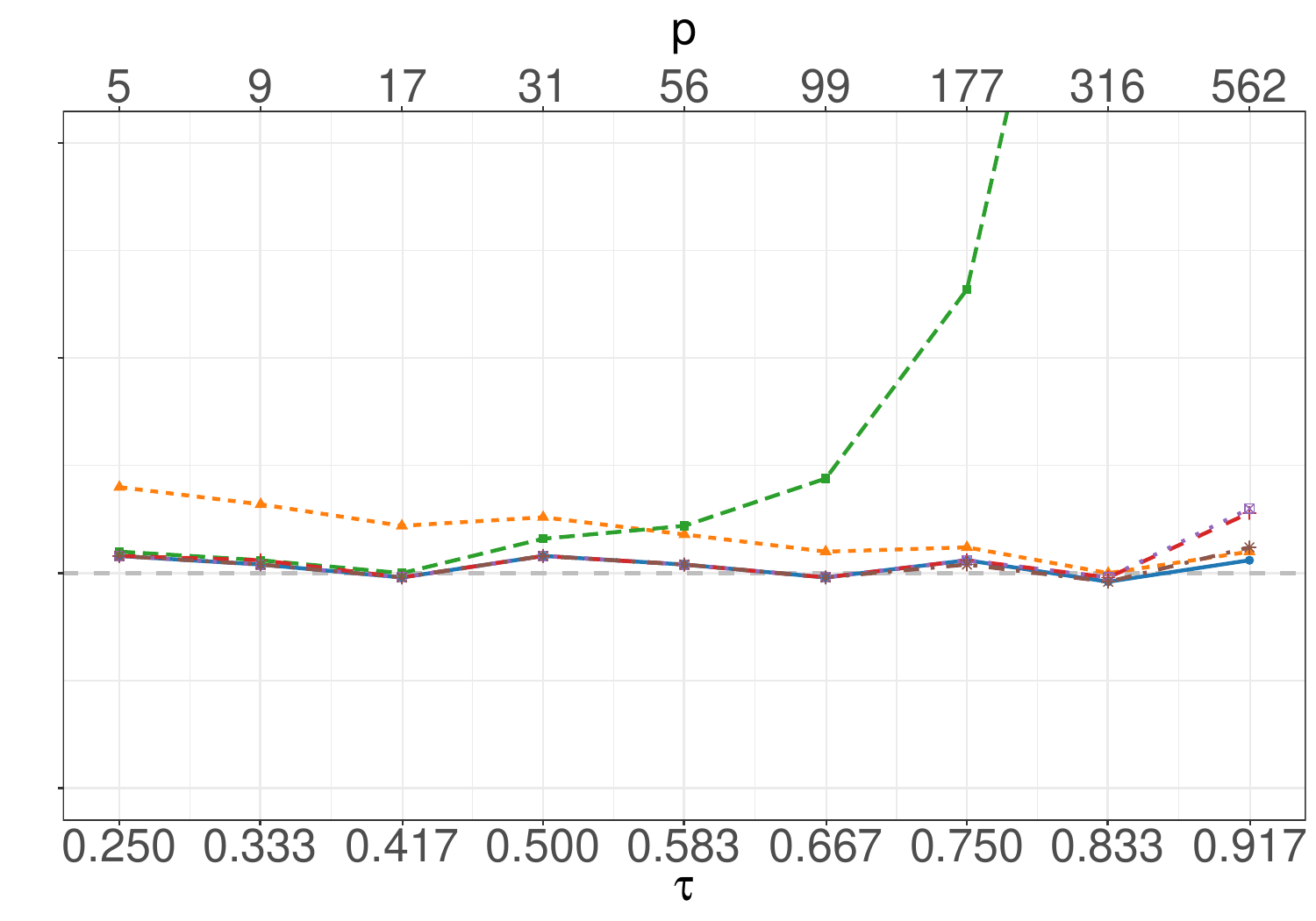}
		\end{minipage}
	}
	\caption{
		Empirical size of  the directional test (DT), central limit theorem test (CLT), likelihood ratio test (LRT), Bartlett correction (BC) and two \citeauthor{skovgaard:2001}'s modifications (Sko1 and Sko2) for hypothesis (8) in the paper, at nominal level $\alpha = 0.05$ given by the gray horizontal line. The left, middle and right panels correspond to $n_i = 100, 500, 1000$, respectively ($g=3$). 
	}
	\label{fig:type I error same}
\end{figure}

\subsection{Empirical results for a  large number $g$ of groups}

\setlength{\tabcolsep}{3.5mm}{
	\begin{table}[H]
		\centering
		\caption{Empirical size of the directional test (DT), central limit theorem test (CLT), likelihood ratio test (LRT), Bartlett correction (BC) and two \citeauthor{skovgaard:2001}'s modifications (Sko1 and Sko2) for hypothesis (8) in the paper with $p = n_i^\tau$, $n_i = 100$ and $g =30$, at nominal level $\alpha = 0.05$}
		\medskip
		{\begin{tabular}{ccccccc}
				\toprule[0.09 em]
				$\tau (p)$ &  DT & CLT & LRT &BC & Sko1 & Sko2\\
				\hline
				0.250 (3) & 0.053 & 0.060 & 0.057 & 0.053 & 0.053 & 0.053 \\ 
				0.333 (4) & 0.050 & 0.057 & 0.055 & 0.050 & 0.050 & 0.050 \\ 
				0.417 (6) & 0.050 & 0.056 & 0.057 & 0.050 & 0.050 & 0.050 \\ 
				0.500 (10) & 0.049 & 0.053 & 0.059 & 0.049 & 0.049 & 0.049 \\ 
				0.580 (14) & 0.051 & 0.056 & 0.064 & 0.052 & 0.051 & 0.051 \\ 
				0.667 (21) & 0.050 & 0.054 & 0.071 & 0.050 & 0.050 & 0.050 \\ 
				0.750 (31) & 0.052 & 0.054 & 0.080 & 0.052 & 0.052 & 0.052 \\ 
				0.833 (46) & 0.048 & 0.051 & 0.095 & 0.049 & 0.048 & 0.048 \\ 
				0.917 (68) & 0.054 & 0.056 & 0.137 & 0.054 & 0.054 & 0.053 \\ 
				\bottomrule[0.09 em]            
			\end{tabular}
		}
		\label{table type I normal same he100 k30}
\end{table}}

\begin{figure}[H]
	\centering
	\captionsetup{font=footnotesize}
	\includegraphics[scale=0.4]{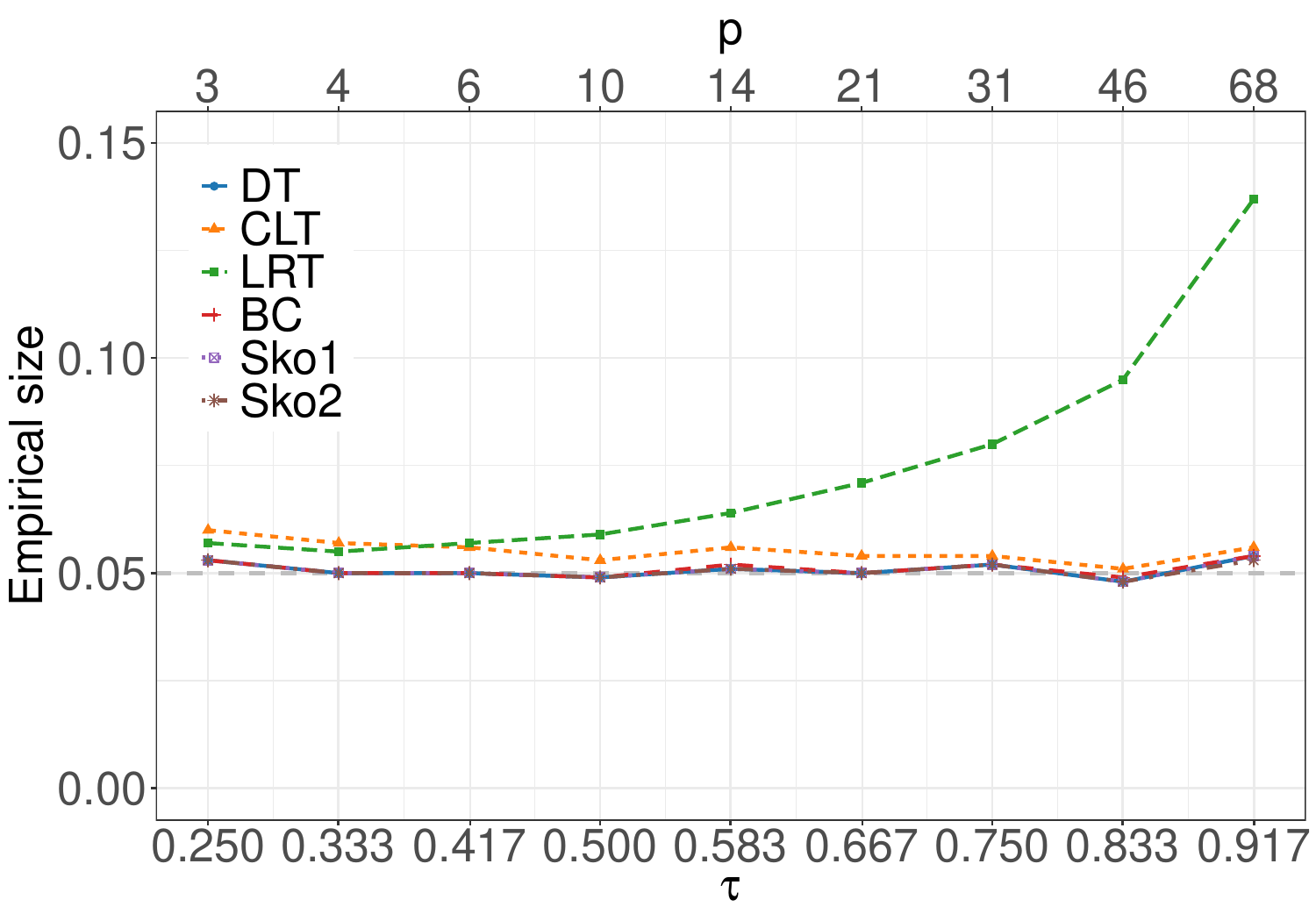}
	\caption{ 	Empirical size of the directional test (DT), central limit theorem test (CLT), likelihood ratio test (LRT), Bartlett correction (BC) and two \citeauthor{skovgaard:2001}'s modifications  (Sko1 and Sko2) for hypothesis (8) in the paper with $n_i = 100$ and $g=30$, at nominal level $\alpha = 0.05$  given by the gray horizontal line.}
	\label{figure:type I error same he100 k30}
\end{figure}

\subsection{Empirical results for the setup of  \citet[][Section A.3]{He2020B}}

In each Monte Carlo experiment, we show  results for $ p = \lfloor n ^ \tau \rfloor$ with $n  = \sum_{i=1}^g n_i$. Under the null hypothesis, we set  $g = 3$, $n_1= n_2 =n_3$, $\tau = j/24$ with $j \in \{6,8,\cdots,22\}$ and $n_i=100$.

\setlength{\tabcolsep}{3.5mm}{
	\begin{table}[H]
		\centering
		\caption{Empirical size of the directional test (DT), central limit theorem test (CLT), likelihood ratio test (LRT), Bartlett correction (BC) and two \citeauthor{skovgaard:2001}'s modifications  (Sko1 and Sko2) for hypothesis (8) in the paper with $p = n^\tau$ and $n = \sum_{i=1}^3 100 = 300$, at  nominal level $\alpha = 0.05$.}
		\medskip
		{\begin{tabular}{ccccccc}
				\toprule[0.09 em]
				$\tau$ ($p$) &  DT & CLT & LRT &BC & Sko1 & Sko2\\
				\hline
				0.250 (4) & 0.048 & 0.064 & 0.052 & 0.047 & 0.046 & 0.046 \\ 
				0.333 (6) & 0.049 & 0.066 & 0.056 & 0.049 & 0.047 & 0.047 \\ 
				0.417 (10) & 0.050 & 0.064 & 0.061 & 0.050 & 0.049 & 0.049 \\ 
				0.500 (17) & 0.053 & 0.065 & 0.074 & 0.053 & 0.051 & 0.051 \\ 
				0.583 (27) & 0.046 & 0.056 & 0.088 & 0.047 & 0.045 & 0.045 \\ 
				0.667 (44) & 0.051 & 0.060 & 0.157 & 0.053 & 0.050 & 0.049 \\ 
				0.750 (72) & 0.048 & 0.054 & 0.347 & 0.054 & 0.052 & 0.047 \\ 
				0.833 (115) & 0.050 & 0.056 & 0.832 & 0.078 & 0.077 & 0.057 \\ 
				0.917 (186) & 0.044 & 0.049 & 1.000 & 0.280 & 0.278 & 0.106 \\ 
				\bottomrule[0.09 em]            
			\end{tabular}
		}
		\label{table type I normal same he100 sum}
\end{table}}

\begin{figure}[H]
	\centering
	\captionsetup{font=footnotesize}
	\includegraphics[scale=0.4]{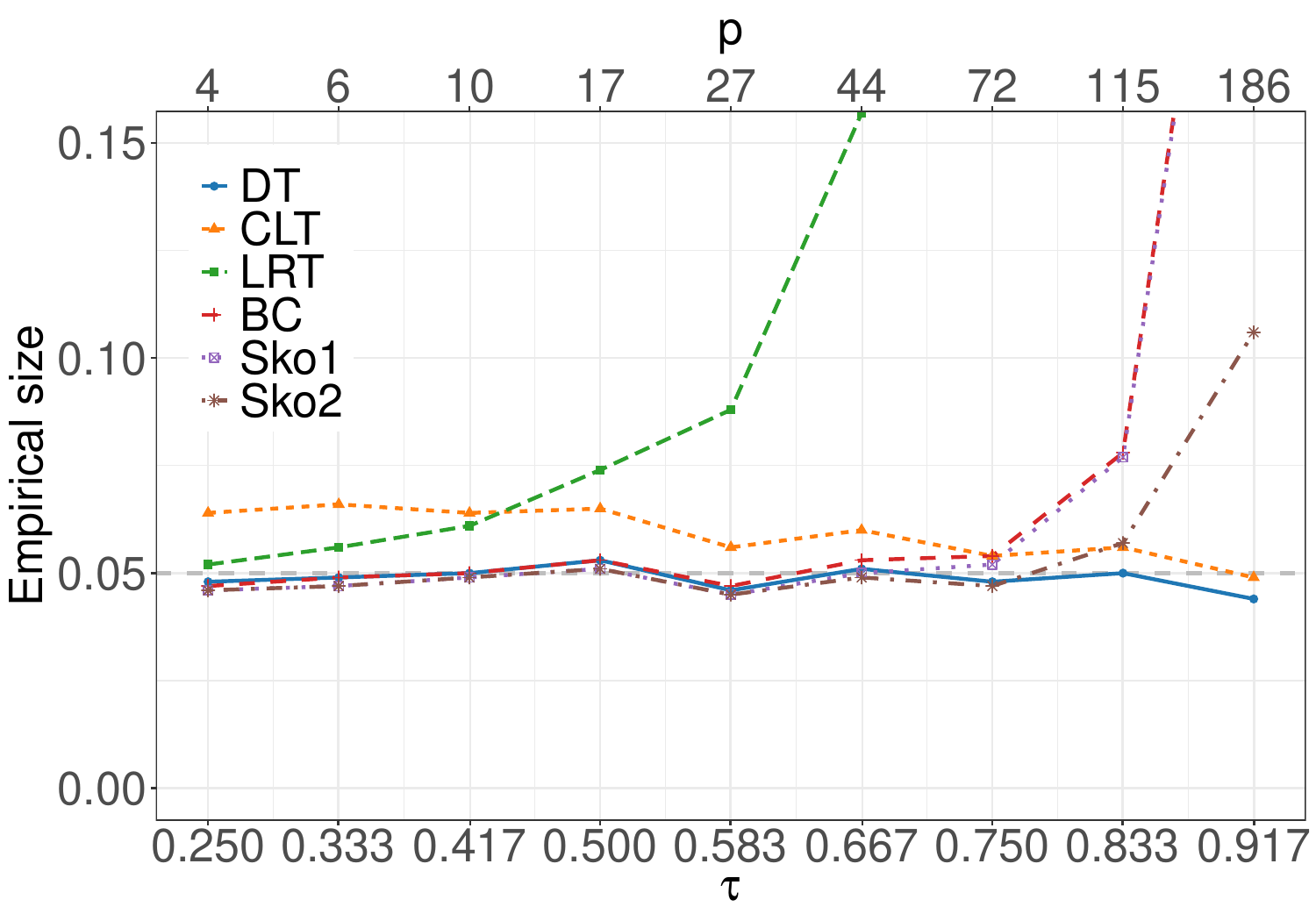}
	\caption{ Empirical size of the directional test (DT), central limit theorem test (CLT), likelihood ratio test (LRT), Bartlett correction (BC) and two \citeauthor{skovgaard:2001}'s modifications  (Sko1 and Sko2) for hypothesis (8) in the paper with $p = n^\tau$ and $n = \sum_{i=1}^3 100 = 300$, at  nominal level $\alpha = 0.05$  given by the gray horizontal line.}
	\label{figure:type I error same he100 sum}
\end{figure}

\setlength{\tabcolsep}{3.5mm}{
	\begin{table}[H]
		\centering
		\caption{Empirical size of  the directional test (DT), central limit theorem test (CLT), likelihood ratio test (LRT), Bartlett correction (BC) and two \citeauthor{skovgaard:2001}'s modifications (Sko1 and Sko2) for hypothesis (8) in the paper with $p = n^\tau$ and $n = \sum_{i=1}^3 500 = 1500$, at nominal level $\alpha = 0.05$.}
		\medskip
		{\begin{tabular}{ccccccc}
				\toprule[0.09 em]
				$\tau$ ($p$) &  DT & CLT & LRT &BC & Sko1 & Sko2\\
				\hline
				0.250 (6) & 0.052 & 0.068 & 0.053 & 0.052 & 0.052 & 0.052 \\ 
				0.333 (11) & 0.051 & 0.063 & 0.053 & 0.051 & 0.050 & 0.050 \\ 
				0.417 (21) & 0.050 & 0.059 & 0.053 & 0.050 & 0.049 & 0.049 \\ 
				0.500 (38) & 0.051 & 0.060 & 0.062 & 0.051 & 0.051 & 0.051 \\ 
				0.583 (71) & 0.051 & 0.059 & 0.080 & 0.051 & 0.050 & 0.050 \\ 
				0.667 (131) & 0.051 & 0.056 & 0.146 & 0.052 & 0.051 & 0.051 \\ 
				0.750 (241) & 0.052 & 0.057 & 0.402 & 0.058 & 0.058 & 0.053 \\ 
				0.833 (443) & 0.052 & 0.056 & 0.976 & 0.081 & 0.083 & 0.063 \\ 
				0.917 (815) & 0.049 & 0.051 & 1.000 & 0.420 & 0.462 & 0.171 \\ 
				\bottomrule[0.09 em]            
			\end{tabular}
		}
		\label{table type I normal same he500 sum}
\end{table}}

\begin{figure}[H]
	\centering
	\captionsetup{font=footnotesize}
	\includegraphics[scale=0.4]{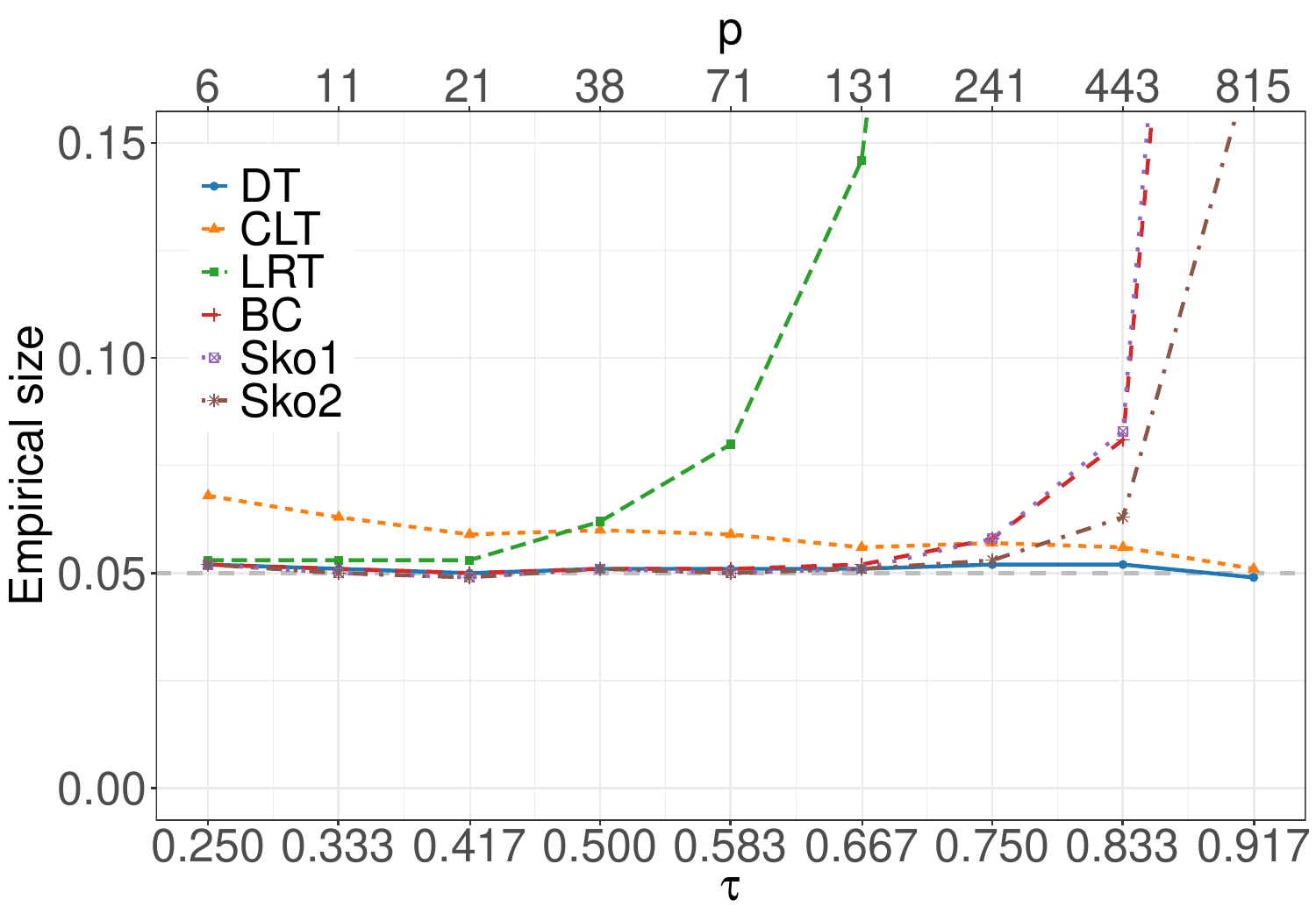}
	\caption{ Empirical size of  the directional test (DT), central limit theorem test (CLT), likelihood ratio test (LRT), Bartlett correction (BC) and two \citeauthor{skovgaard:2001}'s modifications (Sko1 and Sko2) for hypothesis (8) in the paper with $p = n^\tau$ and $n = \sum_{i=1}^3 500 = 1500$, at nominal level $\alpha = 0.05$  given by the gray horizontal line.}
	\label{figure:type I error same he500 sum}
\end{figure}

\subsection{Robustness to misspecification for the high dimensional setup}\label{section:robustness high}

In general, all the approaches examined so far rely on the  normal model and are not guaranteed to be robust under model misspecification. We can assess 
numerically the robustness of the various competitors using simulations. We consider three different distributions for   the true generating process:   multivariate $t$, multivariate skew-normal \citep{azzalini1999statistical} or multivariate Laplace. In more detail,    a multivariate $t$ distribution  with location $0_p$, scale matrix $I_p$ and degrees of freedom $5$, a multivariate skew-normal distribution  with location $1_p$, scale matrix $\Omega = (\omega_{jl})= (0.2)^{|j-l|}$ and shape parameter $1_p$, and a multivariate Laplace distribution with mean vector $1_p$ and identity covariance matrix.
Simulation results are based on $10,000$ replications. 

\begin{figure}[th]
	\centering
	\captionsetup{font=footnotesize}
	\subfigure{
		\begin{minipage}[b]{.32\linewidth}
			\centering
			\includegraphics[scale=0.2]{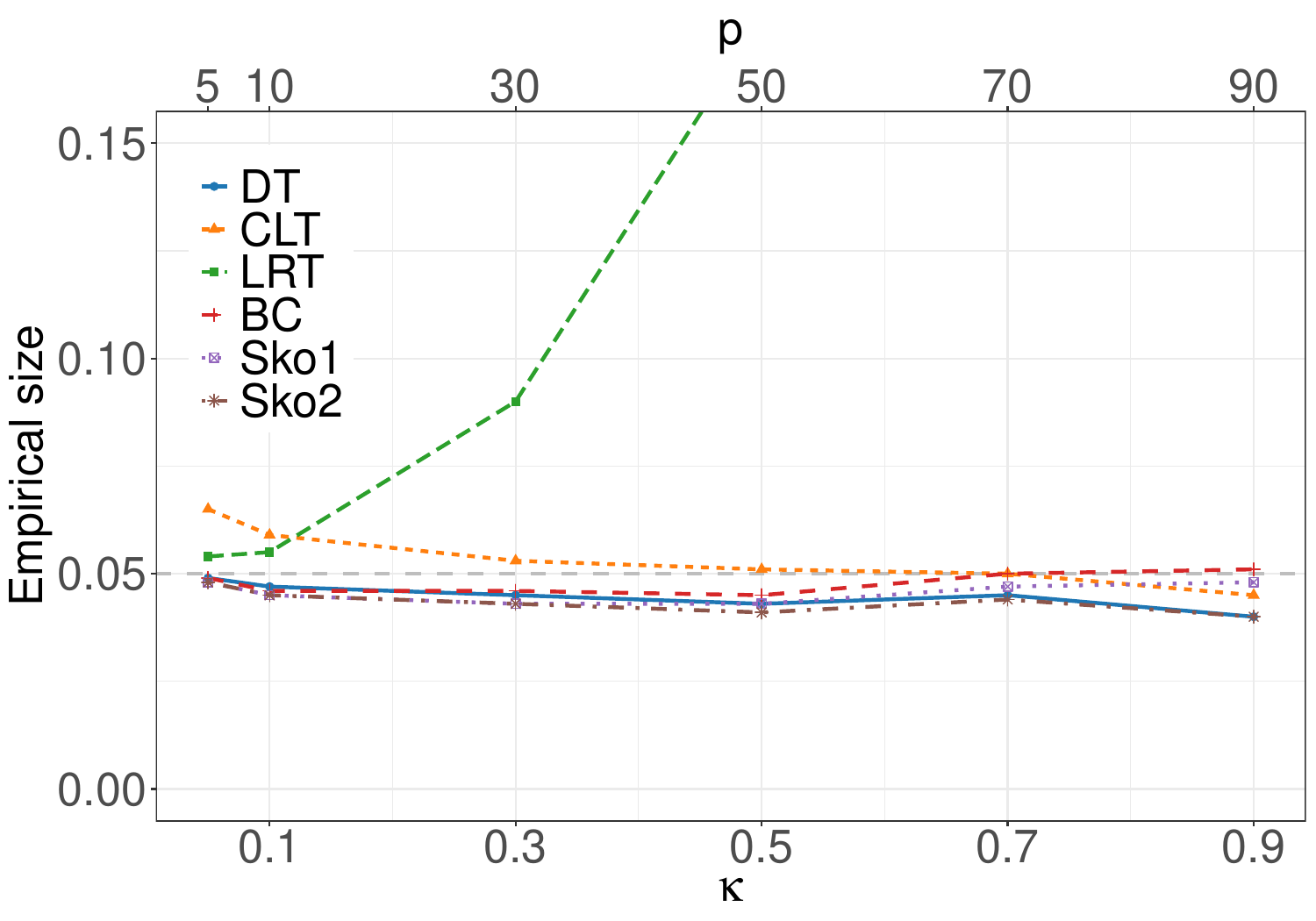}
		\end{minipage}
		\begin{minipage}[b]{.33\linewidth}
			\centering
			\includegraphics[scale=0.2]{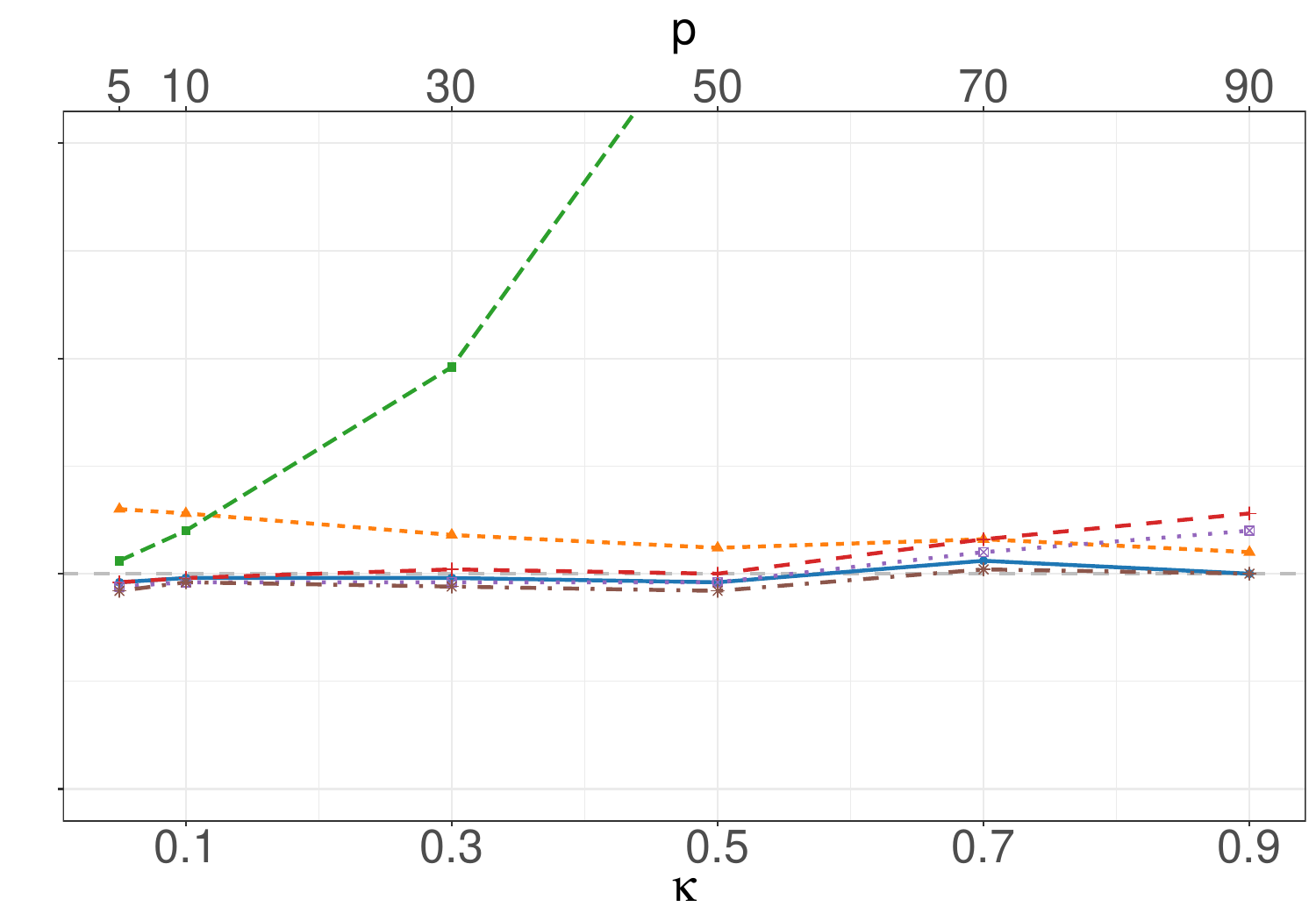}
		\end{minipage}
		\begin{minipage}[b]{.33\linewidth}
			\centering
			\includegraphics[scale=0.2]{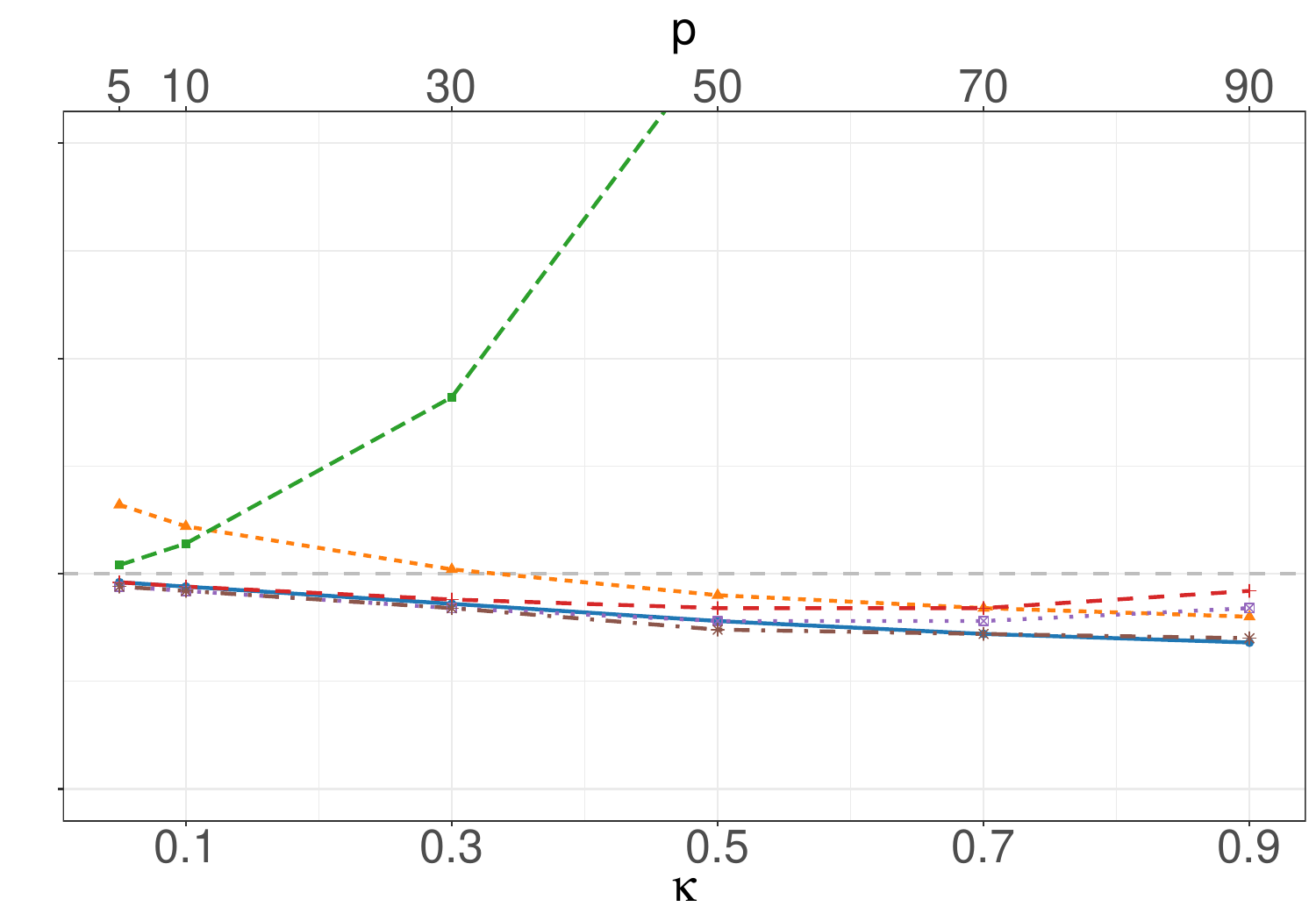}
		\end{minipage}
	}
	\caption{
		Empirical size for the directional test (DT), central limit theorem test (CLT), likelihood ratio test (LRT), Bartlett corrected test (BC) and two \citeauthor{skovgaard:2001}'s modifications (Sko1 and Sko2) for hypothesis (\ref{hypothesis:dirmean}) with $g=3$, at nominal level $\alpha = 0.05$  given by the dashed gray horizontal line. The left, middle and right panels correspond to multivariate $t$, multivariate skew-normal, and multivariate Laplace distributions of the true generating process, respectively, with $n_i = 100$.
	}
	\label{fig:type I error same robustness 100 high}
\end{figure}
\begin{figure}[th]
	\centering
	\captionsetup{font=footnotesize}
	\subfigure{
		\begin{minipage}[b]{.32\linewidth}
			\centering
			\includegraphics[scale=0.2]{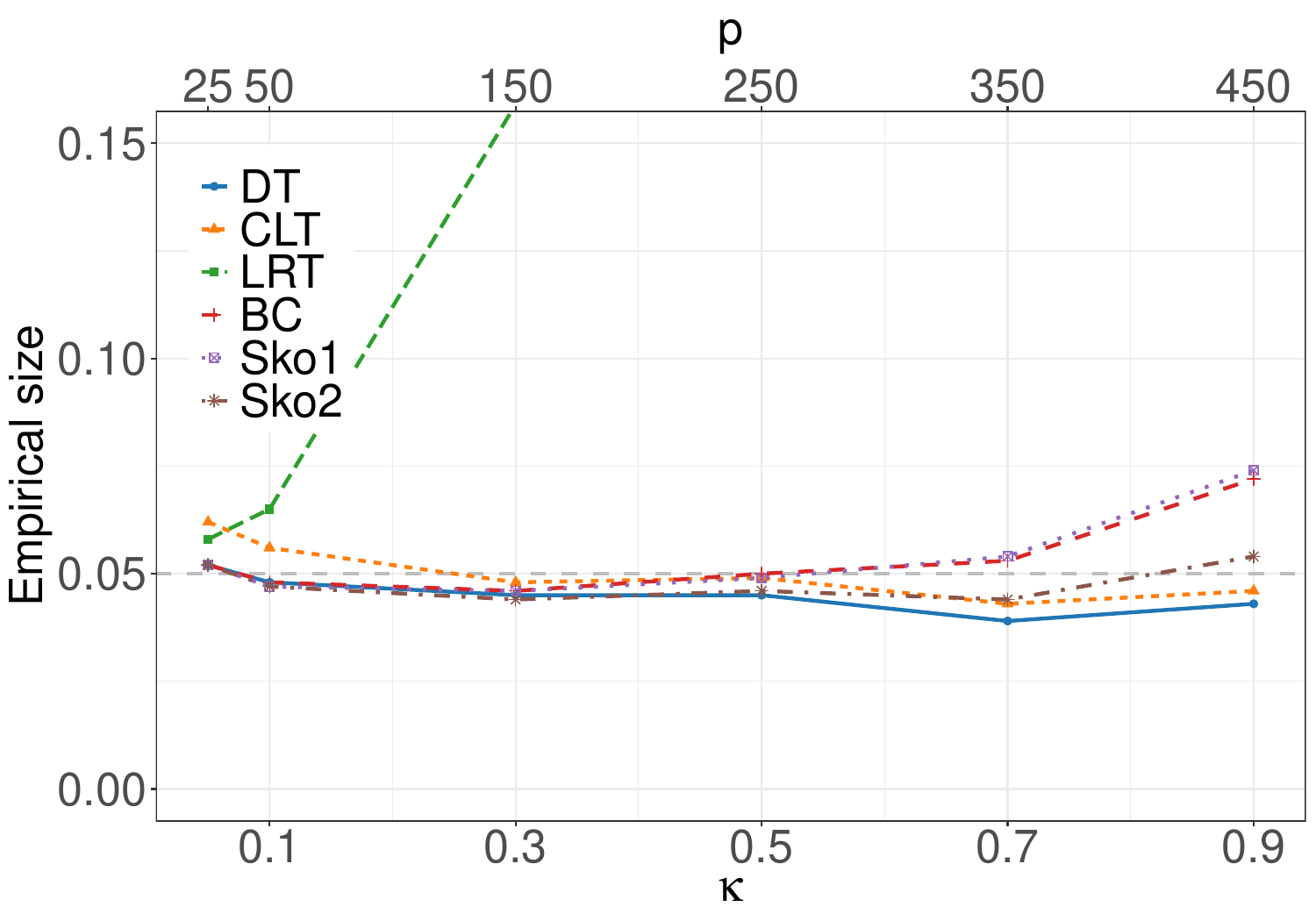}
		\end{minipage}
		\begin{minipage}[b]{.33\linewidth}
			\centering
			\includegraphics[scale=0.2]{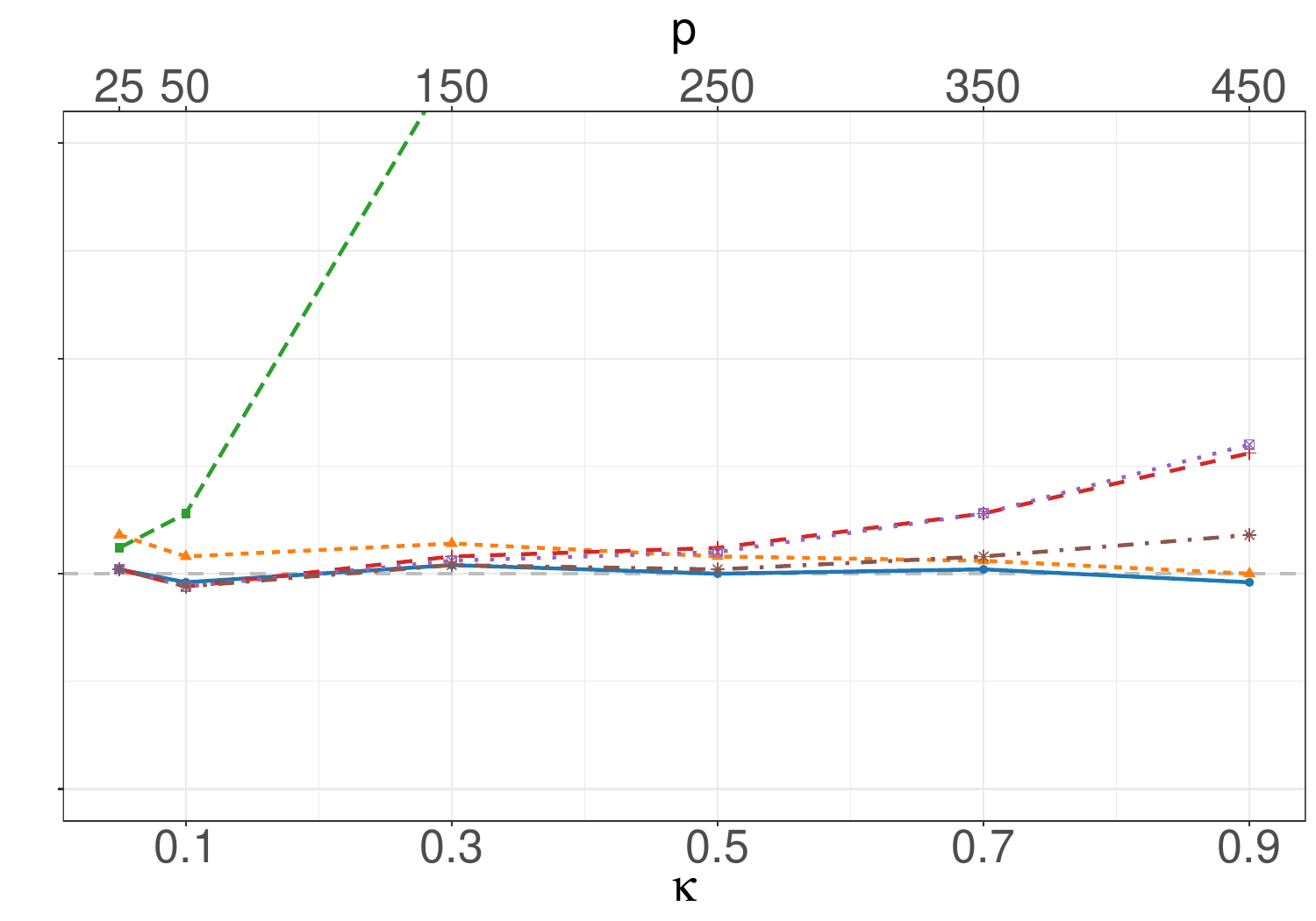}
		\end{minipage}
		\begin{minipage}[b]{.33\linewidth}
			\centering
			\includegraphics[scale=0.2]{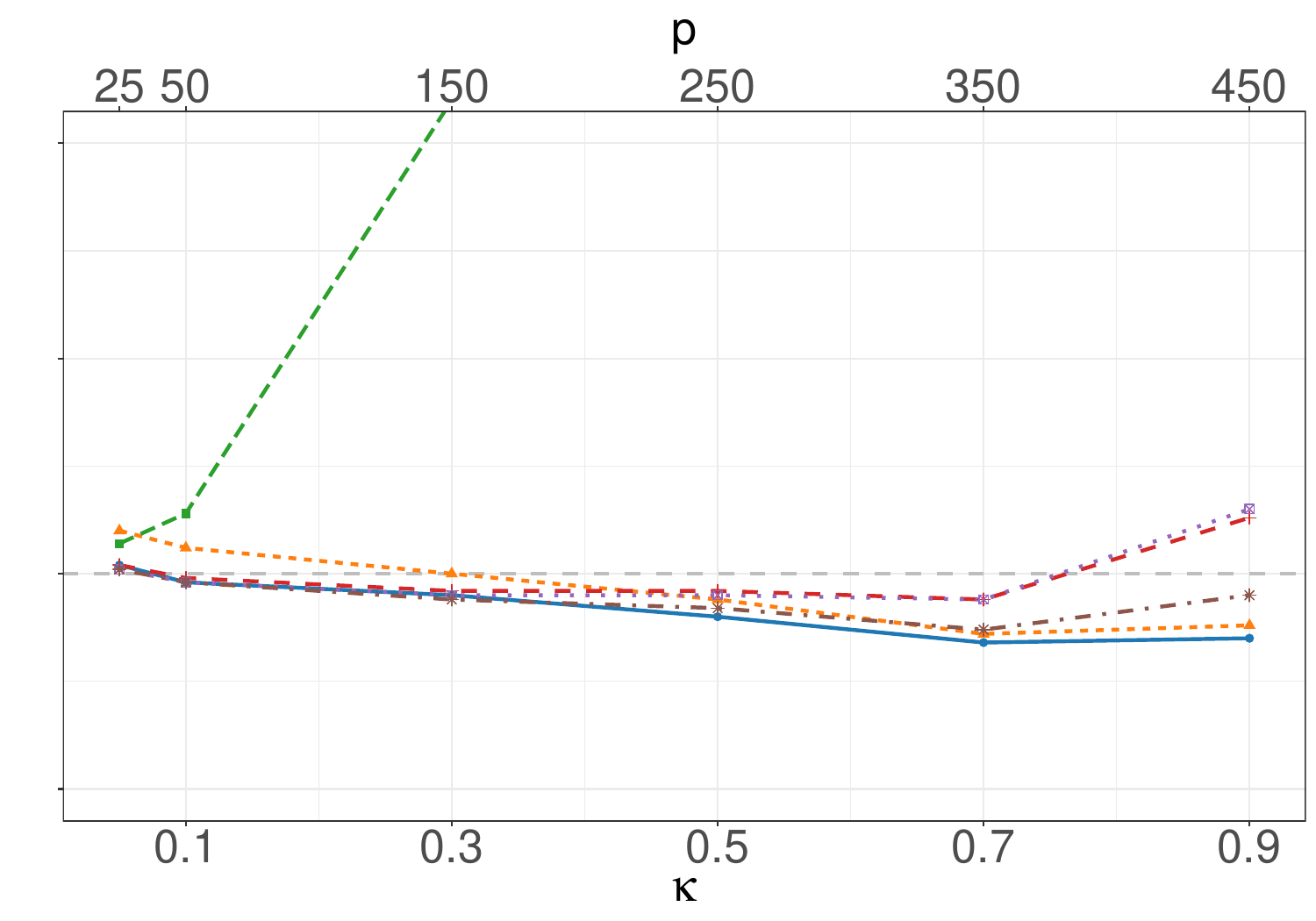}
		\end{minipage}
	}
	\caption{
		Empirical size for the directional test (DT), central limit theorem test (CLT), likelihood ratio test (LRT), Bartlett corrected test (BC) and two \citeauthor{skovgaard:2001}'s modifications  (Sko1 and Sko2, respectively) for hypothesis (\ref{hypothesis:dirmean}) with $g=3$, at nominal level $\alpha = 0.05$  given by the dashed gray horizontal line. The left, middle and right panels correspond to  multivariate $t$, multivariate skew-normal, and multivariate Laplace distributions of the true generating process, respectively, with $n_i = 500$.
	}
	\label{fig:type I error same robustness 500 high}
\end{figure}

For the homoscedastic case, Figures \ref{fig:type I error same robustness 100 high}--\ref{fig:type I error same robustness 500 high}  show the empirical size at the nominal level $\alpha = 0.05$ if the distribution is misspecified.  We see that the directional test still maintains the highest accuracy. 

\subsection{Robustness to misspecification for moderate setup}\label{section:robustness}

In this section we investigate the robustness to misspecification. The true generating processes are multivariate $t$, multivariate skew-normal or multivariate Laplace distributions. Here we setup $ p = n_i^\tau$ with $n_i \in \{100,500\}$. 

In more detail,    a multivariate $t$ distribution  with location $0_p$, scale matrix $I_p$ and degrees of freedom $5$, a multivariate skew-normal distribution  with location $1_p$, scale matrix $\Omega = (\omega_{jl})= (0.2)^{|j-l|}$ and shape parameter $1_p$, and a multivariate Laplace distribution with mean vector $1_p$ and identity covariance matrix.



\begin{figure}[t]
	\centering
	\captionsetup{font=footnotesize}
	\subfigure{
		\begin{minipage}[b]{.3\linewidth}
			\centering
			\includegraphics[scale=0.2]{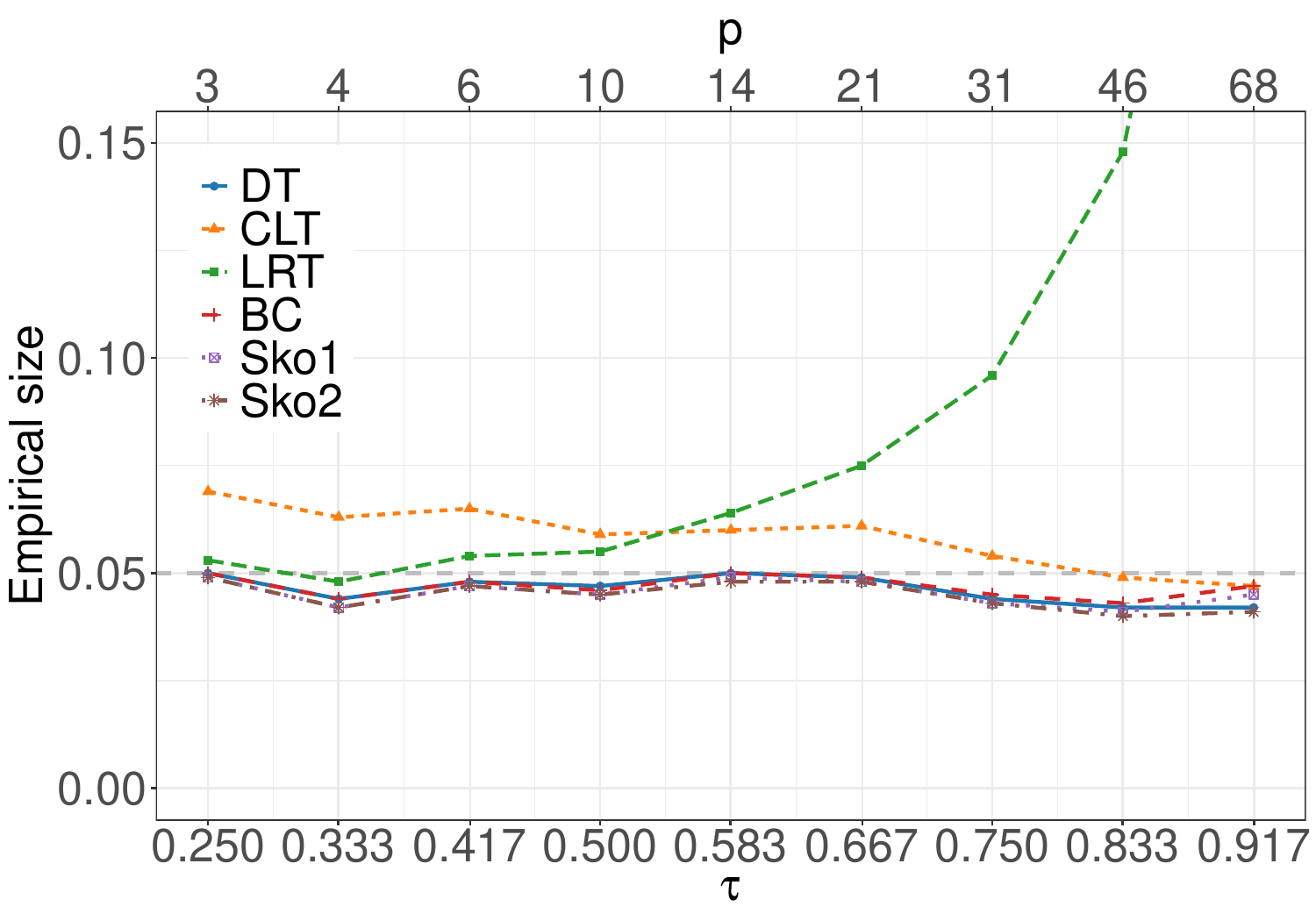}
		\end{minipage}
		\begin{minipage}[b]{.3\linewidth}
			\centering
			\includegraphics[scale=0.2]{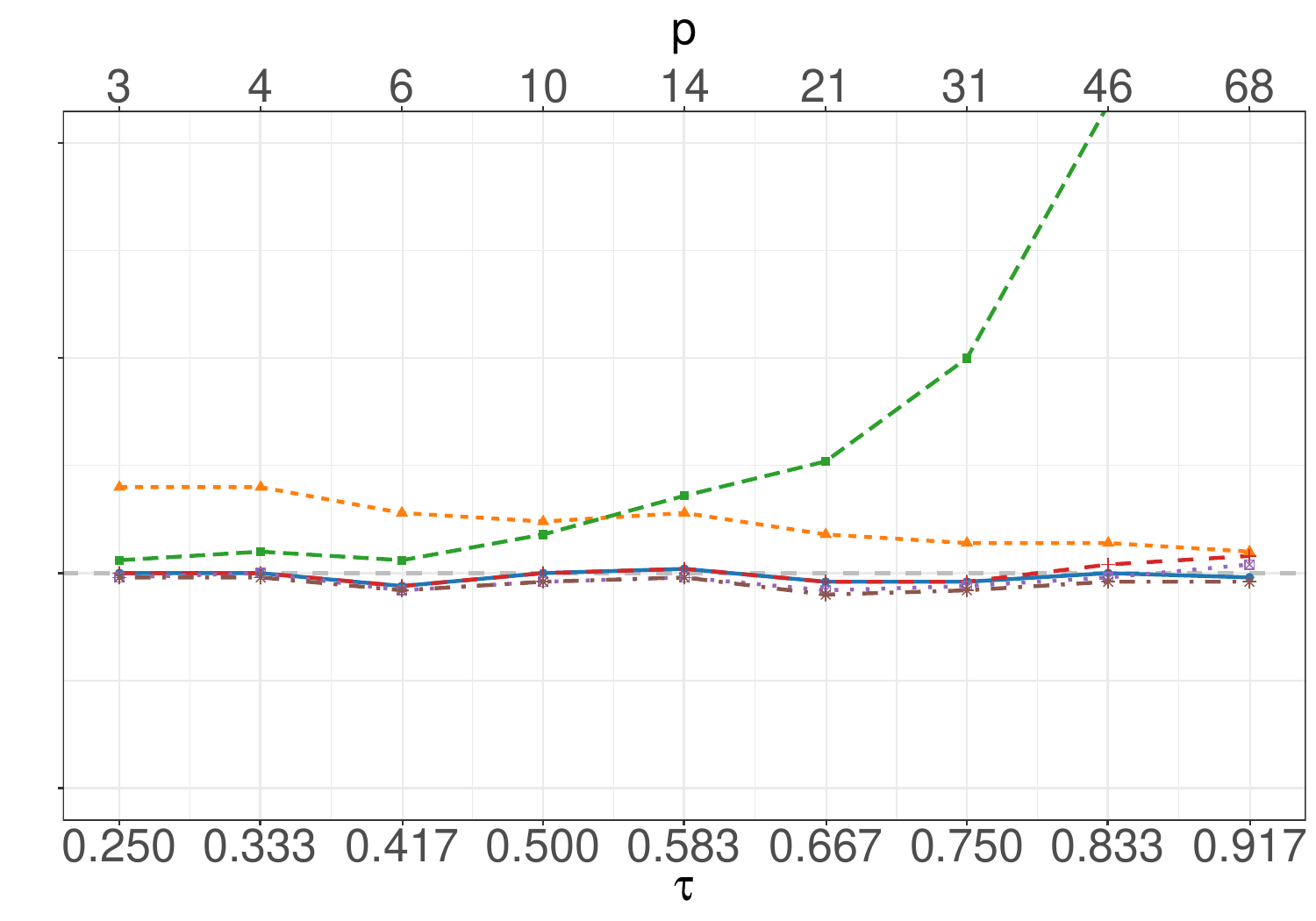}
		\end{minipage}
		\begin{minipage}[b]{.3\linewidth}
			\centering
			\includegraphics[scale=0.2]{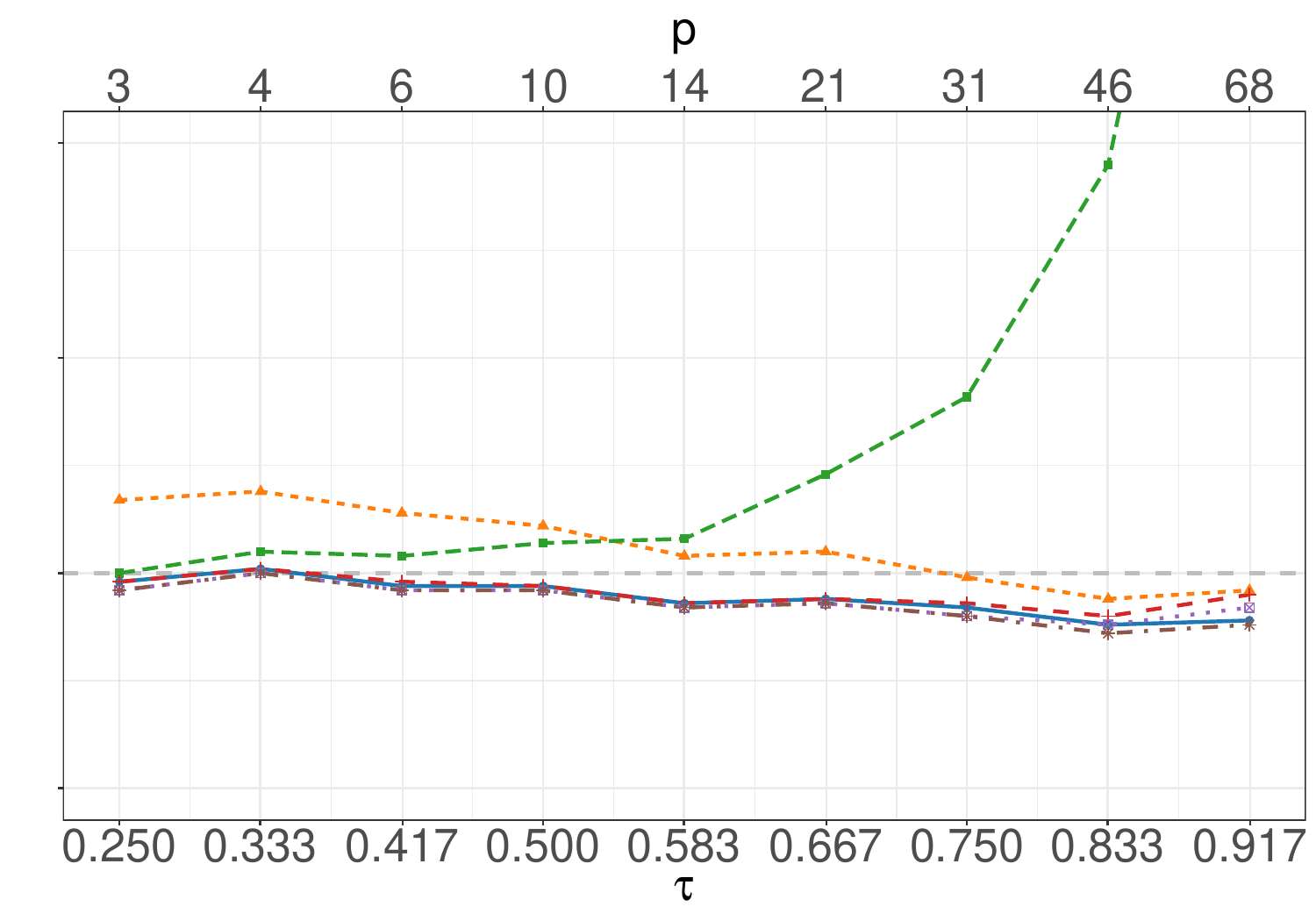}
		\end{minipage}
	}
	\caption{
		Empirical size for the directional test (DT), central limit theorem test (CLT), likelihood ratio test (LRT), Bartlett correction (BC) and two \citeauthor{skovgaard:2001}'s modifications (Sko1 and Sko2) for hypothesis (8) in the paper, at nominal level $\alpha = 0.05$  given by the gray horizontal line. The left, middle and right panels correspond to multivariate $t$, multivariate skew-normal, and multivariate Laplace distributions of the true generating process, respectively, with $n_i = 100$ and $g=3$.
	}
	\label{fig:type I error same robustness 100}
\end{figure}

\begin{figure}[t]
	\centering
	\captionsetup{font=footnotesize}
	\subfigure{
		\begin{minipage}[b]{.3\linewidth}
			\centering
			\includegraphics[scale=0.2]{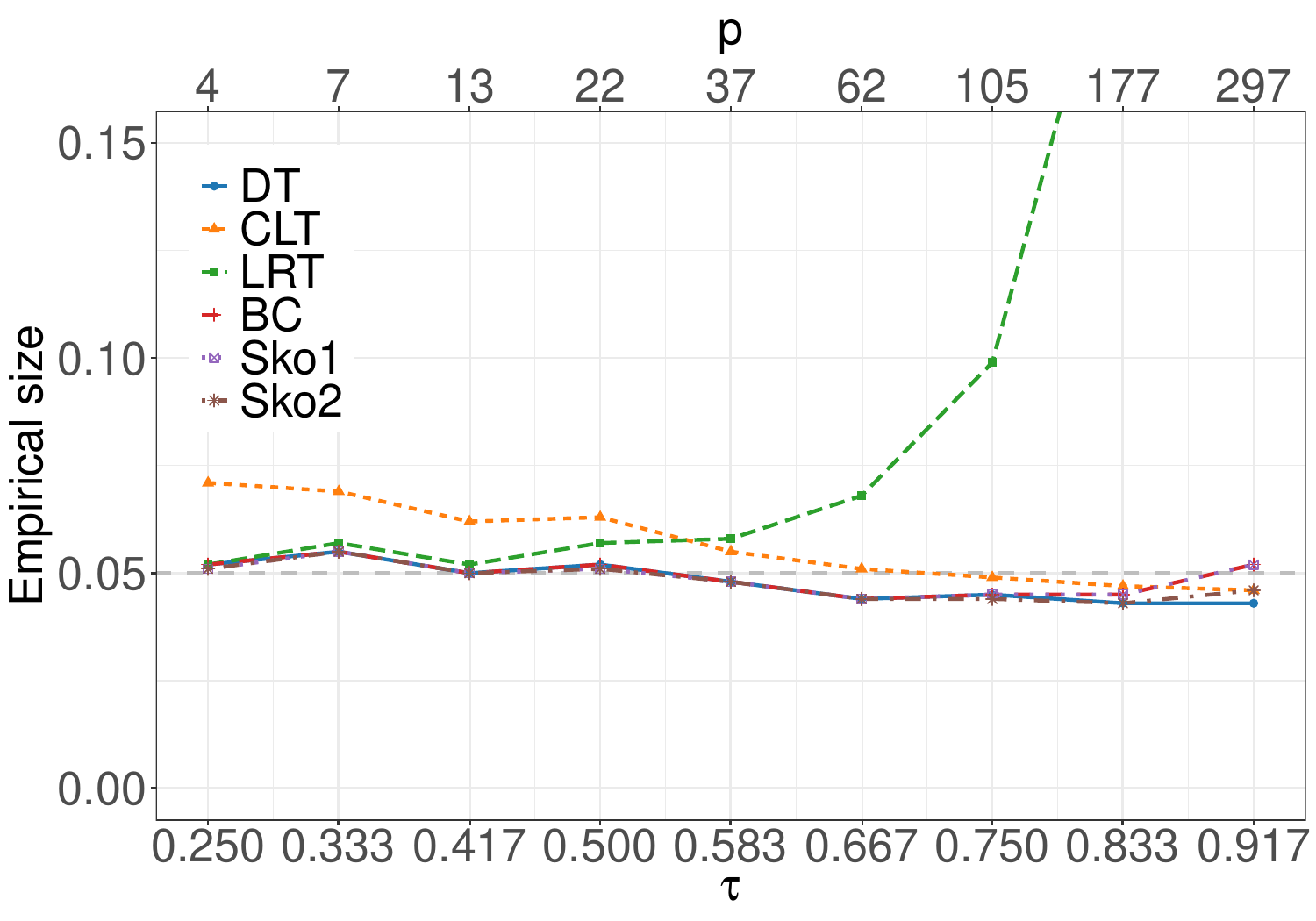}
		\end{minipage}
		\begin{minipage}[b]{.3\linewidth}
			\centering
			\includegraphics[scale=0.2]{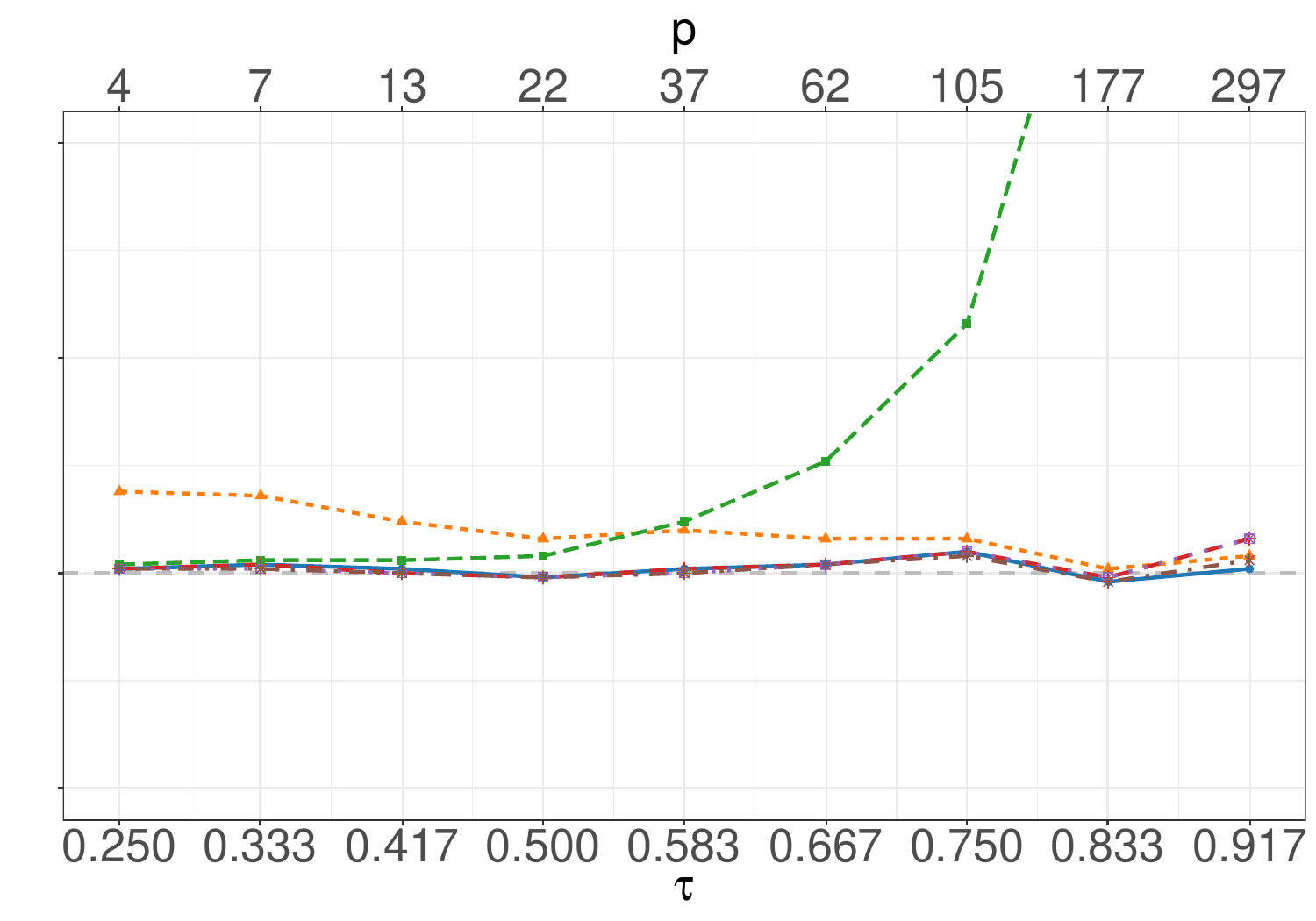}
		\end{minipage}
		\begin{minipage}[b]{.3\linewidth}
			\centering
			\includegraphics[scale=0.2]{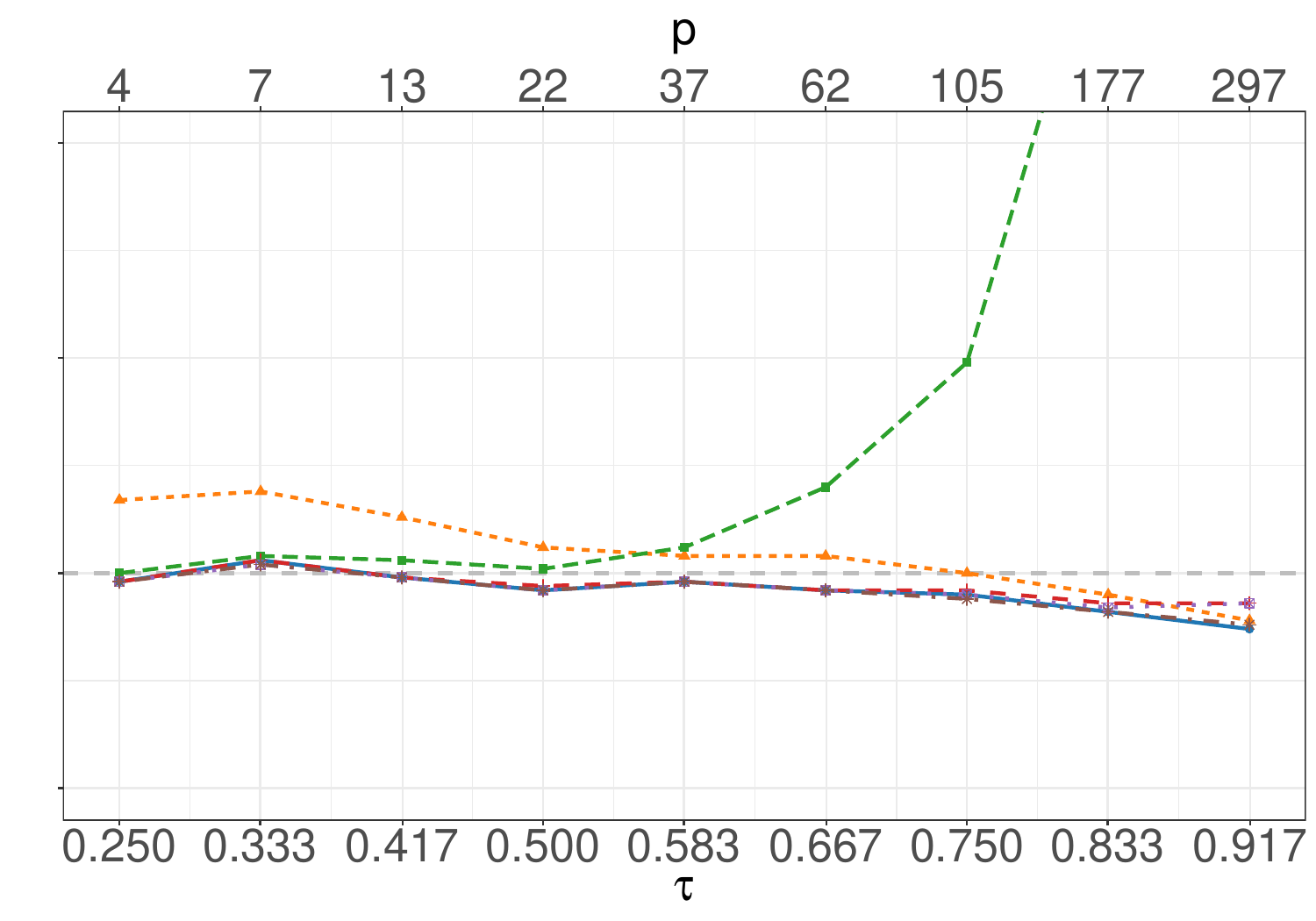}
		\end{minipage}
	}
	\caption{
		Empirical size for the directional test (DT), central limit theorem test (CLT), likelihood ratio test (LRT), Bartlett correction (BC) and two \citeauthor{skovgaard:2001}'s modifications  (Sko1 and Sko2, respectively) for hypothesis (8) in the paper, at nominal level $\alpha = 0.05$  given by the gray horizontal line. The left, middle and right panels correspond to  multivariate $t$, multivariate skew-normal, and multivariate Laplace distributions of the true generating process, respectively, with $n_i = 500$ and $g=3$.
	}
	\label{fig:type I error same robustness 500}
\end{figure}

For hypothesis (8) in the paper, Figures \ref{fig:type I error same robustness 100}--\ref{fig:type I error same robustness 500} and Tables \ref{table type I normal same he robustness 100}--\ref{table type I normal same he robustness 500} show the empirical size at the nominal level $\alpha = 0.05$ if the underling distribution is misspecified.  We see that the directional test still maintains the hightest accuracy.

\setlength{\tabcolsep}{3mm}{
	\begin{table}[H]
		\centering
		\caption{Empirical size of the directional test (DT), central limit theorem test (CLT), likelihood ratio test (LRT), Bartlett correction (BC) and two \citeauthor{skovgaard:2001}'s modifications (Sko1 and Sko2) for hypothesis (8) in the paper with  $p = n_i^\tau$ and $n_i =100$ and $g=3$, at nominal level $\alpha = 0.05$.}
		\medskip
		{\begin{tabular}{cccccccc}
				\toprule[0.09 em]
				True Distribution	& $\tau$ ($p$) &  DT & CLT & LRT &BC & Sko1 & Sko2\\
				\hline
				Multivariate $t$   & 0.250 (3) & 0.050 & 0.069 & 0.053 & 0.050 & 0.049 & 0.049 \\ 
				& 0.333 (4) & 0.044 & 0.063 & 0.048 & 0.044 & 0.042 & 0.042 \\ 
				&0.417 (6) & 0.048 & 0.065 & 0.054 & 0.048 & 0.047 & 0.047 \\ 
				&0.500 (10) & 0.047 & 0.059 & 0.055 & 0.046 & 0.045 & 0.045 \\ 
				&0.583 (14) & 0.050 & 0.060 & 0.064 & 0.050 & 0.049 & 0.048 \\ 
				&0.667 (21) & 0.049 & 0.061 & 0.075 & 0.049 & 0.048 & 0.048 \\ 
				&0.750 (31) & 0.044 & 0.054 & 0.096 & 0.045 & 0.043 & 0.043 \\ 
				&0.833 (46) & 0.042 & 0.049 & 0.148 & 0.043 & 0.041 & 0.040 \\ 
				&0.917 (68) & 0.042 & 0.047 & 0.302 & 0.047 & 0.045 & 0.041 \\
				\hline
				Multivariate skew-normal & 0.250 (3) & 0.050 & 0.070 & 0.053 & 0.050 & 0.049 & 0.049 \\ 
				& 0.333 (4) & 0.050 & 0.070 & 0.055 & 0.050 & 0.050 & 0.049 \\ 
				& 0.417 (6) & 0.047 & 0.064 & 0.053 & 0.047 & 0.046 & 0.046 \\ 
				&0.500 (10) & 0.050 & 0.062 & 0.059 & 0.050 & 0.048 & 0.048 \\ 
				&0.583 (14) & 0.051 & 0.064 & 0.068 & 0.051 & 0.049 & 0.049 \\ 
				&0.667 (21) & 0.048 & 0.059 & 0.076 & 0.048 & 0.046 & 0.045 \\ 
				&0.750 (31) & 0.048 & 0.057 & 0.100 & 0.048 & 0.047 & 0.046 \\ 
				&0.833 (46) & 0.050 & 0.057 & 0.159 & 0.052 & 0.049 & 0.048 \\ 
				&0.917 (68) & 0.049 & 0.055 & 0.312 & 0.054 & 0.052 & 0.048 \\
				\hline	
				Multivariate Laplace & 0.250 (3) & 0.048 & 0.067 & 0.050 & 0.048 & 0.046 & 0.046 \\ 
				& 0.333 (4) & 0.051 & 0.069 & 0.055 & 0.051 & 0.050 & 0.050 \\ 
				&0.417 (6) & 0.047 & 0.064 & 0.054 & 0.048 & 0.046 & 0.046 \\ 
				&0.500 (10) & 0.047 & 0.061 & 0.057 & 0.047 & 0.046 & 0.046 \\ 
				&0.583 (14) & 0.043 & 0.054 & 0.058 & 0.043 & 0.042 & 0.042 \\ 
				&0.667 (21) & 0.044 & 0.055 & 0.073 & 0.044 & 0.043 & 0.043 \\ 
				&0.750 (31) & 0.042 & 0.049 & 0.091 & 0.043 & 0.040 & 0.040 \\ 
				&0.833 (46) & 0.038 & 0.044 & 0.145 & 0.040 & 0.038 & 0.036 \\ 
				&0.917 (68) & 0.039 & 0.046 & 0.301 & 0.045 & 0.042 & 0.038 \\ 
				\bottomrule[0.09 em]            
			\end{tabular}
		}
		\label{table type I normal same he robustness 100}
\end{table}}

\setlength{\tabcolsep}{3mm}{
	\begin{table}[H]
		\centering
		\caption{Empirical size of the directional test (DT), central limit theorem test (CLT), likelihood ratio test (LRT), Bartlett correction (BC) and two \citeauthor{skovgaard:2001}'s modifications  (Sko1 and Sko2) for hypothesis (8) in the paper with $p = n_i^\tau$ and $n_i =500$ and $g=3$, at nominal level $\alpha = 0.05$.}
		\medskip
		{\begin{tabular}{cccccccc}
				\toprule[0.09 em]
				True Distribution	& $\tau$ ($p$) &  DT & CLT & LRT &BC & Sko1 & Sko2\\
				\hline
				Multivariate $t$  & 0.250 (4) & 0.052 & 0.071 & 0.052 & 0.052 & 0.051 & 0.051 \\ 
				& 0.333 (7) & 0.055 & 0.069 & 0.057 & 0.055 & 0.055 & 0.055 \\ 
				&0.417 (13) & 0.050 & 0.062 & 0.052 & 0.050 & 0.050 & 0.050 \\ 
				&0.500 (22) & 0.052 & 0.063 & 0.057 & 0.052 & 0.051 & 0.051 \\ 
				&0.583 (37) & 0.048 & 0.055 & 0.058 & 0.048 & 0.048 & 0.048 \\ 
				&0.667 (62) & 0.044 & 0.051 & 0.068 & 0.044 & 0.044 & 0.044 \\ 
				&0.750 (105) & 0.045 & 0.049 & 0.099 & 0.045 & 0.045 & 0.044 \\ 
				&0.833 (177) & 0.043 & 0.047 & 0.212 & 0.045 & 0.045 & 0.043 \\ 
				&0.917 (297) & 0.043 & 0.046 & 0.618 & 0.052 & 0.052 & 0.046 \\ 
				\hline
				Multivariate skew-normal& 0.250 (4)& 0.051 & 0.069 & 0.052 & 0.051 & 0.051 & 0.051 \\ 
				& 0.333 (7) & 0.052 & 0.068 & 0.053 & 0.052 & 0.051 & 0.051 \\ 
				&0.417 (13) & 0.051 & 0.062 & 0.053 & 0.050 & 0.050 & 0.050 \\ 
				&0.500 (22) & 0.049 & 0.058 & 0.054 & 0.049 & 0.049 & 0.049 \\ 
				&0.583 (37) & 0.051 & 0.060 & 0.062 & 0.051 & 0.050 & 0.050 \\ 
				&0.667 (62) & 0.052 & 0.058 & 0.076 & 0.052 & 0.052 & 0.052 \\ 
				&0.750 (105) & 0.055 & 0.058 & 0.108 & 0.055 & 0.055 & 0.054 \\ 
				&0.833 (177) & 0.048 & 0.051 & 0.218 & 0.049 & 0.049 & 0.048 \\ 
				&0.917 (297) & 0.051 & 0.054 & 0.612 & 0.058 & 0.058 & 0.053 \\ 
				\hline
				Multivariate Laplace & 0.250 (4) & 0.048 & 0.067 & 0.050 & 0.048 & 0.048 & 0.048 \\ 
				& 0.333 (7) & 0.053 & 0.069 & 0.054 & 0.053 & 0.052 & 0.052 \\ 
				&0.417 (13) & 0.049 & 0.063 & 0.053 & 0.049 & 0.049 & 0.049 \\ 
				&0.500 (22) & 0.046 & 0.056 & 0.051 & 0.047 & 0.046 & 0.046 \\ 
				&0.583 (37) & 0.048 & 0.054 & 0.056 & 0.048 & 0.048 & 0.048 \\ 
				&0.667 (62) & 0.046 & 0.054 & 0.070 & 0.046 & 0.046 & 0.046 \\ 
				&0.750 (105) & 0.045 & 0.050 & 0.099 & 0.046 & 0.045 & 0.044 \\ 
				&0.833 (177) & 0.041 & 0.045 & 0.209 & 0.043 & 0.042 & 0.041 \\ 
				&0.917 (297) & 0.037 & 0.039 & 0.609 & 0.043 & 0.043 & 0.038 \\ 
				\bottomrule[0.09 em]            
			\end{tabular}
		}
		\label{table type I normal same he robustness 500}
\end{table}}

\section{Simulation studies for heteroscedastic one-way MANOVA}

This section is studied the performance of directional test for heteroscedastic one-way MANOVA, comparing with LRT, Sko1 and Sko2. In particular, when the number of groups $g=2$, we also consider the $F$-approximation for the Behrens-Fisher test $T^{*2}$. The simulation results are computed via Monte Carlo simulation based on 10,000 replications.

\subsection{Empirical results for the moderate setup}
Groups of size $n_i, i \in \{1,\dots, g\}$,  are  generated from a $p$-variate standard normal distribution $N_p(0_p, \Lambda_i^{-1})$ under the null hypothesis. We use an autoregressive structure for the covariance matrices. i.e. $\Lambda_i^{-1} = (\sigma_{jl})_{p \times p} = (\rho_i^{|j-l|})_{p \times p}$, with the $\rho_i$ chosen to an equally-distance sequence from 0.1 to 0.9 of length $g$. For each simulation experiment, we show  results for $ p = \lceil n_i ^ \tau \rceil$ with $\tau = j/24$, $j \in \{6,7,\cdots,22\}$ and $n_i \in \{100, 500, 1000\}$ and $k=2$.
Note that for $n_i =1000$ the simulations results are based on  $5000$ replications when $j \in \{21, 22\}$ due to the expensive computational cost.

\setlength{\tabcolsep}{3.5mm}{
	\begin{table}[H]
		\centering
		\caption{Empirical size of  the directional test (DT), Behrens-Fisher test (BF) \citep{nel1986}, likelihood ratio test (LRT),  and two \citeauthor{skovgaard:2001}'s modifications (Sko1 and Sko2)  for hypothesis (9) in the paper with $p = n_i^\tau$ and $n_i \in \{ 100, 500, 1000\}$, at nominal level $\alpha = 0.05$.}
		\medskip
		{\begin{tabular}{ccccccc}
				\toprule[0.09 em]
				$n_i$	& $\tau$ ($p$) &  DT & BF & LRT  & Sko1 & Sko2\\
				\hline
				100& 0.250 (4) & 0.053 & 0.054 & 0.057 & 0.053 & 0.053 \\ 
				& 0.333 (5) & 0.050 & 0.052 & 0.056 & 0.050 & 0.050 \\ 
				& 0.417 (7) & 0.049 & 0.051 & 0.057 & 0.050 & 0.050 \\ 
				&0.500 (10) & 0.051 & 0.056 & 0.065 & 0.052 & 0.052 \\ 
				&0.583 (15) & 0.052 & 0.061 & 0.077 & 0.052 & 0.052 \\ 
				&0.667 (22) & 0.049 & 0.065 & 0.094 & 0.051 & 0.050 \\ 
				&0.750 (32) & 0.051 & 0.082 & 0.147 & 0.055 & 0.053 \\ 
				&0.833 (47) & 0.052 & 0.115 & 0.270 & 0.067 & 0.061 \\ 
				&0.917 (69) & 0.064 & 0.183 & 0.594 & 0.112 & 0.092 \\ 
				\hline
				500 & 0.250 (5) & 0.050 & 0.050 & 0.051 & 0.050 & 0.050 \\ 
				& 0.333 (8) & 0.051 & 0.051 & 0.052 & 0.051 & 0.051 \\ 
				&0.417 (14) & 0.051 & 0.052 & 0.055 & 0.051 & 0.051 \\ 
				&0.500 (23) & 0.051 & 0.054 & 0.059 & 0.051 & 0.051 \\ 
				&0.583 (38) & 0.054 & 0.061 & 0.070 & 0.054 & 0.054 \\ 
				&0.667 (63) & 0.051 & 0.068 & 0.089 & 0.052 & 0.052 \\ 
				&0.750 (106) & 0.048 & 0.084 & 0.143 & 0.051 & 0.050 \\ 
				&0.833 (178) & 0.051 & 0.154 & 0.382 & 0.062 & 0.058 \\ 
				&0.917 (298) & 0.060 & 0.392 & 0.923 & 0.137 & 0.106 \\ 
				\hline
				1000 & 0.250 (6) & 0.052 & 0.052 & 0.053 & 0.052 & 0.052 \\ 
				&0.333 (10) & 0.050 & 0.050 & 0.051 & 0.050 & 0.050 \\ 
				&0.417 (18) & 0.046 & 0.048 & 0.048 & 0.046 & 0.046 \\ 
				&0.500 (32) & 0.052 & 0.055 & 0.058 & 0.052 & 0.052 \\ 
				&0.583 (57) & 0.052 & 0.058 & 0.063 & 0.052 & 0.052 \\ 
				&0.667 (100) & 0.050 & 0.066 & 0.083 & 0.050 & 0.050 \\ 
				&0.750 (178) & 0.047 & 0.088 & 0.154 & 0.050 & 0.049 \\ 
				&0.833 (317) & 0.051 & 0.179 & 0.459 & 0.065 & 0.059 \\ 
				&0.917 (563) & 0.064 & 0.528 & 0.987 & 0.155 & 0.112 \\ 	
				\bottomrule[0.09 em]            
			\end{tabular}
		}
		\label{table type I normal full he}
\end{table}}

\begin{figure}[t]
	\centering
	\captionsetup{font=footnotesize}
	\subfigure{
		\begin{minipage}[b]{.9\linewidth}
			\centering
			\includegraphics[scale=0.5]{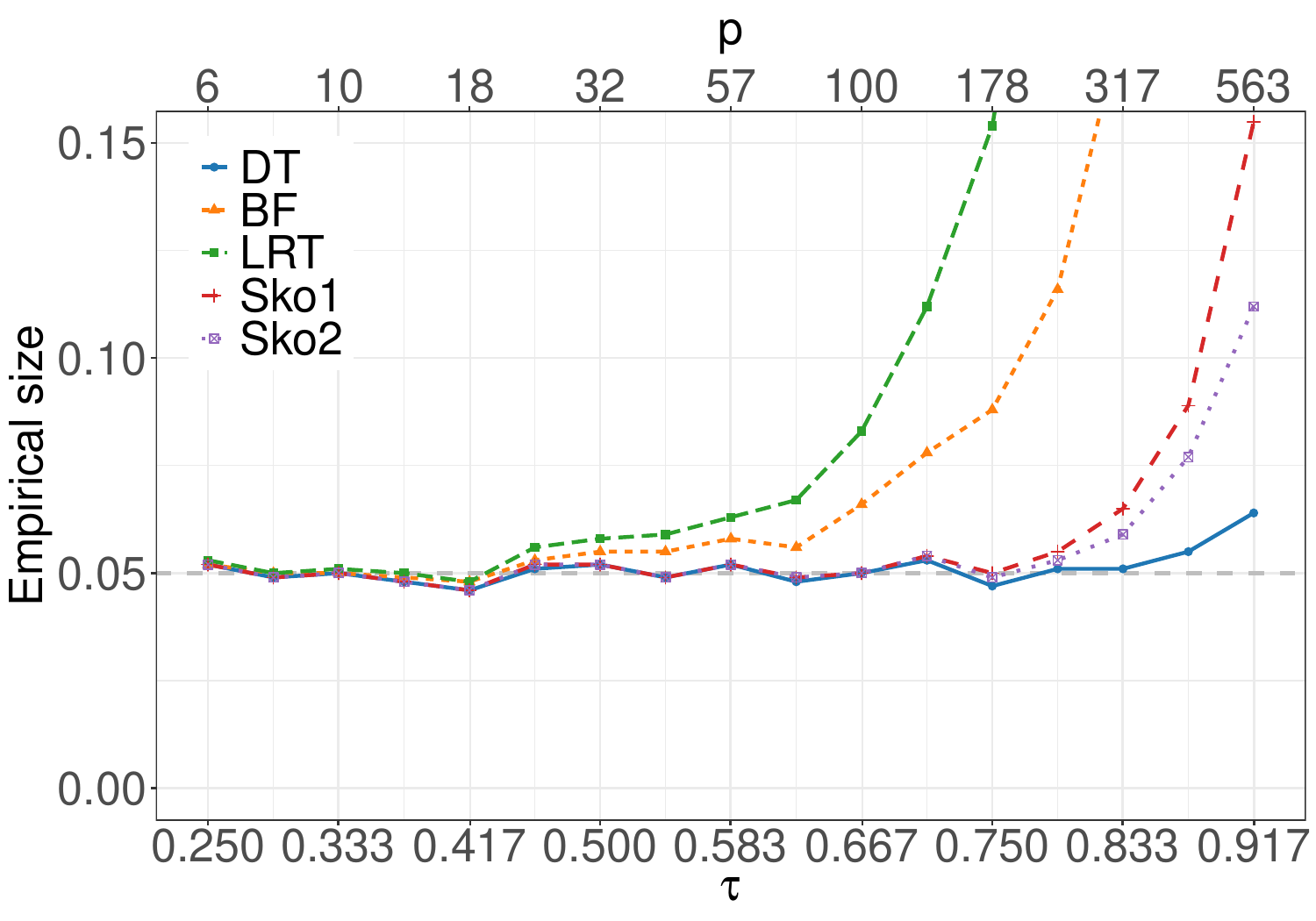}
		\end{minipage}
	}
	
	\caption{
		Empirical size of the directional test (DT), Behrens-Fisher test (BF) \citep{nel1986}, likelihood ratio test (LRT),  and two \citeauthor{skovgaard:2001}'s modifications (Sko1 and Sko2) for hypothesis (\ref{hypothesis:meansdiff}) with {$g=2$}, at nominal level $\alpha = 0.05$  given by the dashed gray horizontal line. The left, middle and right panels correspond to $n_i = 1000$, respectively.
	}
	\label{fig:type I error full}
\end{figure}

\subsection{Empirical results for a large number $g$ of groups}

\setlength{\tabcolsep}{3.5mm}{
	\begin{table}[H]
		\centering
		\caption{Empirical size of the directional test (DT), likelihood ratio test (LRT), and two \citeauthor{skovgaard:2001}'s modifications (Sko1 and Sko2) for hypothesis (9) in the paper with $p = n_i^\tau$ and $n_i = 100$ and $g =30$, at  nominal level $\alpha = 0.05$.}
		\medskip
		{\begin{tabular}{ccccc}
				\toprule[0.09 em]
				$\tau$ ($p$) &  DT  & LRT  & Sko1 & Sko2\\
				\hline 
				0.250 (4)& 0.050 & 0.081 & 0.050 & 0.050 \\ 
				0.333 (5) & 0.050 & 0.091 & 0.050 & 0.050 \\ 
				0.417 (7) & 0.049 & 0.122 & 0.050 & 0.049 \\ 
				0.500 (10) & 0.050 & 0.186 & 0.051 & 0.050 \\ 
				0.583 (15) & 0.054 & 0.391 & 0.058 & 0.054 \\ 
				0.667 (22)& 0.047 & 0.751 & 0.056 & 0.048 \\ 
				0.750 (32) & 0.052 & 0.990 & 0.087 & 0.062 \\ 
				0.833 (47) & 0.049 & 1.000 & 0.223 & 0.095 \\ 
				0.917 (69) & 0.048 & 1.000 & 0.932 & 0.405 \\
				\bottomrule[0.09 em]            
			\end{tabular}
		}
		\label{table type I normal full he100 k30}
\end{table}}

\begin{figure}[H]
	\centering
	\captionsetup{font=footnotesize}
	\includegraphics[scale=0.4]{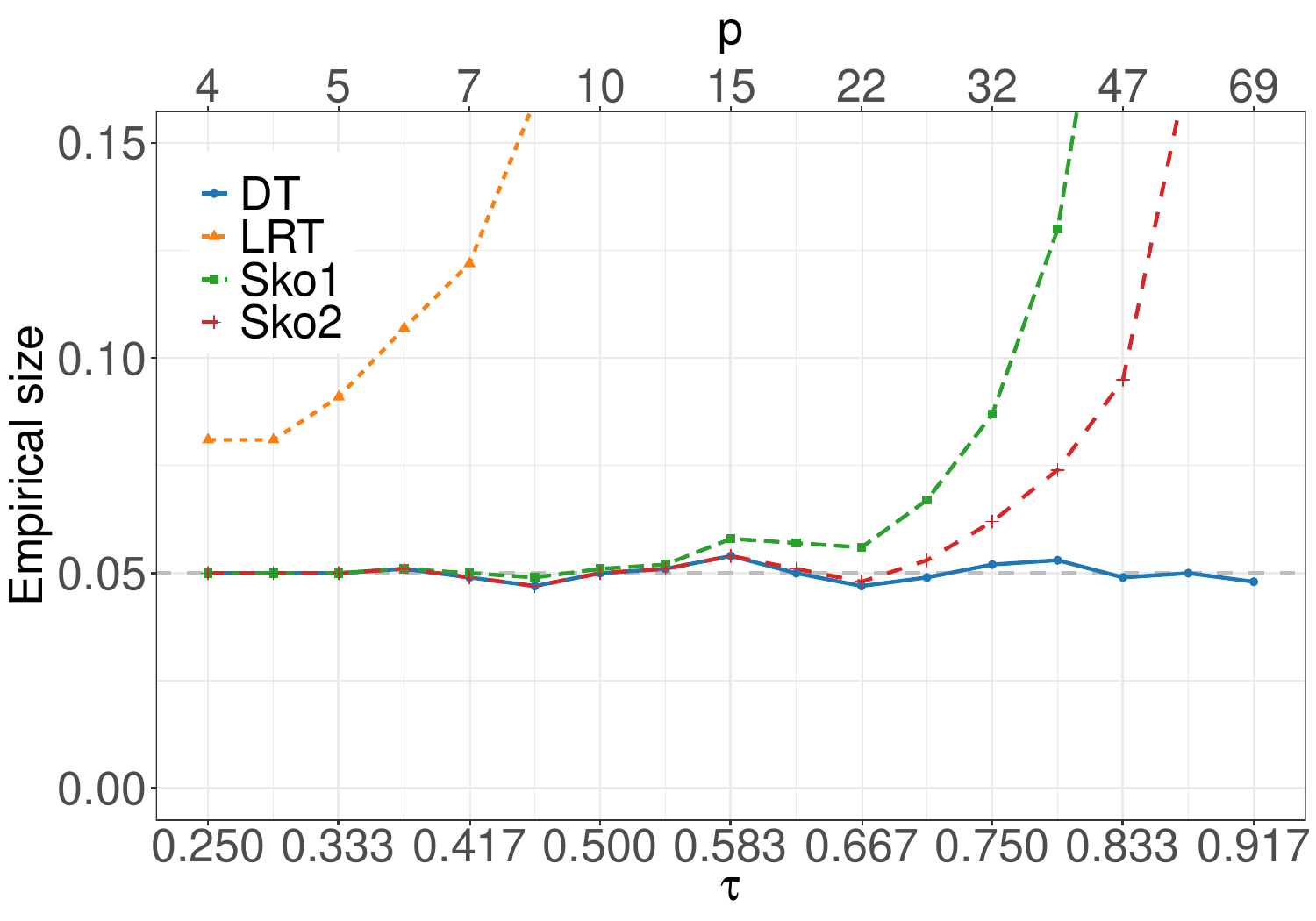}
	\caption{Empirical size of the directional test (DT), likelihood ratio test (LRT), and two \citeauthor{skovgaard:2001}'s modifications  (Sko1 and Sko2) for hypothesis (9) in the paper with $p=n_i^\tau$ and  $n_i = 100$ and $g =30$, at nominal level $\alpha = 0.05$  given by the gray horizontal line.}
	\label{figure:type I error full he100 k30}
\end{figure}

\setlength{\tabcolsep}{3.5mm}{
	\begin{table}[H]
		\centering
		\caption{Empirical size of the directional test (DT), likelihood ratio test (LRT), and two \citeauthor{skovgaard:2001}'s modifications  (Sko1 and Sko2) for hypothesis (9) in the paper with $p = n_i^\tau$ and $n_i = 1000$ and $g=5$, at  nominal level $\alpha = 0.05$}
		\medskip
		{\begin{tabular}{ccccc}
				\toprule[0.09 em]
				$\tau$ ($p$) &  DT  & LRT  & Sko1 & Sko2\\
				\hline
				0.250 (6) & 0.053 & 0.055 & 0.053 & 0.053 \\ 
				0.333 (10) & 0.055 & 0.058 & 0.055 & 0.055 \\ 
				0.417 (18) & 0.048 & 0.054 & 0.048 & 0.048 \\ 
				0.500 (32) & 0.048 & 0.061 & 0.048 & 0.048 \\ 
				0.583 (57) & 0.051 & 0.089 & 0.052 & 0.052 \\ 
				0.667 (100) & 0.054 & 0.160 & 0.055 & 0.054 \\ 
				0.750 (178) & 0.053 & 0.462 & 0.058 & 0.055 \\ 
				0.833 (317) & 0.050 & 0.983 & 0.083 & 0.063 \\ 
				\bottomrule[0.09 em]            
			\end{tabular}
		}
		\label{table type I normal full he1000 k5}
\end{table}}

\begin{figure}[H]
	\centering
	\captionsetup{font=footnotesize}
	\includegraphics[scale=0.4]{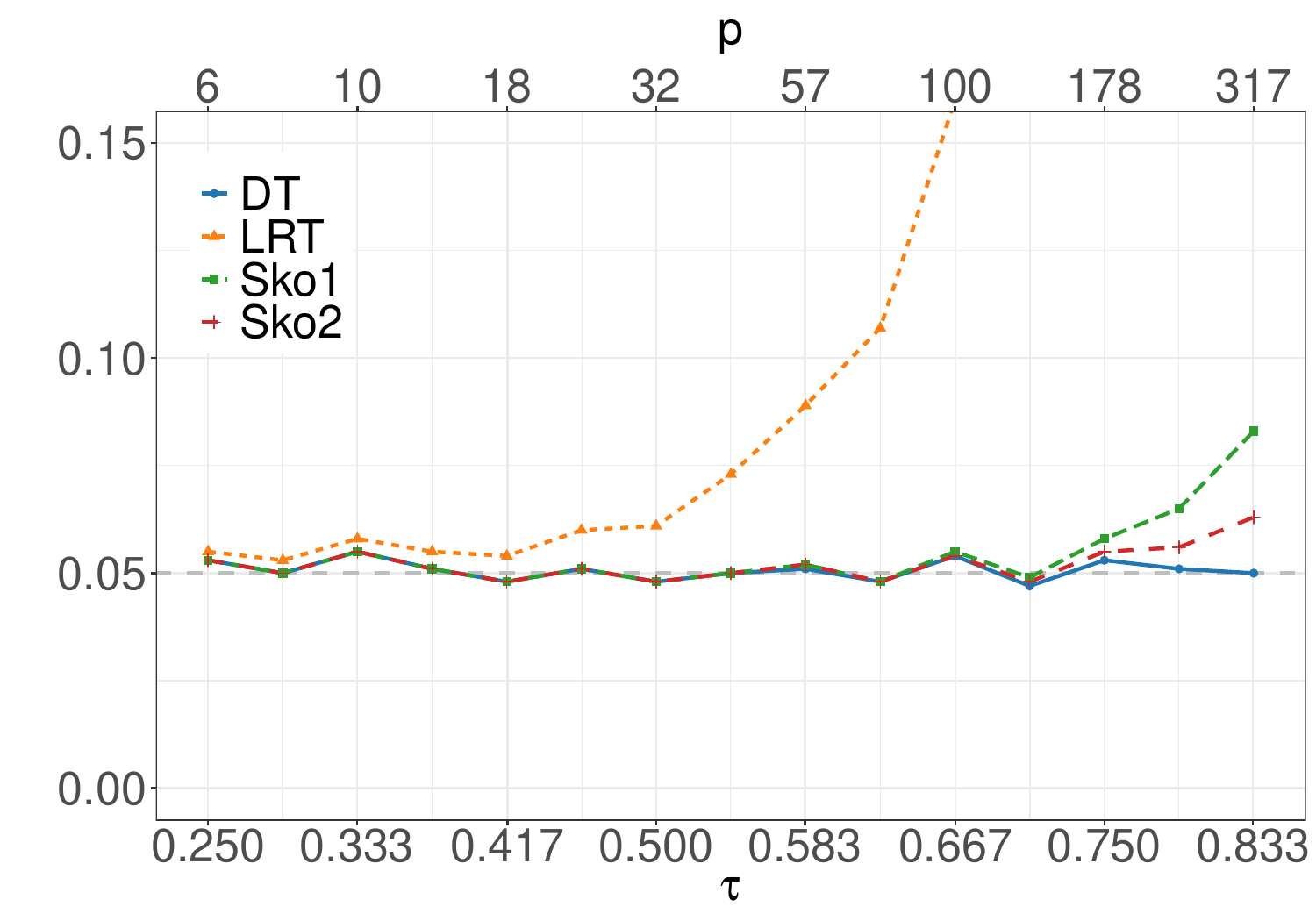}
	\caption{ Empirical size of the directional test (DT), likelihood ratio test (LRT), and two \citeauthor{skovgaard:2001}'s modifications (Sko1 and Sko2) for hypothesis (9) in the paper, with $p = n_i^\tau$ and $n_i = 1000$ and $g =5$, at nominal level $\alpha = 0.05$  given by the gray horizontal line.}
	\label{figure:type I error full he1000 k5}
\end{figure}


\subsection{Robustness to misspecification }


\setlength{\tabcolsep}{3.5mm}{
	\begin{table}[H]
		\centering
		\caption{Empirical size of the directional test (DT), Behrens-Fisher test (BF) \citep{nel1986}, likelihood ratio test (LRT),  and two \citeauthor{skovgaard:2001}'s modifications (Sko1 and Sko2)  for hypothesis (9) in the paper with $p = n_i^\tau$ and $n_i =100$ and $g=2$, at nominal level $\alpha = 0.05$.}
		\medskip
		{\begin{tabular}{ccccccc}
				\toprule[0.09 em]
				Distribution	& $\tau$ ($p$) &  DT & BF & LRT  & Sko1 & Sko2\\
				\hline
				Multivariate $t$ & 0.250 (4) & 0.048 & 0.048 & 0.052 & 0.048 & 0.048 \\ 
				& 0.333 (5) & 0.047 & 0.047 & 0.053 & 0.047 & 0.047 \\ 
				& 0.417 (7) & 0.048 & 0.048 & 0.055 & 0.048 & 0.048 \\ 
				&0.500 (10) & 0.050 & 0.049 & 0.060 & 0.050 & 0.050 \\ 
				&0.583 (15) & 0.045 & 0.044 & 0.064 & 0.045 & 0.045 \\ 
				&0.667 (22) & 0.047 & 0.046 & 0.078 & 0.047 & 0.047 \\ 
				&0.750 (32) & 0.042 & 0.040 & 0.102 & 0.043 & 0.042 \\ 
				&0.833 (47) & 0.041 & 0.038 & 0.166 & 0.044 & 0.042 \\ 
				&0.917 (69) & 0.043 & 0.038 & 0.375 & 0.054 & 0.048 \\ 
				\hline
				Multivariate skew-normal & 0.250 (4) & 0.046 & 0.046 & 0.051 & 0.046 & 0.046 \\ 
				& 0.333 (5) & 0.052 & 0.052 & 0.058 & 0.052 & 0.052 \\ 
				& 0.417 (7) & 0.053 & 0.053 & 0.061 & 0.053 & 0.053 \\ 
				&0.500 (10) & 0.054 & 0.054 & 0.064 & 0.054 & 0.054 \\ 
				&0.583 (15) & 0.049 & 0.049 & 0.070 & 0.050 & 0.049 \\ 
				&0.667 (22) & 0.048 & 0.047 & 0.080 & 0.048 & 0.048 \\ 
				&0.750 (32) & 0.050 & 0.048 & 0.110 & 0.051 & 0.050 \\ 
				&0.833 (47) & 0.052 & 0.051 & 0.187 & 0.056 & 0.053 \\ 
				&0.917 (69) & 0.048 & 0.048 & 0.379 & 0.061 & 0.053 \\   
				\hline
				Multivariate Laplace& 0.250 (4) & 0.049 & 0.049 & 0.054 & 0.049 & 0.049 \\ 
				& 0.333 (5) & 0.044 & 0.044 & 0.050 & 0.044 & 0.044 \\ 
				& 0.417 (7) & 0.045 & 0.045 & 0.051 & 0.045 & 0.045 \\ 
				&0.500 (10) & 0.046 & 0.046 & 0.055 & 0.046 & 0.045 \\ 
				&0.583 (15) & 0.045 & 0.045 & 0.063 & 0.045 & 0.045 \\ 
				&0.667 (22) & 0.044 & 0.043 & 0.075 & 0.044 & 0.044 \\ 
				&0.750 (32) & 0.040 & 0.039 & 0.097 & 0.042 & 0.040 \\ 
				&0.833 (47) & 0.039 & 0.038 & 0.169 & 0.043 & 0.040 \\ 
				&0.917 (69) & 0.034 & 0.032 & 0.372 & 0.044 & 0.038 \\	
				\bottomrule[0.09 em]            
			\end{tabular}
		}
		\label{table type I normal full he robustness 100}
\end{table}}

\setlength{\tabcolsep}{3.5mm}{
	\begin{table}[H]
		\centering
		\caption{Empirical size of  the directional test (DT), Behrens-Fisher test (BF) \citep{nel1986}, likelihood ratio test (LRT),  and two \citeauthor{skovgaard:2001}'s modifications (Sko1 and Sko2)  for hypothesis (9) in the paper with $p = n_i^\tau$ and $n_i =500$ and $g=2$, at nominal level $\alpha = 0.05$.}
		\medskip
		{\begin{tabular}{ccccccc}
				\toprule[0.09 em]
				Distribution	& $\tau$ ($p$) &  DT & BF & LRT  & Sko1 & Sko2\\
				\hline
				Multivariate $t$ & 0.250 (5) & 0.050 & 0.050 & 0.052 & 0.050 & 0.050 \\ 
				& 0.333 (8) & 0.049 & 0.049 & 0.050 & 0.049 & 0.049 \\ 
				&0.417 (14) & 0.047 & 0.047 & 0.050 & 0.047 & 0.047 \\ 
				&0.500 (23) & 0.045 & 0.045 & 0.050 & 0.045 & 0.045 \\ 
				&0.583 (38) & 0.051 & 0.051 & 0.063 & 0.051 & 0.051 \\ 
				&0.667 (63) & 0.044 & 0.044 & 0.070 & 0.044 & 0.044 \\ 
				&0.750 (106) & 0.044 & 0.043 & 0.108 & 0.044 & 0.044 \\ 
				&0.833 (178) & 0.044 & 0.043 & 0.236 & 0.048 & 0.046 \\ 
				&0.917 (298) & 0.047 & 0.045 & 0.705 & 0.064 & 0.053 \\  
				\hline
				Multivariate skew-normal & 0.250 (5) & 0.051 & 0.051 & 0.052 & 0.051 & 0.051 \\ 
				&0.333 (8) & 0.051 & 0.051 & 0.052 & 0.051 & 0.051 \\ 
				&0.417 (14) & 0.055 & 0.055 & 0.059 & 0.055 & 0.055 \\ 
				&0.500 (23) & 0.054 & 0.054 & 0.060 & 0.054 & 0.054 \\ 
				&0.583 (38) & 0.051 & 0.051 & 0.061 & 0.051 & 0.051 \\ 
				&0.667 (63) & 0.049 & 0.049 & 0.074 & 0.049 & 0.049 \\ 
				&0.750 (106) & 0.051 & 0.051 & 0.117 & 0.052 & 0.052 \\ 
				&0.833 (178) & 0.048 & 0.048 & 0.242 & 0.054 & 0.050 \\ 
				&0.917 (298) & 0.052 & 0.051 & 0.699 & 0.070 & 0.058 \\  
				\hline
				Multivariate Laplace& 0.250 (5) & 0.054 & 0.054 & 0.054 & 0.053 & 0.053 \\ 
				& 0.333 (8) & 0.053 & 0.053 & 0.054 & 0.053 & 0.053 \\ 
				&0.417 (14) & 0.050 & 0.050 & 0.052 & 0.050 & 0.050 \\ 
				&0.500 (23) & 0.051 & 0.051 & 0.056 & 0.051 & 0.051 \\ 
				&0.583 (38) & 0.047 & 0.047 & 0.056 & 0.047 & 0.047 \\ 
				&0.667 (63) & 0.046 & 0.046 & 0.070 & 0.047 & 0.046 \\ 
				&0.750 (106) & 0.042 & 0.042 & 0.106 & 0.043 & 0.043 \\ 
				&0.833 (178) & 0.040 & 0.040 & 0.235 & 0.043 & 0.040 \\ 
				&0.917 (298) & 0.033 & 0.032 & 0.708 & 0.051 & 0.040 \\ 	
				\bottomrule[0.09 em]            
			\end{tabular}
		}
		\label{table type I normal full he robustness 500}
\end{table}}

\begin{figure}[th]
	\centering
	\captionsetup{font=footnotesize}
	\subfigure{
		\begin{minipage}[b]{.32\linewidth}
			\centering
			\includegraphics[scale=0.2]{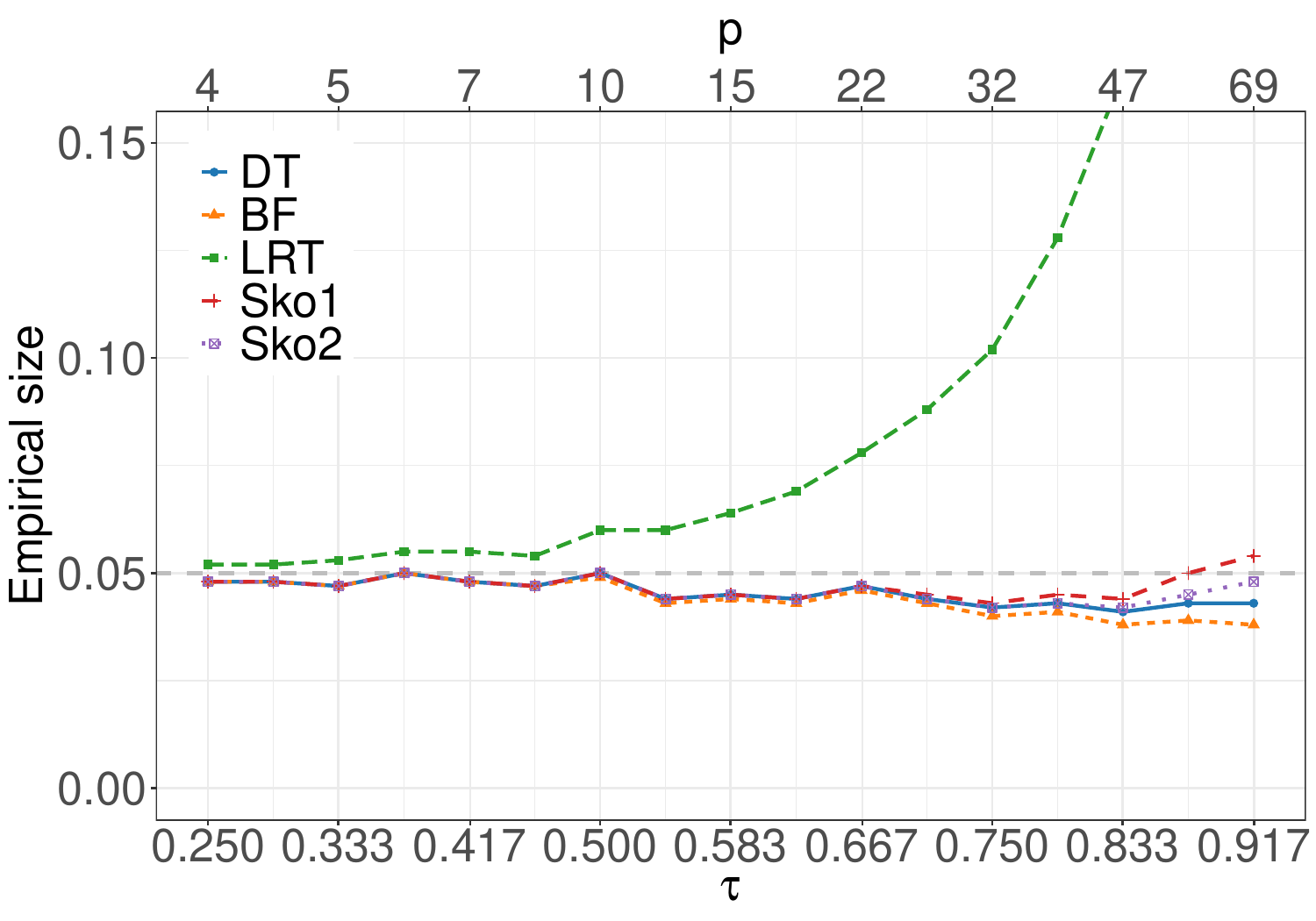}
		\end{minipage}
		\begin{minipage}[b]{.33\linewidth}
			\centering
			\includegraphics[scale=0.2]{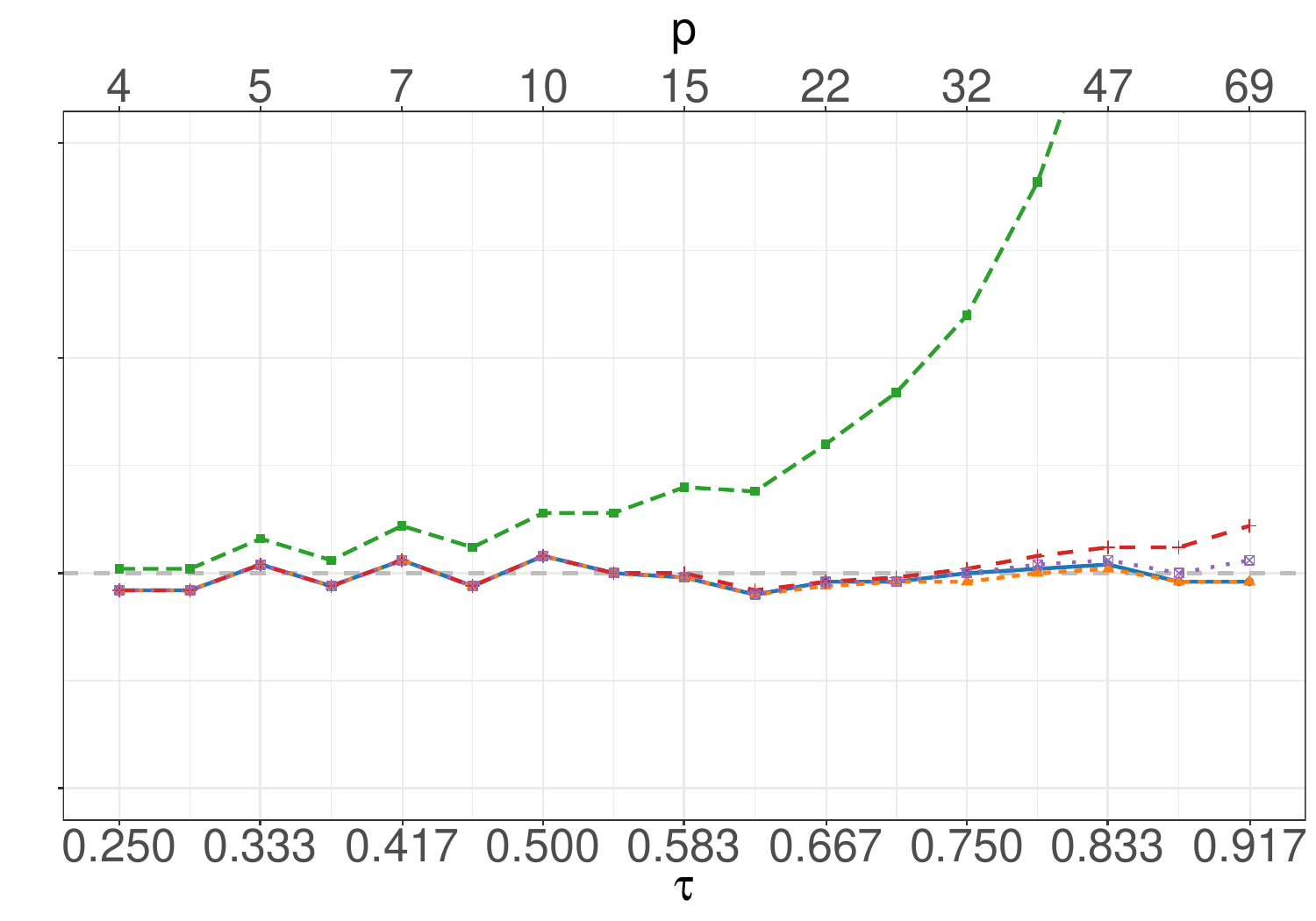}
		\end{minipage}
		\begin{minipage}[b]{.33\linewidth}
			\centering
			\includegraphics[scale=0.2]{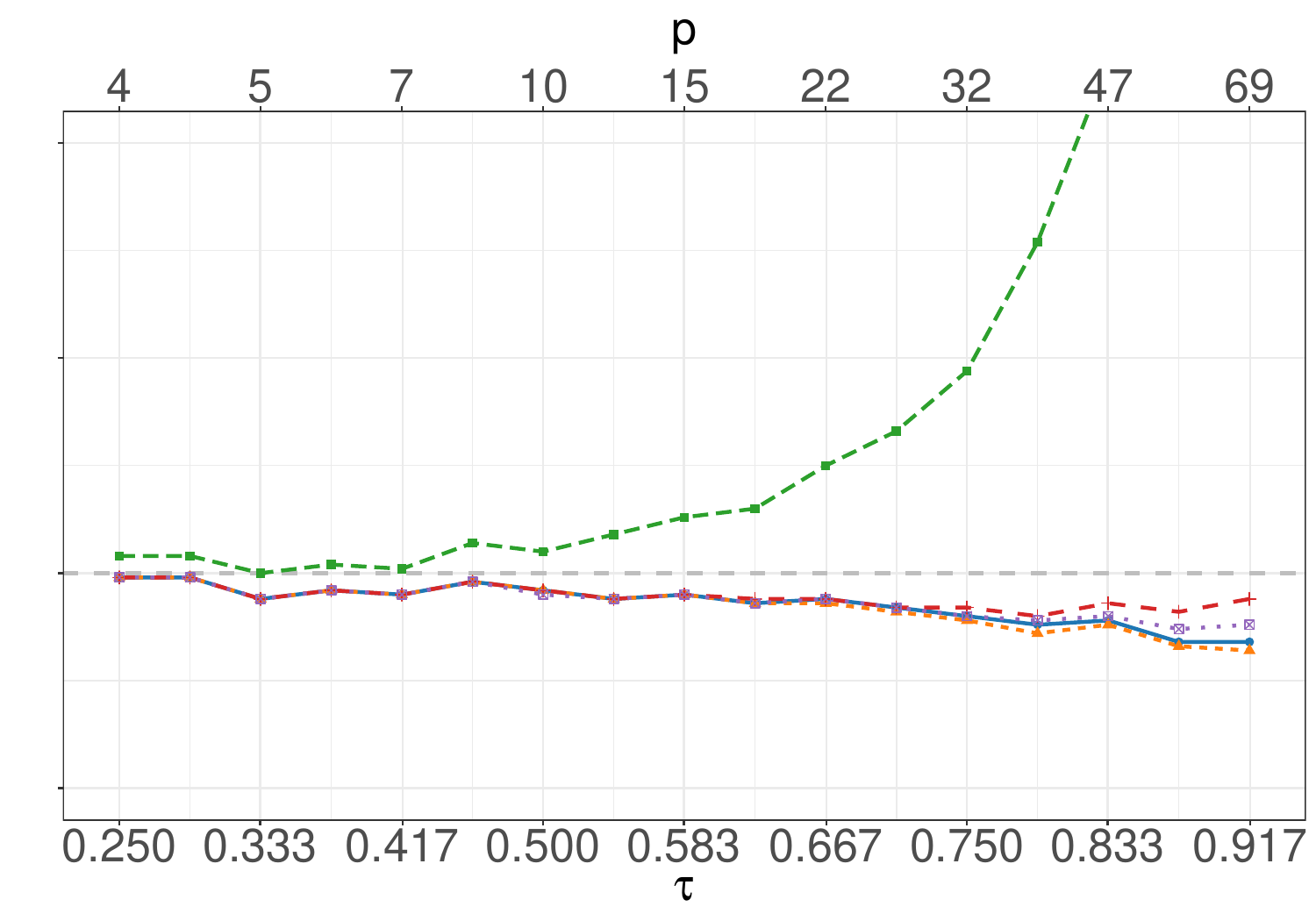}
		\end{minipage}
	}
	\caption{
		Empirical size for the directional test (DT), Behrens-Fisher test (BF), likelihood ratio test (LRT),  and two \citeauthor{skovgaard:2001}'s modifications (Sko1 and Sko2) for hypothesis (9) with $g=2$, at nominal level $\alpha = 0.05$  given by the dashed gray horizontal line. The left, middle and right panels correspond to multivariate $t$, multivariate skew-normal, and multivariate Laplace distributions if the true generating process, respectively, with $n_i = 100$.
	}
	\label{fig:type I error full robustness 100}
\end{figure}

\begin{figure}[t]
	\centering
	\captionsetup{font=footnotesize}
	\subfigure{
		\begin{minipage}[b]{.32\linewidth}
			\centering
			\includegraphics[scale=0.2]{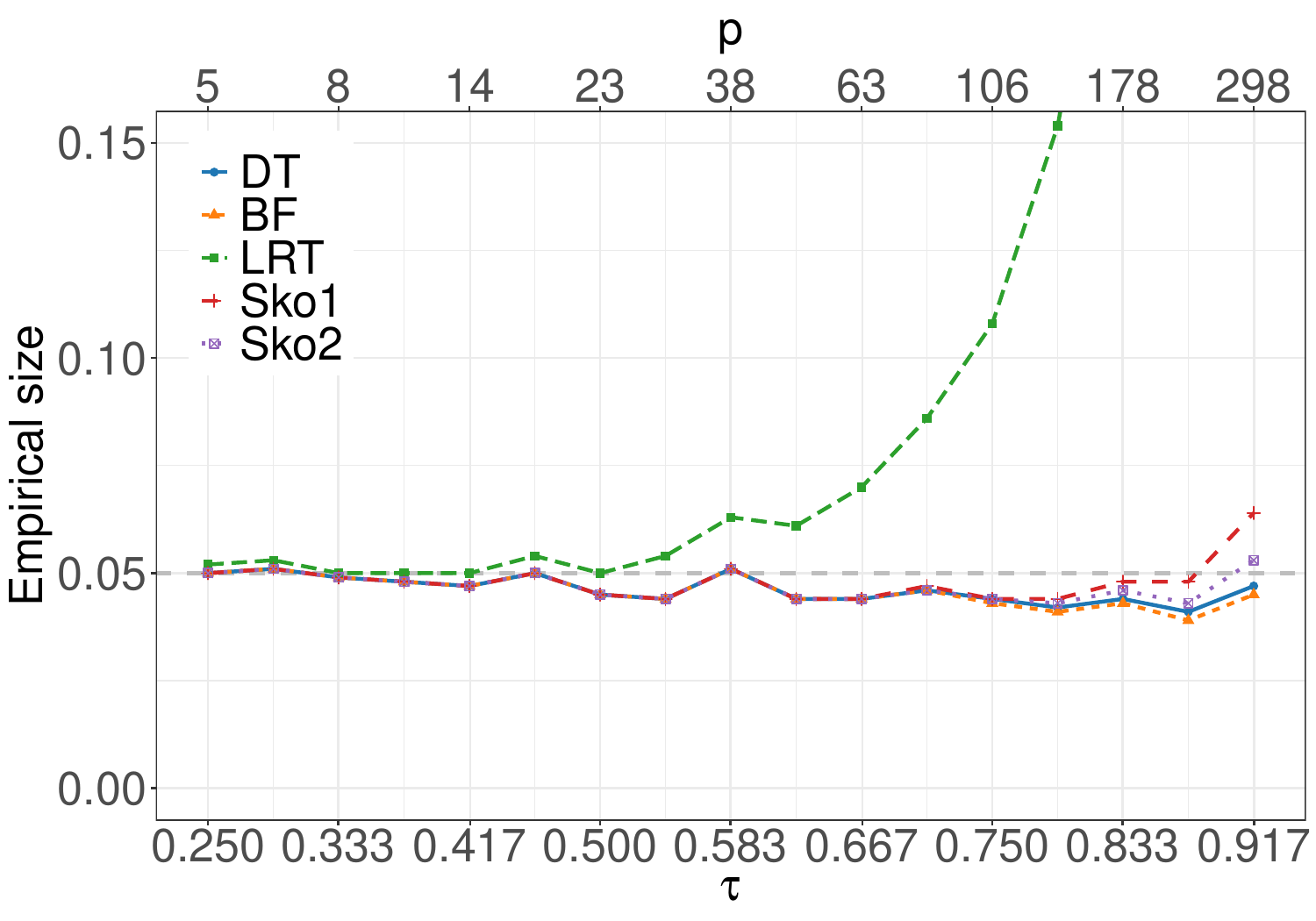}
		\end{minipage}
		\begin{minipage}[b]{.33\linewidth}
			\centering
			\includegraphics[scale=0.2]{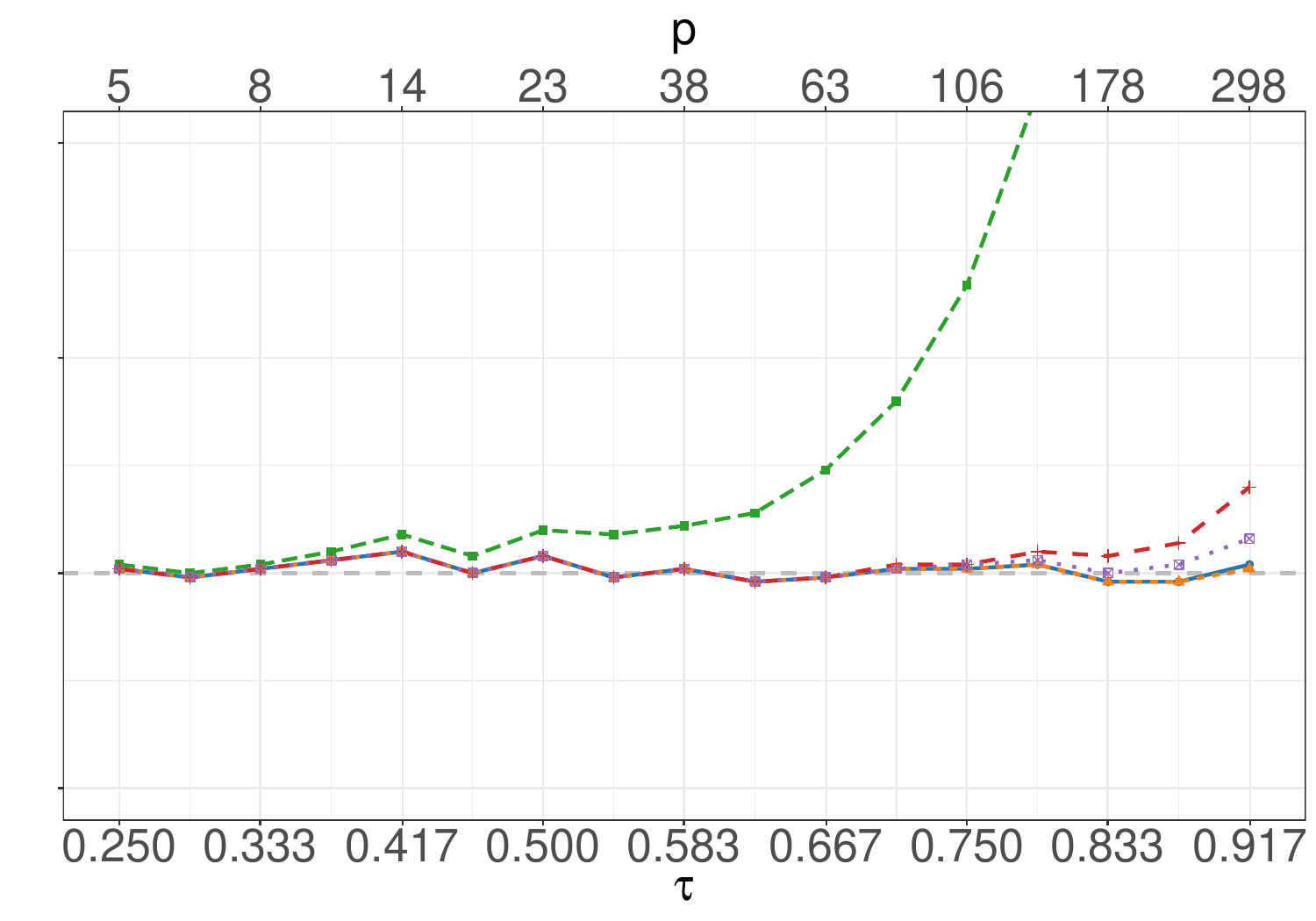}
		\end{minipage}
		\begin{minipage}[b]{.33\linewidth}
			\centering
			\includegraphics[scale=0.2]{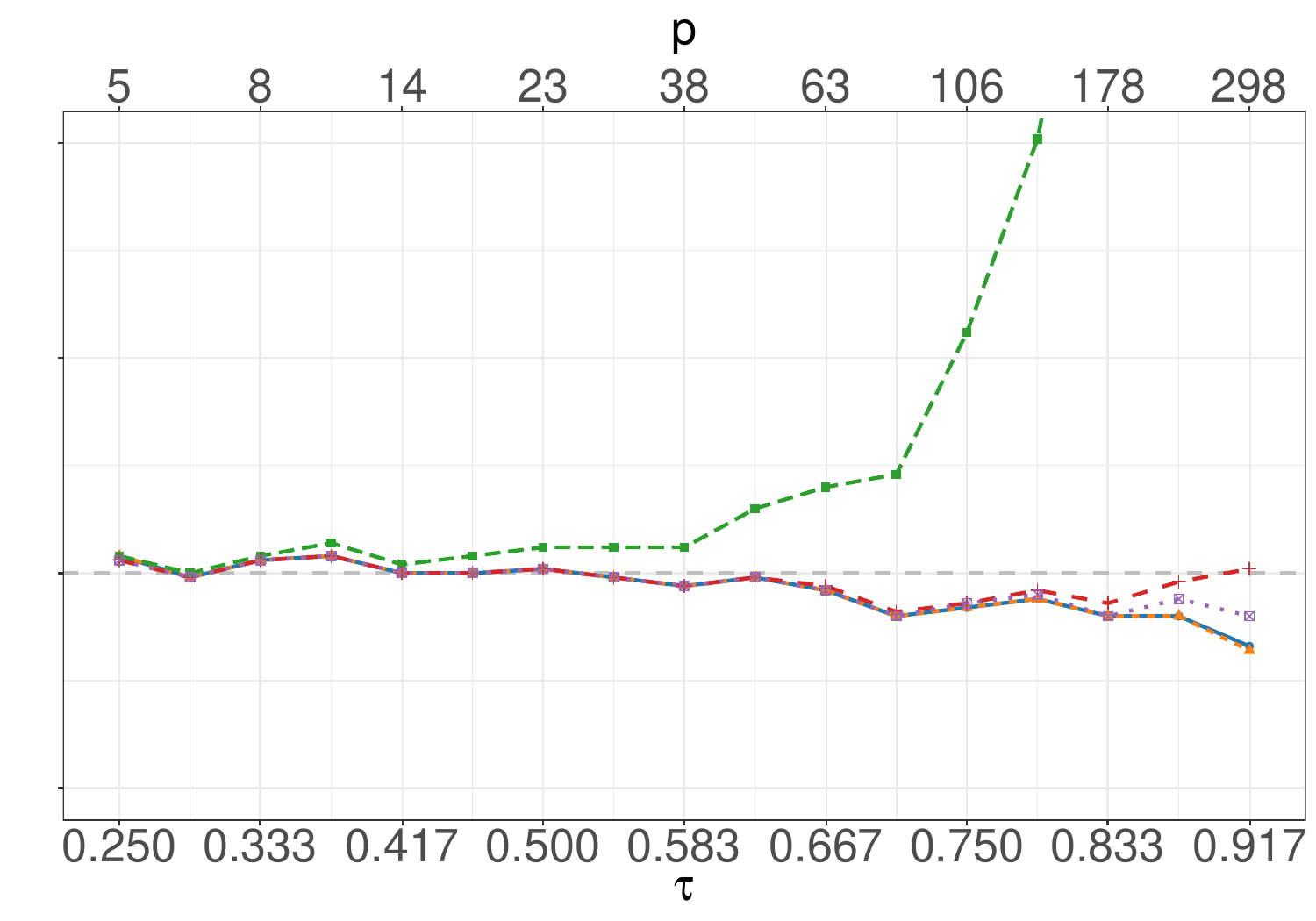}
		\end{minipage}
	}
	\caption{
		Empirical size for the directional test (DT), Behrens-Fisher test (BF), likelihood ratio test (LRT),  and two \citeauthor{skovgaard:2001}'s modifications  (Sko1 and Sko2) for hypothesis (9), at nominal level $\alpha = 0.05$  given by the gray horizontal line. The left, middle and right panels correspond to multivariate $t$, multivariate skew-normal, and multivariate Laplace distributions of the true generating process, respectively, with $n_i = 500$ and $g=2$.
	}
	\label{fig:type I error full robustness 500}
\end{figure}

\section{Application:  Microarray Data of Breast Cancer Patients}

The  dataset studied is taken from  transcript profiles of primary breast tumors to investigate the p53 status. It is based on 251 samples predominantly containing tumor tissue  that underwent sequencing to identify the p53 mutation status through  RNA extraction and subsequent microarray analysis. The RNA was extracted from each sample by using Affymetrix U133 oligonucleotide microarrays, targeting over 30,000 genes \citep{miller2005expression}.   Here we focus on  a  subset of genes corresponding to a signature that is considered important in breast cancer from the research of \cite{bonavita2020antagonistic}, which  stores the appropriately normalized expression values of $p=22$  gene products.  The samples were  divided into three groups based on tumor grade. After excluding two samples due to  the missing  grade variable, we obtain  $g=3$ groups with  $n_1=67$, $n_2=128$ and $n_3=54$ samples for grade I, grade II and grade III, respectively.  Within a high dimensional framework, we have $p/\min(n_i) = 0.407$. Homoscedastic one-way MANOVA leads to the following $p$-values: $1.565 \times 10^{-7}$ for the directional test, $1.966 \times 10^{-12}$ for  the central limit theorem test,  $ 2.627 \times 10^{-8}$ for the likelihood  ratio test, $ 1.915 \times 10^{-7}$ for the Bartlett corrected test, and $1.828 \times 10^{-7}$ and $1.877 \times 10^{-7}$ for the two Skovgaard's  modifications.  All methods clearly indicate to reject the null hypothesis of equal breast tumors gene expression signature in case of homoscedastic groups.  As already seen, the likelihood ratio test results in the smallest $p$-value, whereas its Skovgaard's modified versions and the Bartlett corrected test are closer to the exact directional $p$-value, consistent with our simulations.  To check the reliability  of the various methods in the homoscedastic case, another simulation study on the gene expression data can be conducted. The empirical size  based on $10,000$ replications for the directional, central limit theorem test, likelihood ratio, Bartlett corrected and two Skovgaard's modifications tests  are 0.049,  0.060, 0.087, 0.049 0.048 and  0.047, respectively.  Those of central limit theorem test  and likelihood ratio test are a little inflated. The remaining approaches, including the directional test, appear to be  accurate.

Assuming heteroscedastic groups, for  hypothesis (9) the  directional $p$-value is $2.861 \times 10^{-5}$,  the $p$-value of the likelihood  ratio test is $ 1.928 \times 10^{-7}$, and the $p$-values of the  two Skovgaard's  modifications are  $1.772 \times 10^{-5}$ and $2.079 \times 10^{-5}$.  The same conclusion is hence drawn by all methods, but once again the $p$-value obtained from the likelihood ratio test is relatively smaller than the others.  Such discrepancies are particularly relevant in biostatistical problems, and more generally in medical research, where it is important to look at the specific size of $p$-values to get an idea about the clinical implications and relevance  of the findings for the population of patients of interest.
Even under heterescedasticity a simulation study can be performed. The empirical sizes of the directional test, the likelihood ratio test and the two Skovgaard's modifications  based on $10,000$ replications are  0.051,  0.181, 0.057 and  0.054, respectively. The directional test performs very well, as expected, while the likelihood ratio test seems still not reliable.

\end{document}